%% file: manuscript-numapde-preprint.tex
\pdfoutput=1

\makeatletter
\IfFileExists{./numapde-latex/numapde-packages.sty}{\providecommand*{\input@path}{}\edef\input@path{{./numapde-latex/}\input@path}}{}
\makeatother

\documentclass{numapde-preprint}

\usepackage{numapde-semantic}
\usepackage{numapde-style}
\usepackage{numapde-local}

\IfFileExists{./numapde-bibliography/numapde.bib}{\addbibresource{./numapde-bibliography/numapde.bib}}{\addbibresource{numapde.bib}}

\hypersetup{
	pdftitle={A Discretize-Then-Optimize Approach to PDE-Constrained Shape Optimization},
	pdfauthor={Roland Herzog, Estefanía Loayza-Romero},
	pdfkeywords={discrete shape optimization, mesh quality penalization, Riemannian metric, shape gradient}
}

\title{A Discretize-Then-Optimize Approach to PDE-Constrained Shape Optimization}
\subtitle{}
\shorttitle{Discretize-Then-Optimize Approach to Shape Optimization}

\author{Roland Herzog\thanks{Interdisciplinary Center for Scientific Computing, Heidelberg University, 69120 Heidelberg, Germany (\email{roland.herzog@iwr.uni-heidelberg.de}, \url{https://scoop.iwr.uni-heidelberg.de}, \orcid{0000-0003-2164-6575}).}
\and
Estefanía Loayza-Romero\thanks{Institute for Analysis and Numerics,  University of Münster,   48149 Münster,  Germany (\email{estefania.loayza-romero@uni-muenster.de}, \url{https://www.uni-muenster.de/AMM/num/wirth/people/Loayza/index.html}, \orcid{0000-0001-7919-9259}).}}
\shortauthor{R. Herzog and E. Loayza-Romero}

\dedication{}

\IfFileExists{numapde-uncolorizeHyperrefIfChanges.sty}{\RequirePackage{numapde-uncolorizeHyperrefIfChanges}}{}

\begin{document}
\maketitle

\begin{abstract}
We consider discretized two-dimensional PDE-constrained shape optimization problems, in which shapes are represented by triangular meshes.
Given the connectivity, the space of admissible vertex positions was recently identified to be a smooth manifold, termed the manifold of planar triangular meshes.
The latter can be endowed with a complete Riemannian metric, which allows large mesh deformations without jeopardizing mesh quality; see \cite{HerzogLoayzaRomero:2022:1}.
Nonetheless, the discrete shape optimization problem of finding optimal vertex positions does not, in general, possess a globally optimal solution.
To overcome this ill-possedness, we propose to add a mesh quality penalization term to the objective function. 
This allows us to simultaneously render the shape optimization problem solvable, and keep track of the mesh quality.
We prove the existence of a globally optimal solution for the penalized problem and establish first-order necessary optimality conditions independently of the chosen Riemannian metric.

Because of the independence of the existence results of the choice of the Riemannian metric, we can numerically study the impact of different Riemannian metrics on the steepest descent method.
We compare the Euclidean, elasticity, and a novel complete metric, combined with Euclidean and geodesic retractions to perform the mesh deformation.\end{abstract}


\begin{AMS}
\href{https://mathscinet.ams.org/msc/msc2010.html?t=49Q10}{49Q10}, \href{https://mathscinet.ams.org/msc/msc2010.html?t=49J20}{49J20}, \href{https://mathscinet.ams.org/msc/msc2010.html?t=53Z50}{53Z50}, \href{https://mathscinet.ams.org/msc/msc2010.html?t=35Q93}{35Q93}
\end{AMS}

\input{main.tex}

\printbibliography

\end{document}

%% file: main.tex
\section{Introduction}
\label{section:introduction}

It is well known among practitioners that the numerical solution of shape optimization problems constrained by partial differential equations (PDEs) often exhibits a number of difficulties.
In particular, when the PDE is discretized by a finite element method and the underlying mesh is used to directly represent the shape of the domain to be optimized, one often experiences a degeneracy of the mesh quality as the optimization progresses.
The degeneracy manifests itself in some of the mesh cells thinning in the sense that at least one of its heights approaches zero.

A number of possible solutions to this major obstacle in computational shape optimization have been proposed in the literature.
We do not aim to give a comprehensive overview here but only mention that remeshing \cite{WilkeKokGroenwold:2005:1}, mesh regularization and spatial adaptivity \cite{DoganMorinNochettoVerani:2007:1,MorinNochettoPaulettiVerani:2012:1}, overlapping meshes \cite{DokkenFunkeJohanssonSchmidt:2019:1}, nearly-conformal transformations \cite{IglesiasSturmWechsung:2018:1}, elasticity-based shape gradients \cite{SchulzSiebenbornWelker:2016:1}, and restricted mesh deformations \cite{EtlingHerzogLoayzaWachsmuth:2020:1} have been considered as remedies.
The importance of mesh quality in shape optimization has also been emphasized very recently in \cite{LuftSchulz:2021:2,LuftSchulz:2021:3}, who propose to add certain regularization terms to the objective based on the so-called pre-shape parameterization tracking problem.

In this paper we shed new light on the phenomenon of mesh degeneracy in computational shape optimization.
The discussion is restricted to problems in two dimensions but we believe that much of it extends directly to 3D as well.
We provide evidence that problems in which the nodal positions of the finite element mesh serve as the optimization variables generally possess no solutions in the open set of admissible vertex positions, even when the objective is bounded below.
Moreover, we illustrate by example that as the objective's infimum is approached, the mesh iterates approach the boundary of the set of admissible vertex positions.\footnote{Informally speaking, boundary points are characterized as infeasible meshes in which one or more triangles have a vanishing height, or where a boundary vertex comes into contact with a non-incident boundary edge.}

We therefore must conclude that the class of problems arising from one of the most straightforward and perhaps common approaches for the discretization of PDE-constrained shape optimization problems is generally ill-posed, in the sense that no solution exists.
Consequently, any convergent optimization method, whether of gradient- or Newton-type, will inevitably be led to produce degenerate meshes sooner or later.
We conjecture that this is reason why we often see early stopping and rather loose tolerances for the norm of the gradient in published works.
We are aware of the fact that the ill-posed nature of the problem class under consideration has been noticed previously.
For instance, as \cite{Berggren:2010:1} observed: \enquote{However, in shape optimization, it does not make much sense to optimize the position of each mesh point independently.}
In spite of this observation, we are not aware of a detailed investigation.

\subsubsection*{Our Contributions}
In contrast to the prevailing literature, in this paper we study PDE-constrained shape optimization problems from the discretize-then-optimize perspective.
We characterize the set of admissible vertex positions as an open, connected submanifold of the vector space of all vertex positions, termed the manifold of triangular meshes~$\planarmanifold$.
We demonstrate by example that discretized shape optimization problems generally do not possess solutions in this set.
To overcome this problem we introduce a penalty functional which, briefly speaking, controls the mesh quality.
When added to the shape optimization objective, it renders discrete shape optimization problems well-posed.
Subsequently, such penalized problems can be solved by standard gradient- or Newton-type optimizers.

Gradient algorithms need to be endowed with a (Riemannian) metric. In our context, the Euclidean and an elasticity-based metric are often considered.
A third option, inspired by \cite{HerzogLoayzaRomero:2022:1}, is based on the observation that 
it is generally useful to `precondition' gradient- and quasi-Newton type methods by means of a `base metric' which takes into account some information about the objective.
We therefore devise a Riemannian metric on $\planarmanifold$ which is informed about the penalty part of the objective.
It will be proved in \cref{proposition:NewCompleteMetric} that this metric is complete, and thus, degenerate meshes, which lie on the boundary of $\planarmanifold$, are infinitely far away from any regular mesh in terms of their geodesic distance.
This property is mathematically convenient since an optimization scheme moving along geodesic segments by construction need not explicitly monitor mesh quality and can take arbitrarily large steps.
Unfortunately, the completeness of the metric is also practically difficult to exploit since, unfortunately, the numerical integration of the respective geodesics is prohibitively expensive.
However, we demonstrate that the proposed Riemannian metric is still beneficial to use in gradient methods, even when combined with the inexpensive Euclidean retraction.
Moreover, compared to an elasticity-based metric, the conversion of the shape derivative to the gradient is significantly less expensive since the underlying linear system is governed by a rank-$1$ perturbation of the identity matrix.
It turns out that gradient methods utilizing the complete metric perform well even in the absence of the mesh quality penalty.

The paper is structured as follows.
In \cref{section:manifold_planar_triangular_meshes} we characterize the manifold $\planarmanifold$ of admissible vertex positions.
A new penalty function to control the mesh quality is proposed and analyzed in \cref{section:penalized_shape_optimization}.
We prove that the addition of this penalty function renders a discretized model problem in PDE-constrained shape optimization well posed.
\Cref{section:steepest_descent_method} describes a generic gradient descent method on Riemannian manifolds with Armijo backtracking line search.
We present various numerical experiments with different purposes in \cref{section:numerical_results}.
\Cref{section:conclusions_outlook} offers conclusions and a brief outlook.

\section{Manifold of Planar Triangular Meshes}
\label{section:manifold_planar_triangular_meshes}

As previously mentioned, we study discretized shape optimization problems in which the computational mesh underlying the finite element method serves to represent the sought-after shape, and the coordinates of its vertices serve as optimization variables.
Briefly speaking, a mesh in two space dimensions is a finite collection of non-degenerate triangles such that the intersection of any two triangles is either empty, a common edge, or a common vertex; see for instance \cite[Chapter~3]{QuarteroniValli:1994:1}.
Such meshes are routinely produced by mesh generating software, and the computational verification whether or not a given set of vertices, edges and triangles qualifies as a mesh is relatively straightforward.
In the context of shape optimization, however, a not-so-straightforward question arises: what is the set of possible coordinates the vertices of a mesh be assigned in order for it to remain a mesh with the same connectivity?

In this section we are concerned with precisely this set of admissible vertex positions since it forms the feasible set for discretized shape optimization problems.
As was shown in our previous work \cite{HerzogLoayzaRomero:2022:1}, this set can be characterized as an open, connected submanifold of $\R^{2\times N_V}$, where $N_V$ is the number of vertices.
It is termed the \emph{manifold of planar triangular meshes} and it can be endowed with a complete Riemannian metric.
For convenience, we briefly recall the relevant material from \cite{HerzogLoayzaRomero:2022:1}, to which we refer the reader for a more detailed account.

The conditions which define a mesh can be conveniently described in the language of simplicial complexes.
We work with \emph{abstract} simplicial complexes to describe the connectivity, and we use \emph{geometric} simplicial complexes to formulate conditions on the vertex positions.
An abstract simplicial complex $\Delta$ is a purely combinatorial object, defined over a finite set of vertices, which we denote by $V \coloneqq \{1, 2, \ldots, N_V\}$.
$\Delta$ is a non-empty collection of non-empty subsets of $V$ such that, for all $\sigma \in \Delta$, every non-empty subset of $\sigma$ also belongs to $\Delta$.
Any $\sigma \in \Delta$ is called a face of $\Delta$, and we say that it is of dimension~$k$ (or simply a $k$-face) if its cardinality satisfies $\#\sigma = k+1$.
$0$-faces are vertices, $1$-faces are edges and $2$-faces are triangles.

In order for $\Delta$ to represent the connectivity of a triangular mesh in $\R^2$, we require additional conditions.
We assume that $\Delta$ to be a pure simplicial $2$-complex, \ie, every $\tau \in \Delta$ is contained in some $2$-dimensional~$\sigma \in \Delta$.
Moreover, in order to obtain \emph{rigid} meshes, we require $\Delta$ to be $2$-path connected.
This means that for any two distinct $2$-faces $\sigma, \sigma' \in \Delta$, there exists a finite sequence of $2$-faces, starting in $\sigma_0 = \sigma$ and ending in $\sigma_n = \sigma'$, such that $\sigma_i \cap \sigma_{i+1}$ is a $1$-face for $i = 0,\ldots, n - 1$.
We collect these properties in the following definition.

\begin{definition}
	\label{definition:connectivity_complex}
	Suppose that $\Delta$ is an abstract simplicial $2$-complex with vertex set~$V = \{1, \ldots, N_V\}$.
	We say that $\Delta$ is a \textbf{connectivity complex}, provided that $\Delta$ is pure and $2$-path connected.
\end{definition}

\begin{figure}[htp]
	\centering
	\includegraphics[width=0.25\linewidth]{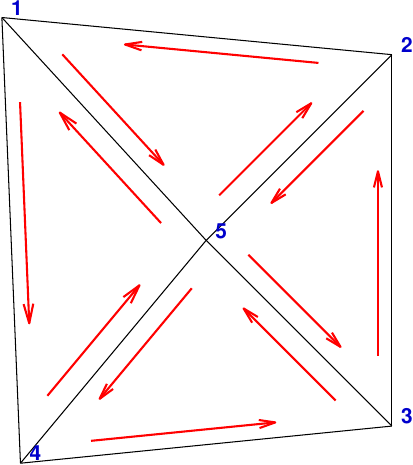}
	\qquad
	\includegraphics[width=0.40\linewidth]{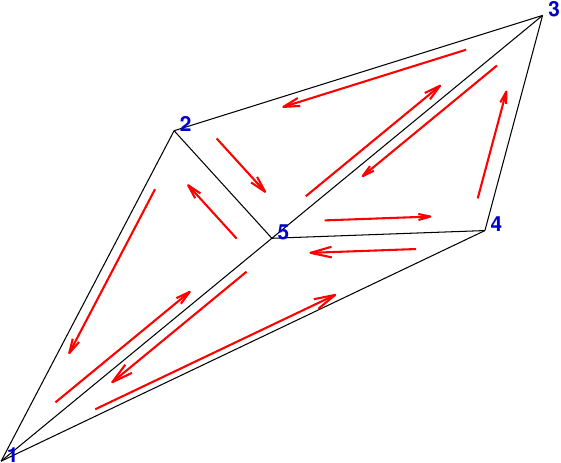}
	\caption{Two admissible oriented meshes with the same consistently oriented connectivity complex~$\Delta$ and different vertex positions~$Q$.}
	\label{fig:points_in_plusmanifold}
\end{figure}

In addition, we recall what it means for a connectivity complex~$\Delta$ to be \textbf{consistently oriented}: the orientations of any two $2$-faces in $\Delta$ sharing a $1$-face induce opposite orientations on that $1$-face; see \cref{fig:points_in_plusmanifold} for an illustration.

An admissible triangular mesh is described using the notion of geometric simplicial complexes, which we briefly recall.
First of all, a simplex~$\sigma$ of dimension~$k$ (a $k$-simplex) in $\R^2$ is the convex hull of $k+1$ affine independent points (called the vertices of $\sigma$), $k = 0,1,2$.
A face of~$\sigma$ is the convex hull of a subset of its vertices.
A geometric simplicial complex $\Sigma$ in $\R^2$ is a non-empty, finite collection of simplices in $\R^2$ which satisfies the following two conditions.
First, every face of any $\sigma \in \Sigma$ is also an element of $\Sigma$.
Second, the non-empty intersection of any two simplices in $\Sigma$ is a face of both.

Let us now connect algebraic and geometric simplicial complexes, \ie, connectivity and geometry.
Suppose that $\Delta$ is a given oriented connectivity complex.
By assigning positions to the vertices of $\Delta$, we can define
\begin{equation}
	\label{eq:collection_of_convex_hulls}
	\Sigma_\Delta(Q)
	\coloneqq
	\setDef[auto]{\conv\{q_{i_0}, \ldots, q_{i_k}\}}{\{i_0, \ldots, i_k\} \in \Delta}
	\subset
	\powerset{\R^2}
	,
\end{equation}
where $\powerset{\R^2}$ denotes the power set of $\R^2$ and $Q \in \R^{2\times N_V}$ records the vertex positions.
Notice that $\Sigma_\Delta(Q)$ is, in general, not an admissible mesh since triangles and edges may intersect in non-admissible ways.
We therefore require that $\Sigma_\Delta(Q)$ forms a geometric simplicial complex.
Moreover, in order to take care of the orientation in $\Delta$, we introduced the signed area of a $2$-face in $\Sigma_\Delta(Q)$ as follows:
\begin{equation}
	\label{eq:signed_area}
	\area{Q}[i_0,i_1,i_2]
	\coloneqq
	\frac{1}{2} \det
	\begin{bmatrix}
		q_{i_1} - q_{i_0}, \; q_{i_2} - q_{i_1}
	\end{bmatrix}
	.
\end{equation}
We point out that triangular mesh generators usually provide orientated meshes for which the signed areas of all $2$-faces are positive.

We are now in the position to define the manifold of planar triangular meshes.
\begin{definition}
	\label{definition:planarmanifold}
	Let $\Delta$ be a connectivity complex as in \cref{definition:connectivity_complex}, which is also oriented.
	We first define the set of admissible oriented meshes with connectivity~$\Delta$, briefly $\plusmanifold$, as the set of points $Q \in \R^{2\times N_V}$ such that the following two conditions hold.
	\begin{enumeratelatin}
		\item
			\label{item:planarmanifold_Sigma_simplicial_complex}
			The collection of convex hulls $\Sigma_\Delta(Q)$, given by \eqref{eq:collection_of_convex_hulls}, is a geometric simplicial complex, whose associated abstract simplicial complex is $\Delta$.
		\item
			\label{item:planarmanifold_signed_area_positive}
			For all $2$-faces of $\Sigma_\Delta(Q)$, the signed area $\area{Q}$ as in \eqref{eq:signed_area} is positive.
	\end{enumeratelatin}
	Moreover, given a reference mesh $\Qref \in \plusmanifold$, the \textbf{manifold of planar triangular meshes} $\planarmanifold$ is the set of points $Q \in \plusmanifold$ such that, in addition:
	\begin{enumeratelatin}[resume]
		\item
			\label{item:planarmanifold_continuous_path}
			There exists a continuous path between $Q$ and $\Qref$, such that all the points along the path satisfy \cref{item:planarmanifold_Sigma_simplicial_complex,item:planarmanifold_signed_area_positive}.
	\end{enumeratelatin}
\end{definition}
We mention that \cref{item:planarmanifold_continuous_path} is a technical condition which ensures that $\planarmanifold$ is path-connected, a property doubtlessly useful for optimization purposes.
The existence of a reference mesh $\Qref$ is not critical.
In fact, in practice one starts with such a reference mesh and can easily extract the underlying oriented connectivity complex~$\Delta$ from it.
Clearly, $\Qref$ belongs to $\planarmanifold$.

It has been proved in \cite[Theoremm~3.11]{HerzogLoayzaRomero:2022:1} that $\planarmanifold$ is a path-connected, open submanifold of $\R^{2 \times N_V}$.
Moreover, a complete Riemannian metric for this manifold has been proposed in the same paper.
Although we will use a different complete metric below more suitable for PDE-constrained shape optimization problems, we briefly recall the original construction.
\cite{Gordon:1973:1} showed that the construction of complete metrics is tied to the definition of a proper, real-valued function, which we refer to as an \emph{augmentation} function.
In our setting, it is the purpose of the augmentation function to penalize impending self-intersections, which can be of interior or exterior nature.
The term which avoids interior self-intersections is based on the heights of a $2$-face $\{i_0,i_1,i_2\}$, which are defined as follows:
\begin{equation}
	\label{eq:triangle_height}
	\height[auto]{Q}{\ell}[i_0,i_1,i_2]
	=
	\frac{2 \, \area{Q}[i_0,i_1,i_2]}{\edgelength{Q}{\ell}[i_0,i_1,i_2]}
	,
	\quad
	\ell = 0, 1, 2
	,
\end{equation}
where
\begin{equation}
	\label{eq:edge_length}
	\edgelength{Q}{\ell}[i_0,i_1,i_2]
	\coloneqq
	\norm{q_{i_{\ell \oplus 1}} - q_{i_{\ell \oplus 2}}}
	,
	\quad
	\ell = 0, 1, 2
\end{equation}
is the Euclidean length of the $\ell$-th edge (the one opposite the $\ell$-th vertex).
Here $\oplus$ denotes addition modulo~$3$.
Moreover, the term which avoids exterior self-intersections is defined in terms of a notion of distance of a vertex~$i_0$ to an edge $\{j_0,j_1\}$, given by
\begin{equation}
	\label{eq:regularized_distance_vertex_edge}
	\regularizedDistance{Q}[i_0][[j_0,j_1]]
	\coloneqq
	\min\setDef[big]{\norm{q_{i_0} - q}_{1,\mu}}{q \in \conv\{q_{j_0}, q_{j_1}\}}
	.
\end{equation}
Here $\norm{\cdot}_{1,\mu}$ denotes a $C^3$-regular estimate (possibly depending on a parameter $\mu > 0$) of the $1$-norm in an orthogonal coordinate frame aligned with the edge $\conv\{q_{j_0}, q_{j_1}\}$.
Moreover, it is required that $0 \le \norm{\cdot}_{1,\mu} \le \norm{\cdot}_1$ holds, so that $\regularizedDistance{Q}$ becomes a non-negative underestimate of the true distance in the $1$-norm; see \cite[Section~4 and Appendix~D]{HerzogLoayzaRomero:2022:1} for details and a concrete example of such a function.

The following augmentation function was proposed in \cite{HerzogLoayzaRomero:2022:1} and proved to be proper:
\begin{definition}
	\label{definition:f_mu}
	Suppose that $\Delta$ and $\Qref$ are as in \cref{definition:planarmanifold}.
	Denote by $V_\partial$ the set of the boundary $0$-faces and by $E_\partial$ the set of boundary $1$-faces.
	Suppose that the $2$-faces in $\Delta$ are numbered from~$1$ to~$N_T$ and that the $k$-th triangle has vertices $i_0^k, i_1^k, i_2^k$.
	Define the augmentation function $\oldpenalty \colon \plusmanifold \to \R$ by
	\begin{equation}
		\label{eq:f_mu}
		\begin{multlined}
			\oldpenalty(Q;\Qref)
			\coloneqq
			\sum_{k=1}^{N_T} \sum_{\ell=0}^2 \frac{\beta_1}{\height[auto]{Q}{\ell}[i^k_0,i^k_1,i^k_2]}
			+
			\sum_{[j_0,j_1]\in E_\partial}\sum_{\substack{i_0 \in V_\partial \\ i_0 \neq j_0,j_1}} \frac{\beta_2}{\regularizedDistance{Q}[i_0][[j_0,j_1]]}
			+
			\frac{\beta_3}{2} \norm{Q - \Qref}_F^2
			.
		\end{multlined}
	\end{equation}
	Here $\norm{\cdot}_F$ denotes the Frobenius norm.
\end{definition}

We close this section by recalling the result which introduces a complete metric on $\planarmanifold$.
\begin{theorem}[\protect{\cite[Theorem~4.12]{HerzogLoayzaRomero:2022:1}}]
	\label{theorem:f_mu_is_proper_and_Riemannian_metric_is_complete}
	Suppose that $\beta_1, \beta_2, \beta_3 > 0$ holds.
	Then the following statements hold.
	\begin{enumeratelatin}
		\item
			\label[statement]{item:f_mu_is_proper}
			The restriction of $\oldpenalty$ defined in \eqref{eq:f_mu} to $\planarmanifold$ is proper.

		\item
			\label[statement]{item:Riemannian_metric_is_complete}
			The manifold $\planarmanifold$, endowed with the Riemannian metric whose components (with respect to the $\vvec$ chart) are given by
			\begin{equation}
				\label{eq:complete_metric_for_planar_triangular_meshes}
				g_{ab}
				=
				\delta_a^b + \dfrac{\partial \oldpenalty}{\partial (\vvec Q)^a} \dfrac{\partial \oldpenalty}{\partial (\vvec Q)^b}
				,
				\quad
				a, b = 1, \ldots, 2 \, N_V
				,
			\end{equation}
			is geodesically complete.
	\end{enumeratelatin}
\end{theorem}
Here $\vvec$ denotes the vectorization operation $\vvec \colon \R^{2 \times N_V} \to \R^{2 \, N_V}$ which stacks $Q \in \R^{2 \times N_V}$ column by column.
Moreover, $\delta_a^b$ denotes the Kronecker delta symbol, representing the Euclidean metric.

\begin{remark}
	\label{remark:rank1perturbationofIdentity}
	Notice that from an algebraic perspective, the matrix of components of the metric given in \eqref{eq:complete_metric_for_planar_triangular_meshes} is a symmetric positive definite rank-$1$ perturbation of the identity matrix.
	This fact simplifies tremendously the solution of linear systems with \eqref{eq:complete_metric_for_planar_triangular_meshes}, which occur, \eg, when converting derivatives to gradients (covectors to vectors).
	For instance, the inverse of the metric coefficients matrix can be explicitly computed using the Sherman-Morrison formula, which results in the following representation of the inverse:
	\begin{equation}
		\label{eq:complete_inverse_metric_for_planar_triangular_meshes}
		g^{ab}
		=
		\delta_a^b - \dfrac{1}{1 + \sum\limits_{c=1}^{2 \, N_V} \paren[auto](){\dfrac{\partial \oldpenalty}{\partial (\vvec Q)^c}}^2} \dfrac{\partial \oldpenalty}{\partial (\vvec Q)^a} \dfrac{\partial \oldpenalty}{\partial (\vvec Q)^b}
		.
	\end{equation}
	Alternatively, an iterative solution of linear systems with $g_{ab}$ can be achieved in a matrix-free way by performing just two iterations of the conjugate gradient method.
\end{remark}

We re-iterate that we will devise a new complete metric for $\planarmanifold$ in \cref{subsection:new_penalty_function} which is more suitable for PDE-constrained shape optimization problems.

\section{Discrete Shape Optimization Problems}
\label{section:penalized_shape_optimization}

\subsection{A Model Problem}
\label{subsection:model_problem}

Throughout, we consider a two-dimensional model problem as in \cite{EtlingHerzogLoayzaWachsmuth:2020:1}.
In continuous form it reads:
\begin{equation}
	\label{eq:example_problem_continuous}
	\text{Minimize} \quad \int_\Omega y \d x
	\quad
	\text{\st}
	\quad
	- \laplace y
	=
	\rhs
	\text{ in } \Omega
	\quad
	\text{\wrt}
	\quad
	\Omega \subset \R^2
	.
\end{equation}
The state~$y$ is subject to Dirichlet boundary conditions $y = 0$ on $\partial \Omega$ and the right-hand side function~$\rhs \colon \R^2 \to \R$ is given.
To discretize it, we represent the unknown domain~$\Omega$ by a mesh with coordinates $Q \in \planarmanifold \subset \R^{2 \times N_V}$ and given oriented connectivity complex, as introduced in \cref{definition:connectivity_complex,definition:planarmanifold}.

We refer to the domain covered by the mesh with vertex coordinates~$Q$ as $\Omega_Q$.
We discretize the PDE in \eqref{eq:example_problem_continuous} by the finite element method.
To this end, let $S^1(\Omega_Q)$ denotes the finite element space of piecewise linear, globally continuous functions, defined over $\Omega_Q$, and let $S^1_0(\Omega_Q)$ denote the subspace of functions with zero Dirichlet boundary conditions.
The discrete version of \eqref{eq:example_problem_continuous} then becomes
\begin{equation}
	\label{eq:discreteProblem}
	\begin{aligned}
		\text{Minimize}
		&
		\quad
		\int_{\Omega_Q} y \d x
		\quad
		\text{\wrt}
		\quad
		Q \in \planarmanifold
		,
		\;
		y \in S^1_0(\Omega_Q)
		\\
		\text{\st}
		&
		\quad
		\int_{\Omega_Q} \nabla y \cdot \nabla v \d x
		=
		\int_{\Omega_Q} \rhs \, v \d x
		\quad
		\text{for all }
		v \in S^1_0(\Omega_Q)
		.
	\end{aligned}
\end{equation}

\subsection{A First Glimpse at the Non-Existence of Solutions}
\label{subsection:first_glimpse_non-existence}

There is a major difference between the continuous and discrete shape optimization problems \eqref{eq:example_problem_continuous} and \eqref{eq:discreteProblem}.
In the former, smooth and bijective re\-pa\-ra\-metri\-za\-tions of the domain $\Omega$ which preserve the boundary do not change the solution of the state equation, nor the value of the objective.
By contrast, the finite element solution of the state equation in the discretized case depends on the positions of \emph{all} vertices, boundary and interior.
Moreover, degenerate meshes usually lead to unrealistically small objective values, whose infimal value is not attained within $\planarmanifold$

\begin{figure}[htp]
	\begin{subfigure}[t]{0.40\textwidth}
		\centering
		\raisebox{10mm}{%
		\includegraphics[width=0.8\textwidth]{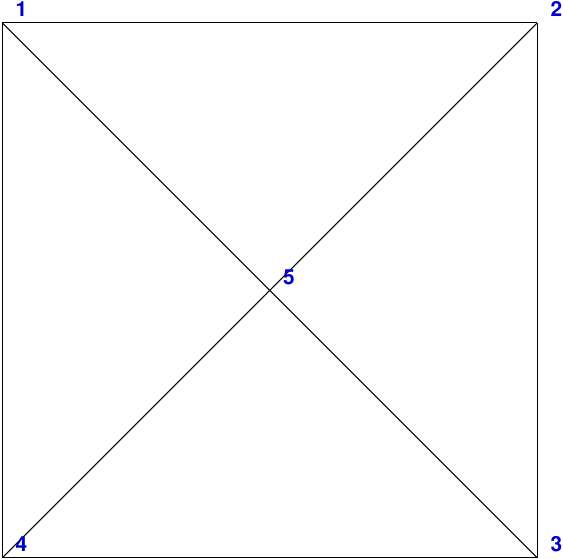}
		}
		\caption{Illustration of the reference mesh $\Qref$.}
		\label{fig:mesh_example_no_solutions}
	\end{subfigure}
	\hfill
	\begin{subfigure}[t]{0.55\textwidth}
		\centering
		\includegraphics[width=\textwidth]{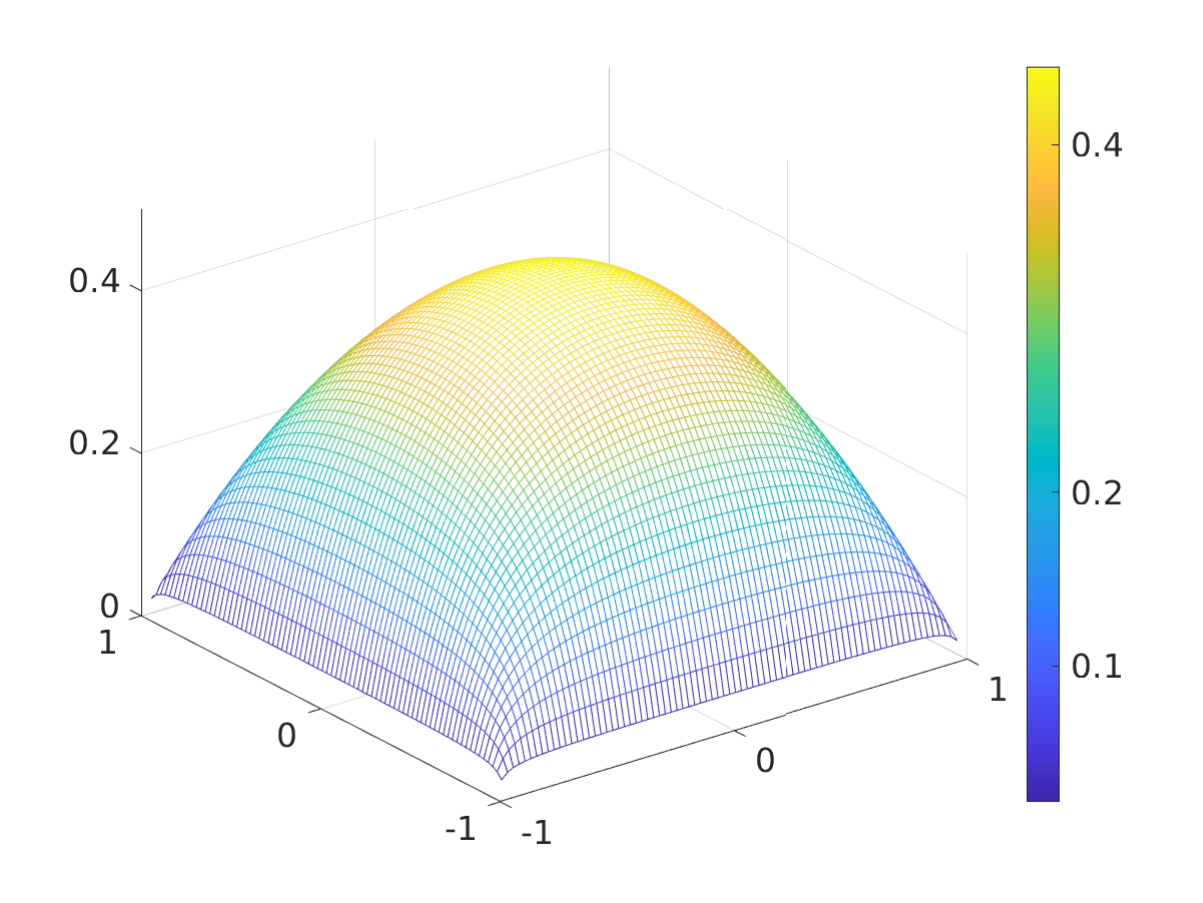}
		\caption{Objective as a function of the nodal position $q_5$.}
		\label{fig:objective_function_example_no_solutions}
	\end{subfigure}
	\caption{Reference mesh and objective function for problem \eqref{eq:example_2_2_problem}.}
	\label{fig:example_no_solutions}
\end{figure}

Let us illustrate this for the simplest possible case.
Consider the reference mesh $\Qref$ covering $[-1,1]^2$ shown in \cref{fig:mesh_example_no_solutions}.
The nodal positions are recorded in $Q = [q_1, q_2, q_3, q_4, q_5] \in \R^{2\times 5}$.
For this experiment, we can even keep the boundary of the shape fixed so that the only remaining unknown is the position of the interior vertex, $q_5$.
It is obvious that $Q \in \planarmanifold$ holds if and only if $q_5 \in (-1,1)^2$.
This leads us to consider the following discrete problem as a particular case of \eqref{eq:discreteProblem},
\begin{equation}
	\label{eq:example_2_2_problem}
	\begin{aligned}
		\text{Minimize}
		&
		\quad
		\int_{\Omega_Q} y \d x
		\quad
		\text{\wrt}
		\quad
		q_5 \in(-1,1)^2
		,
		\;
		y \in S^1_0(\Omega_Q)
		\\
		\text{\st}
		&
		\quad
		\int_{\Omega_Q} \nabla y \cdot \nabla v \d x
		=
		\int_{\Omega_Q} \rhs \, v \d x
		\quad
		\text{for all }
		v \in S^1_0(\Omega_Q)
		.
	\end{aligned}
\end{equation}
For this initial experiment, we fix $r \equiv 1$.
We emphasize that in this scenario, no quadrature error occurs even for the simplest quadrature formula with one evaluation at each cell center.

\Cref{fig:objective_function_example_no_solutions} shows the value of the discrete objective as a function of $q_5$.
It can be observed that the objective takes values arbitrarily close to zero when $q_5$ approaches the boundary of $\Omega_Q$.
To confirm this, consider for instance $q_5 = (0,1-\varepsilon)^\transp$ with a small $\varepsilon > 0$.
It can be easily verified that in this case the linear system representing the PDE in \eqref{eq:discreteProblem} reads $K y = b$ with stiffness matrix $K = \blkdiag(1,1,1,1,4+1/\varepsilon)$ and load vector $b = (0,0,0,0,4/3)^\transp$.
Consequently, the nodal solution vector $y \searrow 0 \in \R^5$ as $\varepsilon \searrow 0$, and thus the value of the objective approaches zero as well.
Similar considerations apply when $q_5$ is anywhere else near the boundary.
Since a location of $q_5$ exactly on the boundary results in a degenerate mesh with $Q \not \in \planarmanifold$, we conclude that the simple problem \eqref{eq:example_2_2_problem} does not have a solution in $\planarmanifold$.
This is in contrast to the continuous problem.
In the continuous setting, due to the fixed boundary, there is no shape to be optimized.
The solution to the state equation on $\Omega = (-1,1)^2$ can be found, \eg, in \cite[Example~1.1.1, p.10]{ElmanSilvesterWathen:2014:1} and the corresponding value of the objective is approximately $0.5622$.

We will later consider in \cref{subsection:unpenalizedAndIdeal} more realistic meshes, a different right-hand side function~$r$ and, of course, impose no constraints which fix the boundary.
However, even this preliminary experiment \eqref{eq:example_2_2_problem} illustrates two fundamental difficulties with discretized shape optimization problems in which the nodal positions serve as the optimization variables.
First, they do not, in general, possess a solution, even if the objective is bounded below.
Second, poor approximations of the state variable can give rise to unreasonably small objective values.
Both observations are related to nearly degenerate finite element meshes.
It is therefore of paramount importance that formulations and solvers for discretized shape optimization problems maintain control over the mesh quality.
Precisely that is the purpose of the penalty function devised in the following subsection.

\subsection{A New Penalty Function and Complete Metric}
\label{subsection:new_penalty_function}

This section proposes a modification of discrete shape optimization problems over the manifold $\planarmanifold$ of planar triangular meshes.
The modification consists in the addition of a penalty function~$\newpenalty$, which renders the resulting problem well-posed in the sense that the existence of a globally optimal solution can be proved.
To the best of the authors' knowledge, we are not aware of existence results for discretized shape optimization problems (in which the vertex positions serve as optimization variables) in the literature.

For the sake of concreteness, we come back to the model problem \eqref{eq:discreteProblem}.
With a penalty term~$\newpenalty$ --- to be specified below --- added, it reads
\begin{equation}
	\label{eq:penalizedProblem}
	\begin{aligned}
		\text{Minimize}
		&
		\quad
		\int_{\Omega_Q} y \d x
		+
		\newpenalty(Q;\Qref)
		\quad
		\text{\wrt}
		\quad
		Q \in \planarmanifold
		,
		\;
		y \in S^1_0(\Omega_Q)
		\\
		\text{\st}
		&
		\quad
		\int_{\Omega_Q} \nabla y \cdot \nabla v \d x
		=
		\int_{\Omega_Q} \rhs \, v \d x
		\quad
		\text{for all }
		v \in S^1_0(\Omega_Q)
		.
	\end{aligned}
\end{equation}

To motivate our choice of penalization, we present a result which guarantees the existence of solutions to an abstract optimization problem in metric spaces.
\begin{proposition}
	\label{proposition:existenceSolutionsProperFunction}
	Suppose that $X$ is a metric space and $\newpenalty \colon X \to \R$ a proper function.
	Moreover, assume that $\newpenalty$ is bounded from below and lower semi-continuous.
	Then the problem
	\begin{equation}
		\label{eq:minimizationProperFunction}
		\text{Minimize}
		\quad
		\newpenalty(x)
		\quad
		\text{\wrt }
		\quad
		x \in X
	\end{equation}
	has at least one globally optimal solution.
\end{proposition}
\begin{proof}
	Let us denote by $\newpenalty_0$ a lower bound for $\newpenalty$.
	We consider a minimizing sequence $\{x^n\}\subset X$, \ie, $\newpenalty(x^n) \searrow \inf \setDef{\newpenalty(x)}{x \in X}$ holds, which implies that the sequence $\{\newpenalty(x^n)\} \subset \R$ is bounded.
	Thus, there exists a constant $K < \infty$ such that $\newpenalty(x^n) \in [\newpenalty_0, K]$ holds for all $n \in \N$.
	Since the interval $[\newpenalty_0, K]$ is compact in $\R$ and thanks to the properness of $\newpenalty$, we know that the set $\newpenalty^{-1}\paren[auto](){[\newpenalty_0, K]}$ is compact in $X$.
	Since $X$ is a metric space, compactness is equivalent to sequential compactness, which in turn implies that we can extract a convergent subsequence from $\{x^n\} \subset \newpenalty^{-1}\paren[auto](){[\newpenalty_0, K]}$, still denoted by $\{x^n\}$.
	Thanks to the lower semi-continuity of $\newpenalty$ and the uniqueness of the limit for $\{\newpenalty(x^n)\}$, we obtain the result.
\end{proof}

Indeed, \cref{proposition:existenceSolutionsProperFunction} is a particular case of a classical result in which one assumes $\newpenalty$ to have at least one non-empty and compact sublevel set.
We formulate a simple corollary tailored to problems of the form \eqref{eq:penalizedProblem}:
\begin{corollary}
	\label{corollary:existenceSolutionpenalizedprob}
	Let $X$ and $\newpenalty$ be as in \cref{proposition:existenceSolutionsProperFunction}.
	Moreover, suppose that $j \colon X \to \R \cup \{\infty\}$ is also bounded from below, lower semi-continuous and not identically equal to $\infty$.
	Then the problem
	\begin{equation}
		\label{eq:minimizationProperFunctionAndRedFunction}
		\text{Minimize}
		\quad
		j(x) + \newpenalty(x)
		\quad
		\text{\wrt }
		\quad
		x \in X
	\end{equation}
	has at least one globally optimal solution.
\end{corollary}

In what follows, $j$ will play the role of the reduced shape functional such as $\int_{\Omega_Q} y \d x$ in \eqref{eq:penalizedProblem}, while $\newpenalty$ denotes the penalty function.
\Cref{corollary:existenceSolutionpenalizedprob} suggests to define the latter so that it is proper on $\planarmanifold$.
Recall, moreover, that the definition of a complete metric on $\planarmanifold$ also relies on a proper function.
Therefore, $\newpenalty$ can serve both purposes at the same time.
We thus require the penalty function~$\newpenalty$ to satisfy the following conditions:
\begin{enumeratelatin}
	\item
		\label[condition]{item:phiProper}
		$\newpenalty \colon \planarmanifold \to \R$ is proper.
	\item
		\label[condition]{item:phiBoundedBelow}
		$\newpenalty \colon \planarmanifold \to \R$ is bounded from below.
	\item
		\label[condition]{item:phiC3}
		$\newpenalty$ is of class~$C^3$.
\end{enumeratelatin}
Moreover, it is desirable for the purpose of shape optimization that $\newpenalty$ satisfies the following additional properties:
\begin{enumeratelatin}[resume]
	\item
		\label[condition]{item:phiInvariantRigidMotions}
		$\newpenalty$ is invariant under rigid body motions (translations and rotations).
	\item
		\label[condition]{item:phiInvariantrefinementsAndScaling}
		$\newpenalty$ is invariant under uniform mesh refinements.
\end{enumeratelatin}

\Cref{item:phiProper,item:phiBoundedBelow} can be used to show the existence of solutions to optimization problems such as \eqref{eq:minimizationProperFunctionAndRedFunction}.
\Cref{item:phiC3} is required for an augmentation function to define a complete metric as in \cref{theorem:f_mu_is_proper_and_Riemannian_metric_is_complete}.
By \cref{item:phiInvariantRigidMotions} we mean the following:
Suppose that $T \colon \R^2 \to \R^2$ is defined by $T(x) = R \, x + b$ with $R \in \text{SO}(2)$ and $b \in \R^2$.
Extend $R$ and $T$ to $\R^{2 \times N_V}$, operating column by column.
Then we ask that $\newpenalty(Q;\Qref) = \newpenalty(T Q;T \Qref)$ holds.
Finally, \cref{item:phiInvariantrefinementsAndScaling} is motivated by applications in PDE-constrained shape optimization.
When every edge of the mesh is bisected and thus every triangle split into four congruent ones, the value of the shape optimization objective~$j$ will remain nearly the same (up to an improvement in the discretization error), and we wish the same to be true for the penalty function~$\newpenalty$.

We mention that the augmentation function~$\oldpenalty$ given in \eqref{eq:f_mu}, which served as the basis of a complete Riemannian metric on $\plusmanifold$ in \cite{HerzogLoayzaRomero:2022:1}, already satisfies \cref{item:phiProper,item:phiBoundedBelow,item:phiC3,item:phiInvariantRigidMotions}.
However it does not satisfy \cref{item:phiInvariantrefinementsAndScaling}.
This motivates us to consider the following function~$\newpenalty$ as an alternative, which satisfies all of the \cref{item:phiProper,item:phiBoundedBelow,item:phiC3,item:phiInvariantRigidMotions,item:phiInvariantrefinementsAndScaling}.
Its construction is based on a well-known triangle quality measure
\begin{equation}
	\label{eq:triangle_quality_measure}
	\frac{(E^0)^2 + (E^1)^2 + (E^2)^2}{4 \sqrt{3}A}
\end{equation}
for the cells in a finite element mesh, first introduced in \cite{BhatiaLawrence:1990:1}; see also \cite[Table~6, Row~4]{Shewchuk:2002:1}.
Here $E^\ell$ ($\ell = 0, 1, 2$) denotes the lengths of the edges, and $A$ refers to the area of a triangular cell.

Our proposal for $\newpenalty$ inherits the terms involving the coefficients $\beta_2$ and $\beta_3$ from $\oldpenalty$ in \eqref{eq:f_mu}.
	However, the $\beta_1$-term, which penalizes small heights and serves to avoid interior self-intersections, is replaced now by a term involving the triangle quality measure.
Since the latter does not take into account the absolute size of a triangle but only its shape, we also add a term which avoids the total area of the mesh going to zero.
Exterior self-intersections, on the other hand, are avoided by a term which agrees with the $\beta_2$-term in \eqref{eq:f_mu}.

\begin{definition}
	\label{definition:phi}
	Suppose that $\Delta$ and $\Qref$ are as in \cref{definition:planarmanifold}.
	Denote by $V_\partial$ the set of the boundary $0$-faces and by $E_\partial$ the set of boundary $1$-faces.
	Their cardinalities are denoted by $\#V_\partial$ and $\#E_\partial$, respectively.
	Suppose that the $2$-faces in $\Delta$ are numbered from~$1$ to~$N_T$ and that the $k$-th triangle has vertices $i_0^k, i_1^k, i_2^k$.
	For parameters $\alpha_j \ge 0$, for $j=1,2,3,4$, define $\newpenalty \colon \planarmanifold \to \R$ as
	\begin{equation}
		\label{eq:penalization}
		\begin{aligned}
			\newpenalty(Q;\Qref)
			&
			\coloneqq
			\sum_{k=1}^{N_T} \frac{1}{N_T} \frac{\alpha_1}{\psi_Q(i_0^k,i_1^k,i_2^k)}
			+
			\frac{\alpha_2}{\sum_{k=1}^{N_T}\area[auto]{Q}[i_0^k,i_1^k,i_2^k] }
			\\
			&
			\quad
			+
			\sum_{[j_0,j_1] \in E_\partial}
			\sum_{\substack{i_0 \in V_\partial \\ i_0 \neq j_0,j_1}} \frac{1}{\#E_\partial \#V_\partial} \frac{\alpha_3}{\regularizedDistance{Q}[i_0][[j_0,j_1]]}
			+
			\frac{\alpha_4}{2} \norm{Q-\Qref}_F^2
		\end{aligned}
	\end{equation}
	with
	\begin{equation}
		\label{eq:Shewchuk64_qualitymeasure}
		\frac{1}{\psi_Q(i_0,i_1,i_2)}
		\coloneqq
		\frac{\paren[big](){\edgelength{Q}{0}[i_0,i_1,i_2]}^2+\paren[big](){\edgelength{Q}{1}[i_0,i_1,i_2]}^2+\paren[big](){\edgelength{Q}{2}[i_0,i_1,i_2]}^2}{4 \sqrt{3} \, \area{Q}[i_0,i_1,i_2]}
		.
	\end{equation}
\end{definition}
Recall that the regularized distance from a vertex to an edge~$\regularizedDistance{Q}$ was defined in \eqref{eq:regularized_distance_vertex_edge}, the edge lengths $\edgelength{Q}{\ell}$ are given in \eqref{eq:edge_length}, the signed area $\area{Q}$ can be found in \eqref{eq:signed_area} and $\norm{\cdot}_F$ is the Frobenius norm.

It is not difficult to see that $\newpenalty$ satisfies \cref{item:phiBoundedBelow,item:phiC3,item:phiInvariantRigidMotions,item:phiInvariantrefinementsAndScaling}.
More details are provided in the following, where we also compare $\newpenalty$ to the proper function~$\oldpenalty$ from \eqref{eq:f_mu} devised in \cite{HerzogLoayzaRomero:2022:1}.
A proof of the properness of $\newpenalty$, \ie, of \cref{item:phiProper}, will follow in \cref{theorem:phiProper}.
\begin{itemize}
	\item
		For any triangle, the function $\nicefrac{1}{\psi_Q}$ is bounded below by $1$, and this bound is attained if and only if the triangle is equilateral.
		This is due to the so-called Weitzenböck inequality; see \cite{AlsinaNelsen:2008:1}.
		Obviously, using this function as a penalty encourages minimizing meshes whose $2$-faces are \emph{as equilateral as possible}.
		Since the terms $\regularizedDistance{Q}$, $\area{Q}$ and $\norm{\cdot}_F$ are always non-negative, $\newpenalty$ is bounded from below, \ie, it satisfies \cref{item:phiBoundedBelow}.

	\item
		$\newpenalty$ is also clearly smooth on $\planarmanifold$ and thus satisfies \cref{item:phiC3}.

	\item
		The invariance of $\newpenalty$ under rigid body motions, \ie, \cref{item:phiInvariantRigidMotions}, follows directly from its definition.

	\item
		The scaling by $N_T,\#E_\partial, \#V_\partial$ is chosen so as to achieve invariance of $\newpenalty$ under uniform mesh refinement, see \cref{item:phiInvariantrefinementsAndScaling} above.

	\item
		The term associated with $\alpha_2$ penalizes small total areas of the entire mesh.
		Even in the continuous case, the inclusion of such a term into the objective makes sense in order to avoid domains shrinking to a point becoming optimal.
\end{itemize}

In order to prove that $\newpenalty$ is a proper function on $\planarmanifold$, the following result is essential.
It shows that on any non-empty sublevel set of $\newpenalty$, the edge lengths $\edgelength{Q}{\ell}$ and the reciprocals of the heights $\nicefrac{1}{\height{Q}{\ell}}$ are uniformly bounded, independently of the vertex positions $Q \in \planarmanifold$.
\begin{proposition}
	\label{proposition:edgelengths_heights_bounded_sublevelset}
	Suppose that $\Delta$ and $\Qref$ are as in \cref{definition:planarmanifold}.
	Consider $\newpenalty$ defined in \eqref{eq:penalization} with $\alpha_j > 0$, $j = 1,\ldots,4$.
	Let $\cN_b$ be a non-empty sublevel set of $\newpenalty$, \ie,
	\begin{equation}
		\label{eq:sublevel_set}
		\cN_b
		\coloneqq
		\setDef[auto]{Q \in \planarmanifold}{\newpenalty(Q;\Qref) \le b}
		=
		\newpenalty(\cdot;\Qref)^{-1}((-\infty,b])
		.
	\end{equation}
	Then there exist constants $c, C, D > 0$ such that the edge lengths and heights satisfy
	\begin{align}
		\label{eq:bounded_edge_lengths}
		c
		\le
		\edgelength{Q}{\ell}(i_0^k,i_1^k,i_2^k)
		&
		\le
		C
		,
		\\
		\label{eq:bounded_reciprocal_heights}
		\frac{1}{\height{Q}{\ell}(i_0^k,i_1^k,i_2^k)}
		&
		\le
		D
	\end{align}
	for all $Q \in \cN_b$, all $k = 1,\ldots,N_T$ and all $\ell = 0,1,2$.
	The constants $c, C, D$ are independent $k$ and $\ell$.
\end{proposition}

The proof of \cref{proposition:edgelengths_heights_bounded_sublevelset} builds on the fact that $\Delta$ is a connectivity complex in the sense of \cref{definition:connectivity_complex}, and in particular it uses the $2$-path connectedness of $\Delta$.
We encourage the reader to check the proof in \cref{appendix:proofs_varphi} for details.

Now, we are ready to prove the properness of $\newpenalty$ from \eqref{eq:penalization}, by relating it to the function~$\oldpenalty$ given in \eqref{eq:f_mu}, for which properness has already been proved; see \cref{theorem:f_mu_is_proper_and_Riemannian_metric_is_complete}.

\begin{theorem}
	\label{theorem:phiProper}
	Suppose that $\Delta$ and $\Qref$ are as in \cref{definition:planarmanifold}.
	Consider the functions~$\oldpenalty$ from \eqref{eq:f_mu} and $\newpenalty$ from \eqref{eq:penalization} with all coefficients $\beta_j$ and $\alpha_j$ strictly positive.
	Then for any sublevel set $\cN_b$ of $\newpenalty$ as in \eqref{eq:sublevel_set}, there exists a constant $B > 0$ such that $\cN_b \subset \oldpenalty(\cdot;\Qref)^{-1}\paren[auto](){[0,B]}$.
	Therefore, $\newpenalty$ is proper.
\end{theorem}
\begin{proof}
	Let us consider vertex positions $Q \in \cN_b$.
	From \cref{proposition:edgelengths_heights_bounded_sublevelset} and the definition of $\newpenalty$, we obtain the following estimates:
	\begin{align*}
		\sum_{k=1}^{N_T} \sum_{\ell = 0}^2 \frac{1}{\height{Q^n}{\ell}(i_0^k,i_1^k,i_2^l)}
		&
		\le
		\frac{3N_TD}{\alpha_1}
		,
		\\
		\sum_{[j_0,j_1] \in E_\partial}
		\sum_{\substack{i_0 \in V_\partial \\ i_0 \neq j_0,j_1}} \frac{1}{\regularizedDistance{Q^n}[i_0][[j_0,j_1]]}
		&
		\le
		\frac{b\#E_\partial \#V_\partial}{\alpha_3}
		,
		\\
		\frac{1}{2} \norm{Q^n-\Qref}_F^2
		&
		\le
		\frac{b}{\alpha_4}
		.
	\end{align*}
	Recalling the definition of $\oldpenalty$ from \eqref{eq:f_mu}, we also have
	\begin{equation*}
		\oldpenalty(Q;\Qref)
		\le
		3 N_T D \, \frac{\beta_1}{\alpha_1}
		+ b \, \#E_\partial \, \#V_\partial \, \frac{\beta_2}{\alpha_3}
		+ b \, \frac{\beta_3}{\alpha_4}
		\eqqcolon
		B
		.
	\end{equation*}
	Since $Q \in \cN_b \subset \planarmanifold$ holds, we also know $\oldpenalty(Q;\Qref) \ge 0$, which in turn implies $Q \in \oldpenalty(\cdot;\Qref)^{-1}\paren[auto](){[0,B]}$.

	To show the properness of $\newpenalty$, consider any compact subset~$K$ of $\R$.
	We need to verify that $\newpenalty(\cdot;\Qref)^{-1}(K)$ is compact in $\planarmanifold$.
	In case $\newpenalty(\cdot;\Qref)^{-1}(K)$ is empty, nothing is to be shown.
	Otherwise, we can find an interval $(-\infty,b]$ such that $\newpenalty(\cdot;\Qref)^{-1}(K) \subset \cN_b = \newpenalty(\cdot;\Qref)^{-1}((-\infty,b])$ holds.
	In the rest of the proof we are going to show that $\cN_b$ is compact.
	Since $\newpenalty$ is continuous on $\planarmanifold$, this then implies that $\newpenalty(\cdot;\Qref)^{-1}(K)$ is a closed subset of a compact set, and thus also compact.

	Let us now prove that $\cN_b$ is compact in $\planarmanifold$.
	Since the latter is a metric space (endowed here with the Euclidean metric of $\R^{2 \times N_V}$), compactness is equivalent to sequential compactness.
	Hence, we consider a sequence $\{Q^n\} \subset \cN_b$.
	Thanks to the first part of the proof, $Q^n$ also belongs to $\oldpenalty(\cdot;\Qref)^{-1}([0,B])$.
	Owing to the properness of $\oldpenalty$ (\cref{theorem:f_mu_is_proper_and_Riemannian_metric_is_complete}), we know that $\oldpenalty(\cdot;\Qref)^{-1}\paren[auto](){[0,B]}$ is sequentially compact.
	Therefore, we can extract from $\{Q^n\}$ a subsequence, denoted again by $\{Q^n\}$, which converges to some $Q^*$ in $\planarmanifold$.
	Thanks to the continuity of $\newpenalty$ on $\planarmanifold$, $Q^* \in \cN_b$ holds, which shows the desired sequential compactness of $\cN_b$.
\end{proof}

\begin{remark}
	\label{remark:addCutOffAugmentation}
	Similar to \cite[Remark~4.13]{HerzogLoayzaRomero:2022:1}, we can add $C^3$ cut-off functions to various terms in $\newpenalty$ while maintaining the properness of the function.
	For instance, \cref{theorem:phiProper} remains true when the function $\newpenalty$ given in \eqref{eq:penalization} is replaced by
	\begin{multline}
		\label{eq:penalization_with_cut-off}
			\newpenalty(Q;\Qref)
			\coloneqq
			\sum_{k=1}^{N_T} \frac{1}{N_T} \frac{\alpha_1}{\psi_Q(i_0^k,i_1^k,i_2^k)}
			+
			\frac{\alpha_2}{\sum_{k=1}^{N_T}\area{Q}(i_0^k,i_1^k,i_2^k)}
			\\
			+
			\sum_{[j_0,j_1]\in E_\partial}
			\sum_{\substack{i_0 \in V_\partial \\ i_0 \neq j_0,j_1}} \frac{\alpha_3}{\#E_\partial \#V_\partial} \chi \paren[auto](){\frac{1}{\regularizedDistance{Q}[i_0][[j_0,j_1]]}}
			+
			\frac{\alpha_4}{2} \norm{Q-\Qref}_F^2
			.
	\end{multline}
	Here $\chi$ is a cut-off function of class~$C^3$ which satisfies $\chi(s) = 0$ on some interval $[0,\underline{s}]$ and $\chi = s$ for $s \ge \overline{s}$.
	Similar cut-off functions could be added to any of the three remaining terms in \eqref{eq:penalization_with_cut-off} as well.
\end{remark}

The properness of $\newpenalty$ guarantees the existence of solutions to the penalized discrete shape optimization  model problem \eqref{eq:penalizedProblem}.
The proof of this result is presented in \cref{proposition:existenceSolutions} under the customary assumption of a hold-all domain.
We define the latter by requiring that all nodal positions belong to a certain box, \ie,
\begin{equation}
	\label{eq:hold-all}
	D
	\coloneqq
	\setDef{Q = [q_1,\ldots, q_{N_V}]\in \planarmanifold}{q_i \in [\underline{a},\overline{a}]\times [\underline{b},\overline{b}] \text{ for all } i = 1,\ldots, N_V}
\end{equation}
for some constants $\underline{a} < \overline{a}$ and $\underline{b} < \overline{b}$.
Notice that this implies that the mesh $\Omega_Q$ itself lies inside $[\underline{a},\overline{a}] \times [\underline{b},\overline{b}]$.
\begin{proposition}
	\label{proposition:existenceSolutions}
	Let $\newpenalty$ be as in \eqref{eq:penalization} or \eqref{eq:penalization_with_cut-off} with $\alpha_j > 0$, $j = 1,\ldots,4$.
	Suppose, moreover, that $\Qref$ belongs to the hold-all~$D$ as in \eqref{eq:hold-all}.
	Denote by $I_D(Q)$ the characteristic function of $D$.
	Finally, suppose that $\rhs$ belong to $L^\infty([\underline{a},\overline{a}]\times [\underline{b},\overline{b}])$.
	Then the problem
	\begin{equation}
		\label{eq:penalizedProblem_withCharacteristicFunction}
		\begin{aligned}
			\text{Minimize}
			&
			\quad
			\int_{\Omega_Q} y \d x
			+
			I_D(Q)
			+
			\newpenalty(Q;\Qref)
			\quad
			\text{\wrt}
			\quad
			Q \in \planarmanifold
			,
			\;
			y \in S^1_0(\Omega_Q)
			\\
			\text{\st}
			&
			\quad
			\int_{\Omega_Q} \nabla y \cdot \nabla v \d x
			=
			\int_{\Omega_Q} \rhs \, v \d x
			\quad
			\text{for all }
			v \in S^1_0(\Omega_Q)
		\end{aligned}
	\end{equation}
	has at least one globally optimal solution in $\planarmanifold$.
\end{proposition}
\begin{proof}
	By virtue of \cref{corollary:existenceSolutionpenalizedprob,theorem:phiProper} it is enough to show that the function $\int_{\Omega_Q} y \d x + I_D(Q)$ is bounded from below, lower semi-continuous and not identically equal to $\infty$.
	First we note that $I_D$ is lower semi-continuous since $D$ is closed in $\R^{2 \times N_V}$ and thus closed in $\planarmanifold$.
	On the other hand, the continuity of $\int_{\Omega_Q} y \d x$ follows from the continuity of the mass matrix and the inverse of the stiffness matrix associated with the weak formulation of the partial differential equation, as a function of the vertex coordinates~$Q \in \planarmanifold$.
	Notice, moreover, that $j$ is everywhere finite on $\planarmanifold$ and $I_D$ is not identically equal to $\infty$ since $\Qref \in D$.

	Thanks to the definition of the characteristic function, it remains to be proved that $\int_{\Omega_Q} y \d x$ is bounded from below on $D$.
	Using $L^2(\Omega_Q) \subset L^1(\Omega_Q)$ and Poincaré's inequality, one can obtain the following estimate:
	\begin{align*}
		\int_{\Omega_Q} y \d x
		\ge
		&
		- \norm{y}_{L^1(\Omega_Q)}
		\ge
		- \abs{\Omega_Q}^{\nicefrac{1}{2}} \norm{y}_{L^2(\Omega_Q)}
		\ge
		- \abs{\Omega_Q}^{\nicefrac{1}{2}} \diam(\Omega_Q) \norm{\nabla y}_{L^2(\Omega_Q)}
	\end{align*}
	where $\abs{\Omega_Q}$ stands for the volume of $\Omega_Q$ and $\diam(\Omega_Q)$ is the diameter of $\Omega_Q$.
	From the weak formulation of the state equation and under similar arguments as before, it is easy to see that
	\begin{equation*}
		\norm{\nabla y}_{L^2(\Omega_Q)}
		\le
		\norm{r}_{L^\infty(B)}
		\,
		\abs{\Omega_Q}^{\nicefrac{1}{2}}
		\,
		\diam(\Omega_Q)
	\end{equation*}
	where we abbreviate $B \coloneqq [\underline{a},\overline{a}]\times[\underline{b},\overline{b}]$.
	Altogether this implies that
	\begin{equation*}
		\int_{\Omega_Q} y \d x
		\ge
		- \norm{r}_{L^\infty(B)}
		\,
		\abs{\Omega_Q}
		\,
		\diam(\Omega_Q)^2
		.
	\end{equation*}
	Moreover, it is easy to see that $Q\in D$ implies $\Omega_Q \subset B$, thus giving
	\begin{equation*}
		\int_{\Omega_Q} y \d x
		\ge
		- \norm{r}_{L^\infty(B)}
		\,
		\abs{B}
		\,
		\diam(B)^2
		,
	\end{equation*}
	which concludes the proof.
\end{proof}

\begin{remark}
	Depending on the specific form of the objective~$j$, it may be possible to obtain an existence result even with one or several of the coefficients $\alpha_j$ in \eqref{eq:penalization} equal to zero.
	For instance, suppose that the $j \colon \planarmanifold \to \R$ is such that there exists $A_0 > 0$ and $\varepsilon > 0$ for which $\area{Q} < A_0$ implies $j(Q) \ge j^* + \varepsilon$, where $j^*$ is the infimum of $j$ on $\planarmanifold$.
	Then, the second term in \eqref{eq:penalization} can be omitted, \ie, $\alpha_2$ can be chosen equal to zero.

	Moreover, if the objective function $j(Q)$ is  bounded below on $\planarmanifold$, such as a quadratic tracking-type or compliance-type objective, the existence of solutions follows from \cref{corollary:existenceSolutionpenalizedprob} and there is no need to impose a hold-all domain.
\end{remark}

\begin{figure}[htp]
	\begin{center}
		\begin{subfigure}[t]{0.45\textwidth}
			\centering
			\includegraphics[width=\textwidth]{Figures/example_2_2.pdf}
			\caption{Objective function without penalization as a function of the nodal position $q_5$.}
		\end{subfigure}
		\hfill
		\begin{subfigure}[t]{0.45\textwidth}
			\centering
			\includegraphics[width=\textwidth]{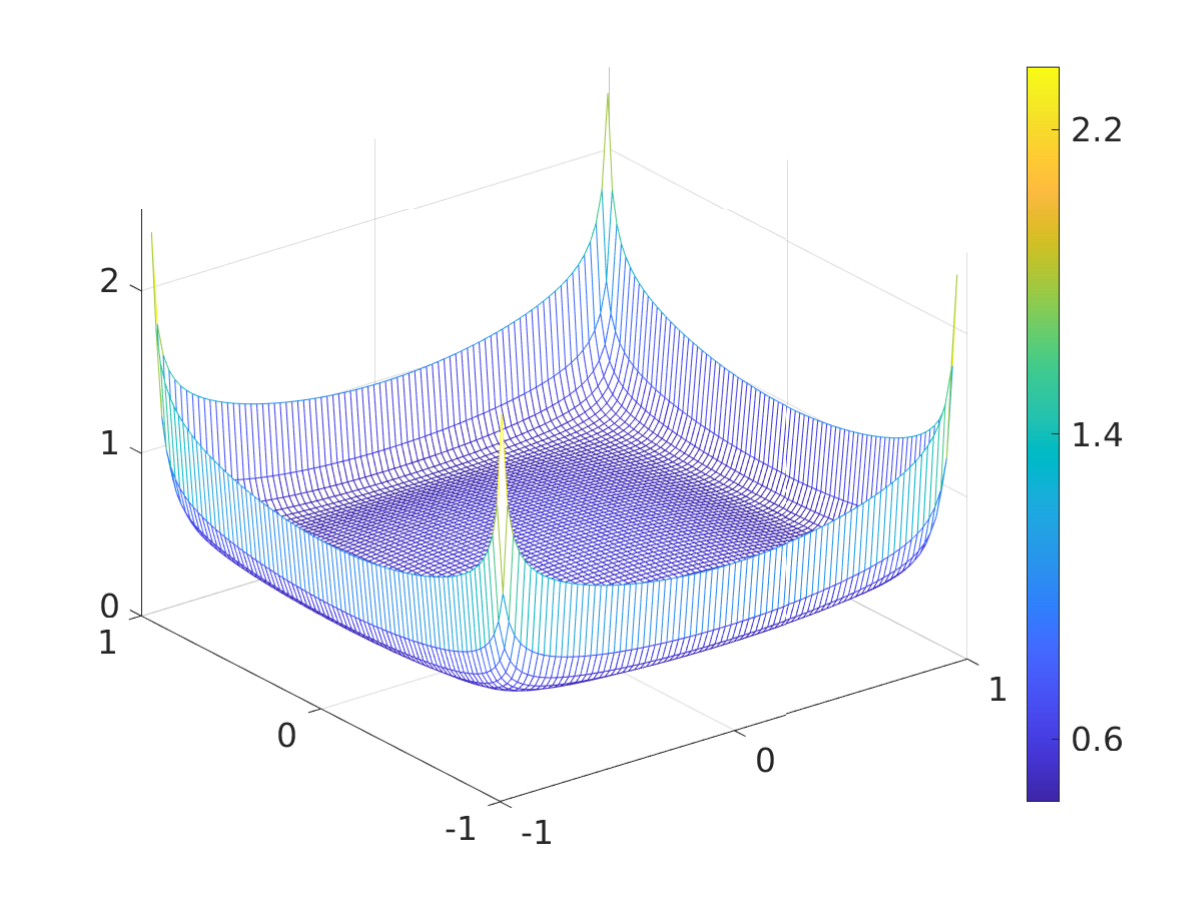}
			\caption{Objective function with penalization ($\alpha_1 = 0.1$, $\alpha_2 = 0.01$, $\alpha_3 = 0$ and $\alpha_4 = 0.01$) as a function of the nodal position $q_5$.}
			\label{fig:example_2_2_penalized}
		\end{subfigure}
		\caption{Transformation of the unpenalized problem \eqref{eq:example_2_2_problem} lacking a solution (left) into a penalized problem which has a solution (right), using the penalty function~$\newpenalty$.}
	\end{center}
\end{figure}

We now revisit example \eqref{eq:example_2_2_problem}, which served as a counterexample to the existence of solution for discrete shape optimization problems in \cref{subsection:first_glimpse_non-existence}.
With the penalty~$\newpenalty$ added, the existence of a solution now follows from \cref{proposition:existenceSolutions}.
The definition of a hold-all is actually not required since the boundary is fixed.
For the same reason, the boundary self-intersection term in $\newpenalty$ is not necessary, \ie, $\alpha_3$ can be set to zero.
To confirm the existence of a solution for this simple example, \cref{fig:example_2_2_penalized} shows a comparison of the objectives with and without penalization, the former with parameters $\alpha_1 = 0.1$, $\alpha_2 = 0.01$, $\alpha_3 = 0$ and $\alpha_4 = 0.01$.
As in \cref{fig:mesh_example_no_solutions}, the right-hand side is chosen as $\rhs \equiv 1$ in \eqref{fig:example_no_solutions}.

We end this section by introducing the first-order necessary optimality conditions for the penalized model problem \eqref{eq:penalizedProblem}.
A stationary point $Q^* \in \planarmanifold$ of $j + \newpenalty \colon \planarmanifold\rightarrow \R$ is characterized by vanishing directional derivatives, \ie,
\begin{equation}
	\label{eq:necessary_optimality_conditions}
	\d_Q \paren[auto][]{j + \newpenalty}[Q^*][V]
	=
	0
	\quad
	\text{for all }
	V \in \tangent{Q^*}[\planarmanifold],
\end{equation}
where $\tangent{Q^*}[\planarmanifold]$ denotes the tangent space to $\planarmanifold$ at $Q^*$, which here agrees with $\R^{2 \times N_V}$.
Using any Riemannian metric $\riemannian{\cdot}{\cdot}[Q^*]$ on $\planarmanifold$, we can define the gradient via
\begin{equation}
	\label{eq:definition_riemannian_gradient}
	\riemannian[auto]{\grad \paren[auto][]{j + \newpenalty}}{V}[Q^*]
	=
	\d_Q \paren[auto][]{j + \newpenalty}[Q^*][V]
	\quad
	\text{ for all }
	V \in \tangent{Q^*}[\planarmanifold]
\end{equation}
and, equivalently to \eqref{eq:necessary_optimality_conditions}, write
\begin{equation}
	\label{eq:necessary_optimality_conditions_using_gradient}
	\grad \paren[auto][]{j + \newpenalty}(Q^*)
	=
	0
	.
\end{equation}

For the model problem \eqref{eq:penalizedProblem} at hand, the derivative $\d_Q j$ can be characterized using the adjoint equation.
This leads to the following formulation of the first-order necessary optimality conditions.
\begin{proposition}
	\label{proposition:optimality_system}
	Let $Q$ be a locally optimal solution to \eqref{eq:penalizedProblem} with associated state $y$.
	Then, there exists a unique adjoint state $p \in S^1_0(\Omega)$ such that the following system of equations is satisfied:
	\begin{subequations}
		\label{eq:optimality_system}
		\begin{align}
			\label{eq:stateEquation}
			\quad
			&
			\int_{\Omega_Q} \nabla y \cdot \nabla e_a \d x
			-
			\int_{\Omega_Q} \rhs \, e_a \d x
			=
			0
			&
			\text{for all } a = 1, \ldots, N_V
			,
			\\
			\label{eq:adjointEquation}
			\quad
			&
			\int_{\Omega_Q}\nabla p \cdot \nabla e_b \d x
			+
			\int_{\Omega_Q} e_b \d x
			=
			0
			&
			\text{for all } b = 1, \ldots, N_V
			,
			\\
			\nonumber
			\quad
			&
			\int_{\Omega_Q} y \div \bV_i \d x
			+
			\mrep{\int_{\Omega_Q} (\nabla y)^\transp \paren[auto][]{\div \bV_i - D\bV_i - D\bV_i^\transp} \nabla p \d x}{}
			\\
			\label{eq:designEquation}
			&
			\quad
			- \int_{\Omega_Q} \div(\rhs \bV_i) \, p \d x
			+
			\frac{\partial \newpenalty (Q)}{\partial (\vvec Q)^i}
			=
			0
			&
			\text{for all } i = 1,\ldots, 2 \, N_V
			.
		\end{align}
	\end{subequations}
	Here $\{e_a\}_{a = 1}^{N_V}$ is the standard nodal finite element basis of $S^1_0(\Omega_Q)$.
	The vector fields $\bV_i$ are defined as follows:
	\begin{equation*}
		\begin{cases}
			\bV_i = (e_{(i+1)/2},0)^\transp & \text{if $i$ is odd},
			\\
			\bV_i = (0,e_{i/2})^\transp & \text{if $i$ is even}.
		\end{cases}
	\end{equation*}
\end{proposition}
We refer to the individual equations in \eqref{eq:optimality_system} as state, adjoint, and design equations, respectively.
The proof of this proposition can be done via standard techniques and is omitted.
We refer the reader to \cite[Theorem~2, p.105]{Pironneau:1984:1} or \cite[Section~4, p.192]{SouliZolesio:1993:1}.

\section{Steepest Descent Method on \texorpdfstring{$\planarmanifold$}{the Manifold of Planar Triangular Meshes}}
\label{section:steepest_descent_method}

In this section we briefly describe a general steepest descent method for the solution of the model problem \eqref{eq:penalizedProblem} on $\planarmanifold$.
The description of the method is kept generic since we wish to conduct numerical experiments for various choices of the Riemannian metric and the retraction later on in \cref{section:numerical_results}.
Clearly, higher-order optimization methods such as quasi-Newton or Newton methods are known to be advantageous with respect to their local convergence properties.
However, a quasi-Newton method would require an implementation of the parallel transport or, more generally, a vector transport associated with the chosen retraction.
By contrast, a Newton method would require the evaluation of the second-order covariant derivative of the penalized objective,
Both of these topics are outside the scope of the present paper.

Let us recall that a Riemannian metric defines a notion of covariant derivative and thus a notion of acceleration along curves.
Consequently, every Riemannian metric defines a notion of geodesics, which can be thought of as acceleration-free curves.
Retractions can be thought of as generalizations of geodesic curves, which agree with respect to their initial points and initial velocities.
A formal definition can be found, for instance, in~\cite[Definition~4.1.1]{AbsilMahonySepulchre:2008:1} and \cite[Definition~3.47]{Boumal:2023:1}.

\begin{algorithm2e}[ht]
	\SetAlgoLined
	\KwData{reference mesh $\Qref \in \plusmanifold \subset \R^{2 \times N_V}$ with oriented connectivity complex $\Delta$}
	\KwData{Armijo parameter $\sigma$}
	\KwData{maximum number of iterations~$N_{\max}$}
	\KwData{Riemannian metric $\riemannian{\cdot}{\cdot}[Q]$ on $\planarmanifold$}
	\KwData{retraction $\retract{Q}$ on $\planarmanifold$}
	\KwResult{approximate stationary point of the problem \eqref{eq:penalizedProblem} on $\planarmanifold$}

	\While{stopping criterion is not satisfied and $n < N_{\max}$}{%
		set $Q_0 \coloneqq \Qref$ and $n \coloneqq 0$

		compute the state $y$ by solving \eqref{eq:stateEquation}

		compute the adjoint state $p$ by solving \eqref{eq:adjointEquation}

		evaluate the derivative $\d_Q [j+\newpenalty](Q^n) \in \cotangent{Q^n}[\planarmanifold]$ via the left-hand side of \eqref{eq:designEquation}

		find the negative gradient $d^n \in \tangent{Q^n}[\planarmanifold]$ by solving the linear system
		\begin{equation*}
			\riemannian{d^n}{V}[Q^n]
			=
			- \d_Q [j+\newpenalty](Q^n)[V]
			\quad
			\text {for all }V \in \tangent{Q^n}[\planarmanifold]
		\end{equation*}

		find a step size $s_n$ via Armijo backtracking, satisfying
		\begin{equation}
			\label{eq:Armijo_condition}
			(j + \newpenalty)\paren[auto](){\retract{Q^n}(s_n \, d^n)}
			\le
			(j + \newpenalty)(Q^n) + \sigma \, s_n \d_Q [j + \newpenalty](Q^n)[d^n]
		\end{equation}

		update $Q^{n+1} \coloneqq \retract{Q^n}(s_n \, d^n)$

		set $n \coloneqq n+1$
	}
	\KwRet{$Q^{n+1}\in \planarmanifold$, an approximate stationary point of $j + \newpenalty$}

	\caption{General formulation of the steepest descent method on $\planarmanifold$ for \eqref{eq:penalizedProblem}.}
	\label{algorithm:general_formulation_steepest_method}
\end{algorithm2e}

As already mentioned, there exist various reasonable choices of Riemannian metrics and retractions, which we describe in what follows.
The most obvious choice is to use the Euclidean metric, which is possible because $\planarmanifold$ is an open sub-manifold of $\R^{2\times N_V}$.
In this metric, the conversion of the derivative $\d_Q$ to the gradient $\grad$ is trivial, and the geodesics are straight lines, \ie,
\begin{equation}
	\label{eq:euclideanRetraction}
	\retractEuc{Q}(V)
	=
	Q + V
\end{equation}
holds.
As a drawback, this metric is not complete and one has to take extra care not to take too large line search steps, which would lead to degenerate meshes.

Another option, proposed for instance by \cite{SchulzSiebenbornWelker:2016:1}, is to choose the Riemannian metric representing the bi-linear form associated with the Lamé system of linear elasticity.
For a point $Q\in\planarmanifold$ and tangent vectors $V,\widetilde{V}\in\tangent{Q}[\planarmanifold] = \R^{2 \times N_V}$, we represent this Riemannian metric as
\begin{equation}
	\label{eq:linear_elasticity}
	\riemannianLinElas{V}{\widetilde{V}}[Q]
	\coloneqq
	V^\transp \K \, \widetilde{V}
	+
	\delta \, V^\transp \M \, \widetilde{V}
	.
\end{equation}
The matrix $\K$ is the finite element stiffness matrix, for piecewise linear elements over the mesh defined by $Q$, associated with the linear elasticity operator
\begin{equation*}
	2 \mu \int_{\Omega_Q} \bvarepsilon(\bu) \dprod \bvarepsilon(\bv) \d x + \lambda \int_{\Omega_Q} \trace(\bvarepsilon(\bu)) \, \trace(\bvarepsilon(\bv)) \d x
	,
\end{equation*}
depending on the Lamé constants $\lambda,\mu$.
The parameter $\delta > 0$ is a damping parameter and it is required to ensure that the metric is positive definite, since we do not have a clamping boundary.
Moreover, $\M$ is the mass matrix.

It is relatively straightforward to see that \eqref{eq:linear_elasticity} is indeed a Riemannian metric on $\planarmanifold$, we refer the reader to~\cite[Theorem~5.1.1]{LoayzaRomero:2022:1} for a formal proof of this statement.
Indeed, with $\mu > 0$, $\lambda + \mu > 0$ and $\delta > 0$, \eqref{eq:linear_elasticity} is symmetric and positive definite.
Moreover, the metric coefficients vary smoothly along $\planarmanifold$.
The latter can be verified by considering the \emph{local stiffness and mass matrices} and noticing that only the transformation from the reference element to the world element depends on the node positions and it does so in a smooth manner since the connectivity remains fixed.
In practice, the elasticity metric is usually combined with the Euclidean retraction \eqref{eq:euclideanRetraction}, which may result again in a restriction of the sizes steps in order to avoid a degenerate mesh.

In this paper, we put particular emphasis on the use of a new Riemannian metric derived from the penalty function~$\newpenalty$ in \eqref{eq:penalization}.
Since the latter is proper by \cref{theorem:phiProper}, the following result is a direct consequence of \cite[Theorem~1]{Gordon:1973:1}:
\begin{proposition}
	\label{proposition:NewCompleteMetric}
	Suppose that $\alpha_j > 0$ holds, $j = 1,\ldots,4$.
	Then the manifold $\planarmanifold$, endowed with the Riemannian metric whose components (\wrt the $\vvec$ chart) are given by
	\begin{equation}
		\label{eq:phiCompleteMetric}
		g_{ab}^\textup{complete}
		=
		\delta^b_a + \frac{\partial \newpenalty}{\partial (\vvec Q)^a}\frac{\partial \newpenalty}{\partial (\vvec Q)^b}
		,
		\quad
		a, b = 1, \ldots, 2 \, N_V
		,
	\end{equation}
	is geodesically complete.
\end{proposition}
This new metric for discretized shape optimization problems compares as follows to the Euclidean and elasticity metrics recalled above.
First, exploiting the fact that \eqref{eq:phiCompleteMetric} is merely a rank-$1$ perturbation of the identity matrix, the solution of the linear system to obtain the respective gradient of the penalized objective from its derivative is very efficient; compare \cref{remark:rank1perturbationofIdentity}.
Second, we can, in principle, follow the geodesic with respect to this metric in negative gradient direction in the Armijo line search procedure.
In other words, we can use the exponential map as the retraction.
Due to the completeness of the metric, no artificial restriction of the step sizes is then required in order to avoid degenerate meshes, \ie, in order to remain on $\planarmanifold$.

The exponential map $\exponential{Q} \colon \tangent{Q}[\planarmanifold] \to \planarmanifold$ at the point $Q \in \planarmanifold$ is defined as
\begin{equation}
	\label{eq:exponential_map}
	V
	\mapsto
	\exponential{Q}{V}
	\coloneqq
	\geodesic<p>{Q}{V}(1)
	,
\end{equation}
where $\geodesic<p>{Q}{V}(t)$ denotes the geodesic, starting at $Q$ with initial velocity~$V$, evaluated at time~$t$.
In spite of the simplicity of the metric \eqref{eq:phiCompleteMetric}, the geodesic equation must be solved numerically, \eg, as described in \cite[Section~5]{HerzogLoayzaRomero:2022:1}.
In practice, as confirmed by our experiments in \cref{section:numerical_results}, this step in \cref{algorithm:general_formulation_steepest_method} is prohibitively expensive.
However, even when combined with the Euclidean retraction, the new metric \eqref{eq:phiCompleteMetric} performs very favorably in practice, at lower numerical cost than the elasticity metric.

\section{Numerical Experiments}
\label{section:numerical_results}

\begin{figure}[htp]
	\centering
	\includegraphics[width=0.6\textwidth]{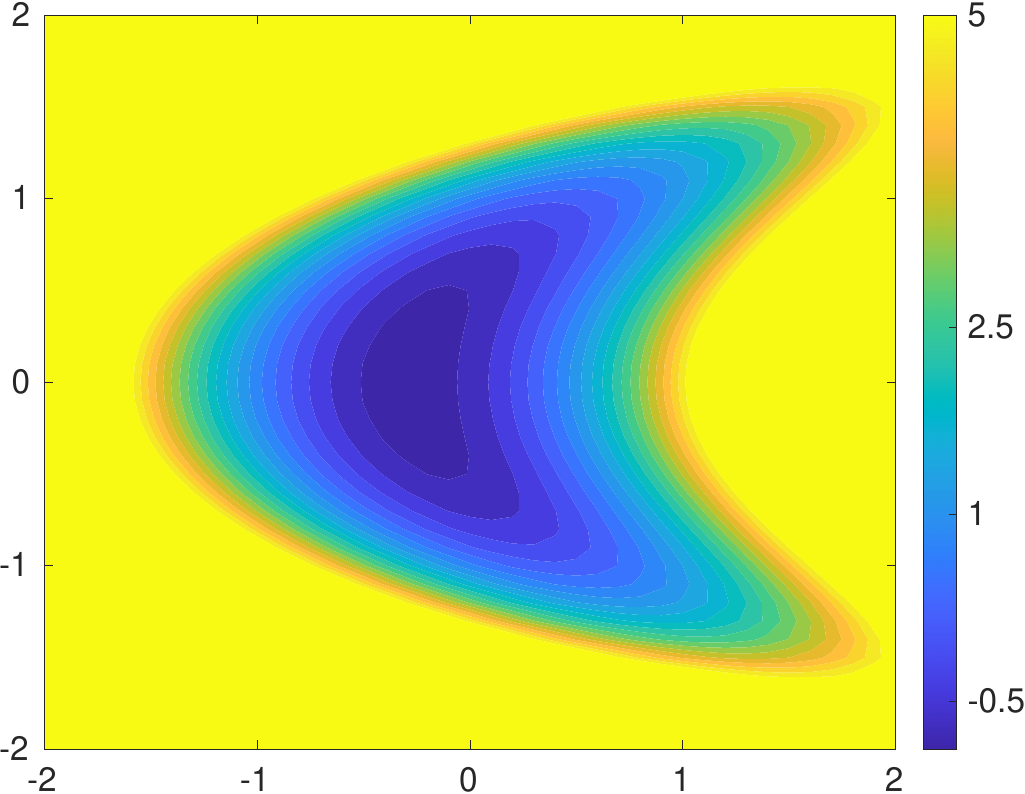}
	\caption{Contour plot of $\rhs$.}
	\label{fig:contourPlotRHS}
\end{figure}

This section aims to compare the performance of different combinations of Riemannian metrics and retractions within the steepest descent method, given in \cref{algorithm:general_formulation_steepest_method}, for the solution of a discretized, PDE-constrained shape optimization problem.
For the first three out of four experiments, we stick to the model problem \eqref{eq:penalizedProblem} with right-hand side $\rhs(x_1,x_2) = 2.5 \, (x_1 + 0.4 - x_2^2)^2 + x_1^2 + x_2^2 - 1$, as previously used in \cite{EtlingHerzogLoayzaWachsmuth:2020:1}.
As was observed in \cite{BartelsWachsmuth:2020:1}, the main motivation for this choice is its simple interpretation of the expected solution.
Recall that our goal is to minimize $\int_{\Omega_Q} y \d x$ and notice that the  sublevel set $\setDef{x \in \R^2}{\rhs(x) \le 0}$ is connected.
Due to the maximum principle, we can therefore expect to find an optimal shape close to this sublevel set, at least in the continuous setting where a maximum principle is available.
In the discrete setting, the maximum principle hinges upon the condition of non-obtuse triangles, which is not guaranteed a-priori.
Indeed, we did find obtuse triangles in most our experiments to occur.
\Cref{fig:contourPlotRHS} shows a contour plot of $\rhs$ for comparison with the optimal shapes obtained.

The variants we compare in what follows are termed \EucEuc, \ElasEuc, \CompEuc and \CompComp.
The first component of the name refers to the metric used for the evaluation of the shape gradient; see \eqref{eq:definition_riemannian_gradient}.
The three choices indicate the Euclidean metric, the elasticity metric \eqref{eq:linear_elasticity} and the new complete metric \eqref{eq:phiCompleteMetric}.
Their precise parameters are specified further below.
The second component of the name refers to the choice of the retraction, which is either Euclidean \eqref{eq:euclideanRetraction} or the exponential map \eqref{eq:exponential_map}, evaluated via numerical integration as detailed in \cite[Section~5]{HerzogLoayzaRomero:2022:1}.

This section is structured as follows.
In \cref{subsection:implementation_details}, we describe the implementation details used throughout the numerical experiments.
Four experiments are then conducted to explore various points.
First, we consider problem \eqref{eq:penalizedProblem} without a penalty term in \cref{subsection:unpenalizedAndIdeal}.
We confirm that, as expected, this problem then does not possess a solution.
Consequently, this leads any gradient descent method, regardless of the metric employed, to ultimately produce degenerate meshes in the pursuit of smaller and smaller objective values.
However, the variants \ElasEuc, and \CompComp still produce \enquote{good} iterates along the way, albeit at different iteration counts, while \EucEuc breaks down early.

Our second experiment in \cref{subsection:experiment_diff_penalization_parameters} targets the penalized problem, for which the existence of a solution can be proved.
It turns out that now, as expected, the gradient descent method finds this solution regardless of the metric chosen, yet at different iteration numbers.
Computationally, we observe that the new metric, combined with the Euclidean retraction (\CompEuc), outperforms \EucEuc and also \ElasEuc for the problem under consideration.

The results this far indicate that the typically ill-posed problem of minimizing a discrete shape optimization objective may be tackled either by early stopping or by the addition of a penalty term.
The penalty approach may be criticized since it requires the user to make a somewhat arbitrary choice of parameters $\alpha_1$, $\alpha_2$, $\alpha_3$, $\alpha_4$.
In our third experiment in \cref{subsection:unpenalizedDiffMeshLevels}, we therefore revisit the first strategy and compare the two most promising candidates, \ElasEuc and \CompEuc, using finer meshes than before.
Once again, it turns out that the use of the new metric may maintain better-quality meshes and requires less time per iteration compared to the elasticity metric.

We found numerically that the hold-all domain assumption required for the proof of \cref{proposition:existenceSolutions} did not require to be enforced.

Finally, we consider in \cref{subsection:experimentOptimalBridge} a classical compliance minimization problem.
We solve the penalized problem, comparing again the two most promising variants \ElasEuc and \CompEuc.
Again, the numerical results indicate that using the proposed complete metric, one can obtain similar results as with the \ElasEuc variant, but in less time and with meshes of better quality.

\subsection{Implementation Details}
\label{subsection:implementation_details}

Our implementation is achieved in \matlab, using the \lstinline!initmesh! function of the PDE toolbox for the generation of all initial meshes and the code provided by \cite{Koko:2016:1,Koko:2016:2} to assemble the elasticity stiffness and mass matrices required for the elasticity metric \eqref{eq:linear_elasticity}.
All experiments were performed on a computer with an Intel Core~i7-7500 CPU with 2.7~GHz and 16GiB~RAM.

\subsubsection*{Initialization of the Armijo Backtracking Procedure:}
As already described in \cref{algorithm:general_formulation_steepest_method}, we use Armijo's condition \eqref{eq:Armijo_condition} in order to guarantee sufficient decreasing of the (penalized) objective function through a backtracking procedure.
It is well-known that the steepest descent method is not scale invariant and therefore relies on a judicious choice of the initial line search step size.
We use the technique presented in \cite[p.59]{NocedalWright:2006:1}, \ie, the candidate for the initial step size in iteration~$n$ is given by
\begin{equation*}
	\overline{s}_n
	=
	s_{n-1} \frac{\d_Q [j+\newpenalty](Q^{n-1})[d^{n-1}]}{\d_Q [j+\newpenalty](Q^n)[d^n]}
	.
\end{equation*}
This candidate step size gets overwritten in the initial iteration or when $\overline{s}_n$ becomes too small.
We use the rule
\begin{equation}
	\label{eq:definition_sbar}
	s^\text{initial}_n
	=
	\begin{cases}
		\frac{1}{\riemanniannorm{d^n}[Q^n]}
		&
		\text{if $n=0$ or $\overline{s}_n \, \riemanniannorm{d^n}[Q^n] < 10^{-4}$}
		,
		\\
		\overline{s}_n
		&
		\text{otherwise}
	\end{cases}
\end{equation}
for this purpose.
We denote by $ \riemanniannorm{d^n}[Q^n]$ the Riemannian norm at the point $Q^n$,  \ie, $\riemanniannorm{d^n}[Q^n] = \sqrt{\riemannian{d^n}{d^n}}$.
Should a trial step size fail to satisfy the Armijo condition \eqref{eq:Armijo_condition}, we repeatedly multiply it by a factor $\tau \in (0,1)$ specified further below.

We recall that some of the variants of the algorithm involve the Euclidean retractions \eqref{eq:euclideanRetraction}.
In this case, mesh vertices move independently of each other and thus extra care needs to be taken regarding the trial step sizes in order to avoid degenerate meshes.
We proceed as follows.
When the Euclidean retraction is used, we consider the Armijo condition \eqref{eq:Armijo_condition} failed for the trial step size~$s$ as long as the distance a vertex would travel is relatively large compared to the heights of any of its incident triangles.
More precisely, we treat the Armijo condition as failed as long as
\begin{equation}
	\label{eq:modified_initial_s_euclidean_retraction}
	s \, \euclideanriemanniannorm{d_{i^k_\ell}}[Q]
	\ge
	0.5 \, \height{Q}{\ell}(i^k_0,i^k_1,i^k_2)
	\quad
	\text{for any $k=1,\ldots,N_T$ and any $\ell = 0,1,2$}
\end{equation}
holds.
Here $\euclideanriemanniannorm{d_{i^k_\ell}}[Q]$ denotes the Euclidean norm of the subvector of the negative gradient direction~$d$ pertaining to the $\ell$-th vertex of the $k$-th triangle, and $\height{Q}{\ell}$ is the corresponding height, see \eqref{eq:triangle_height}.
(For the purpose of readability, we temporarily dropped the iteration index~$n$ here.)

\subsubsection*{Armijo Backtracking with the Exponential Map:}
In the experiment in \cref{subsection:unpenalizedAndIdeal}, we use \CompComp as one of the variants of \cref{algorithm:general_formulation_steepest_method}.
As opposed to all other variants using the Euclidean retraction, the geodesic equation with respect to the metric \eqref{eq:phiCompleteMetric} must be integrated numerically, which is expected to be expensive.
In order to avoid repeated evaluations of the geodesic in case the Armijo's condition \eqref{eq:Armijo_condition} happens to fail for the initial trial step size, we make use of the re-scaling lemma; see, \eg, \cite[Lemma~5.18, p.127]{Lee:2018:1}.
In our context, this lemma states that for every initial data $Q\in\planarmanifold$ and $d\in \tangent{Q}[\planarmanifold]$, trial step size and backtracking parameter $\tau > 0$, we have $\geodesic<p>{Q}{\tau \, s \, d}(1) = \geodesic<p>{Q}{s \, d}(\tau)$.
When integrating the initial trial geodesic with velocity $s^\textup{initial} d$ until $t = 1$, our implementation of the numerical integrator thus stores the values at $t \in \{\tau, \tau^2, \ldots\}$.
This can be conveniently achieved by setting $\tau = 0.5$ and using a number of time steps divisible by a sufficiently large power of~2.

\subsubsection*{Parameter Choices:}
We keep the following parameters fixed for all experiments.
For the Armijo line search, we use the acceptance and backtracking parameters $\sigma = 10^{-4}$ and $\tau = 0.5$.
The linear elasticity metric given in \eqref{eq:linear_elasticity} uses Lamé constants given by
\begin{equation}
	\label{eq:Lame_parameters}
	\mu = \frac{E}{2(1+\nu)}
	,
	\quad
	\lambda = \frac{E \, \nu}{(1+\nu)(1-2\nu)}
	,
	\quad
	\delta = 0.2 \, E
	,
\end{equation}
with Young's modulus $E = 1$ and Poisson ratio $\nu = 0.4$.

As a measure of the quality of the meshes generated, we monitor the function
\begin{equation}
	\label{eq:quality_measure}
	\Theta(Q)
	=
	\sum_{k=1}^{N_T} \frac{1}{N_T} \frac{1}{\psi_Q(i_0^k,i_1^k,i_2^k)}
	,
\end{equation}
which is part of the penalty function's definition \eqref{eq:penalization}, where $\nicefrac{1}{\psi}$ is given by \eqref{eq:Shewchuk64_qualitymeasure}.
We remind the reader that $\Theta(Q) \ge 1$ holds, and $1$ constitutes the best value while bad quality meshes correspond to large values of $\Theta$.

We also recall that the penalty function~$\newpenalty$ serves two purposes: its addition to the objective renders the penalized problem well-posed, and it forms the basis for the complete metric \eqref{eq:phiCompleteMetric}.
For flexibility, we allow two different sets of parameters $\alpha_j$, $j = 1, \ldots, 4$ for both occurrences.
They are denoted as $\alphaPenalization_j$ and $\alphaMetric_j$, respectively.
For the problems under consideration, we do not run the risk of exterior self-intersections so we set $\alphaPenalization_3 = \alphaMetric_3 = 0$ for all experiments.
This can be justified using a thresholding function as in \cref{remark:addCutOffAugmentation}.
The remaining parameters are specified in each of the following sections as needed.

\subsubsection*{Derivative-Gradient Transformation:}
The evaluation of the gradient \eqref{eq:definition_riemannian_gradient} requires the solution of a linear system whenever the metric is not the Euclidean one.
In case of the linear elasticity metric \eqref{eq:linear_elasticity}, we assemble the stiffness and mass matrices using the code provided by \cite{Koko:2016:1,Koko:2016:2}.
The subsequent solve of the linear system was achieved using the default sparse direct solver of \matlab.
For the moderate size of the experiments conducted, a more sophisticated strategy such as a geometric multigrid method does not pay off.

In case of the complete metric \eqref{eq:phiCompleteMetric}, we exploit the fact that the associated matrix is is a rank-$1$ perturbation of the identity matrix.
We therefore solve the linear system \eqref{eq:definition_riemannian_gradient} using two iterations of the conjugate gradient method without preconditioning, which is sufficient for convergence; see \eg, \cite[eq.~(2.11), p.76]{ElmanSilvesterWathen:2014:1}.
Our implementation is matrix-free, \ie, we provide only matrix-vector products.
The most expensive part of this process is the evaluation of the first-order derivatives of the penalty function~$\newpenalty$.

\subsubsection*{Definition of Unsuccessful Experiments:}
As a precautionary measure, we keep track of several indicators during the iteration of the gradient descent method \cref{algorithm:general_formulation_steepest_method}.
In particular, we verify that each search direction $d^n$ is indeed a descent direction, \ie, $\d_{Q^n} [j + \newpenalty](Q^n)[d^n] < 0$ holds.
Moreover, we make sure that the signed areas \eqref{eq:signed_area} of all triangles remain positive for all iterates, which is a requirement for them to belong to the manifold $\planarmanifold$.
As expected, these indicators were never found to be violated.

It can happen, however, that a close-to-degenerate mesh enforces very small trial step sizes due to \eqref{eq:modified_initial_s_euclidean_retraction} when the Euclidean retraction is used.
Indeed, we declare a gradient descent run unsuccessful and stop as soon as a trial step size becomes smaller than $10^{-7}$.
In our experiments we only observed this in case of the \EucEuc variant.

\subsubsection*{Stopping Criteria:}
Choosing a stopping criterion is a delicate task.
This is especially true in case of the unpenalized problem, which may not possess solutions, and early stopping (before the norm of the gradient becomes too small) becomes essential.
Since the attempt to approximate the infimum results in degenerate meshes, using any criterion involving the value of the objective alone will also not be suitable.
As a compromise, we therefore settle on a fixed number of iterations for the experiments in \cref{subsection:unpenalizedAndIdeal,subsection:unpenalizedDiffMeshLevels}, which concern the unpenalized problem.

For the penalized problems in \cref{subsection:experiment_diff_penalization_parameters,subsection:experimentOptimalBridge}, which do have a solution, we can use a more classic approach.
Since we compare different metrics, which entail different ways to measure the norm of the gradient, the gradient norm does not allow a fair comparison.
We therefore resort to measuring the absolute change of the values of the penalized objective function over a span of $5$~past iterations, and use it as an stopping criterion.
This results in stopping as soon as
\begin{equation}
	\label{eq:stoppingCriterionPenalized}
	\max_{m = 1,\ldots, 5} \paren[big]\{\}{(j+ \newpenalty)(Q^{n-m})-(j + \newpenalty)(Q^n)}
	<
	\textup{tol}
	.
\end{equation}
This is motivated by a condition proposed in \cite[Section~6.15, p.1324]{Laurain:2018:1}.

\subsection{Experiment~1: Lack of Solutions for the Unpenalized Problem}
\label{subsection:unpenalizedAndIdeal}

As was argued in \cref{subsection:first_glimpse_non-existence}, discretized shape optimization problems in which the vertex positions serve as optimization variables can not be expected to possess a solution.
Here we confirm this observation for our model problem \eqref{eq:penalizedProblem} without a penalty, \ie, we set $\alphaPenalization_j = 0$ for all $j=1,2,3,4$.

Consequently, this leads any gradient descent method, regardless of the metric employed, to ultimately produce degenerate meshes in the pursuit of smaller and smaller objective values.
We also trace back the specific nature of the degeneracy observed to an exploitation of the quadrature formula for the problem at hand.

\begin{figure}[htp]
	\begin{center}
		\captionsetup[subfigure]{labelformat=empty}
		\begin{subfigure}{0.9\textwidth}
			\begin{center}
				\raisebox{5mm}{%
					\includegraphics[width=0.3\textwidth]{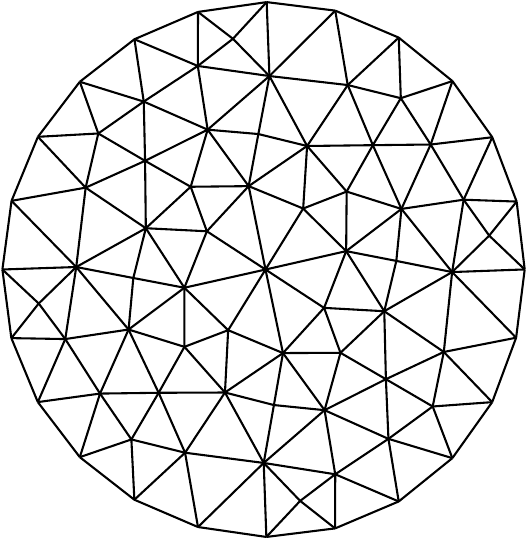}
				}
				\quad
				\includegraphics[width=0.30\textwidth]{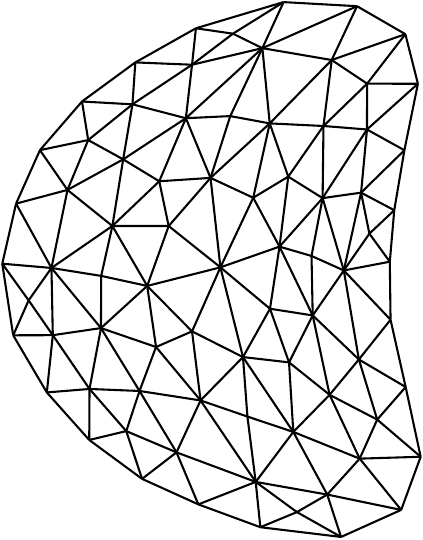}
				\quad
				\includegraphics[width=0.30\textwidth]{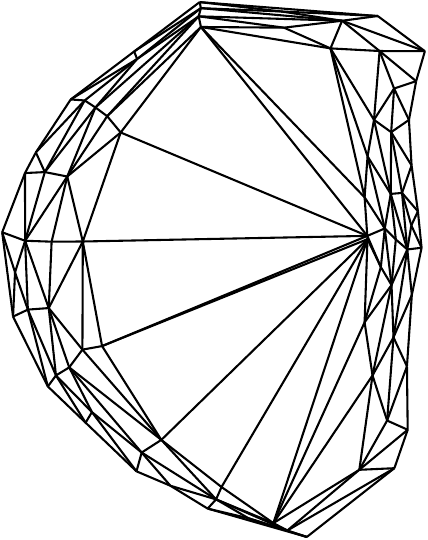}
				\caption{Iterates from left to right: $15$ (too early), $150$ (good), $1000$ (too late) for variant \ElasEuc.}
				\label{fig:experimentNoPenalizationvariant2}
			\end{center}
		\end{subfigure}
		\\
		\begin{subfigure}{0.9\textwidth}
			\begin{center}
				\includegraphics[width=0.30\textwidth]{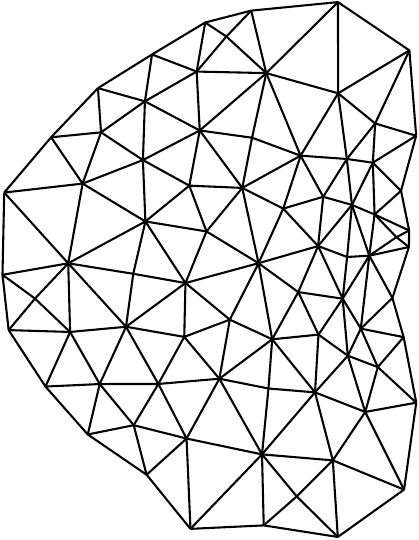}
				\quad
				\includegraphics[width=0.30\textwidth]{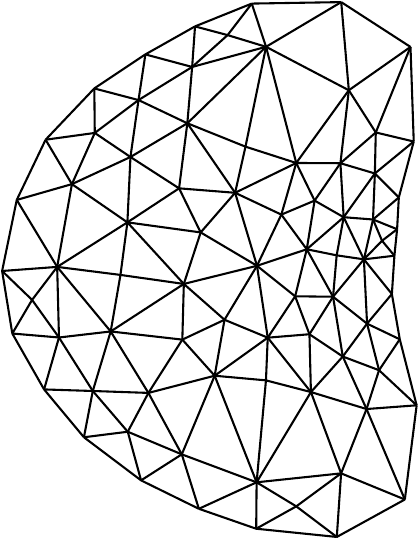}
				\quad
				\includegraphics[width=0.30\textwidth]{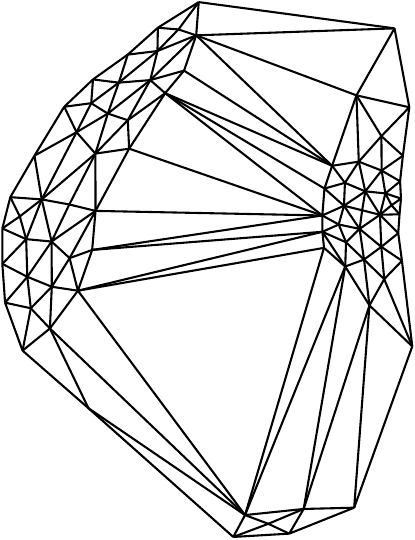}
				\caption{Iterates from left to right: $5$ (too early), $15$ (good), $1000$ (too late) for variant \CompComp.}
				\label{fig:experimentNoPenalizationvariant4}
			\end{center}
		\end{subfigure}
		\\
		\begin{subfigure}{0.45\textwidth}
			\begin{center}
				\pgfplotsset{compat=1.14, width=\textwidth}
				\begin{tikzpicture}
					\begin{axis}[
					title = Objective,
						legend style={ draw=none, at={(0,-0.2)},anchor=north west}
						]
						\addplot[color=black, dashed, very thick] table [x=iter,y=Obj, col sep=comma]{matlab/dataResults/experimentNoPenalization/meshLevel0/metricType_1_retractionType_1_DATA.txt};
						\addlegendentry{\EucEuc}
						\addplot[color=blue, dotted, very thick] table [x=iter,y=Obj, col sep=comma]{matlab/dataResults/experimentNoPenalization/meshLevel0/metricType_2_retractionType_1_DATA.txt};
						\addlegendentry{\ElasEuc}
						\addplot[color=magenta, very thick] table [x=iter,y=Obj, col sep=comma]{matlab/dataResults/experimentNoPenalization/meshLevel0/metricType_3_retractionType_2_DATA.txt};
						\addlegendentry{\CompComp}
						\end{axis}
					\end{tikzpicture}
					\label{fig:experimentNoPenalizationObjective}
				\end{center}
		\end{subfigure}
		\hfill
		\begin{subfigure}{0.45\textwidth}
			\begin{center}
				\pgfplotsset{compat=1.14, width=\textwidth}
				\begin{tikzpicture}
					\begin{axis}[
					title = Mesh Quality,
						legend style={ draw=none,at={(0,-0.2)},anchor=north west},
						ymax = 10
						]
						\addplot[color=black, dashed, very thick] table [x=iter,y=mshQua, col sep=comma]{matlab/dataResults/experimentNoPenalization/meshLevel0/metricType_1_retractionType_1_DATA.txt};
						\addlegendentry{\EucEuc}
						\addplot[color=blue, dotted, very thick] table [x=iter,y=mshQua, col sep=comma]{matlab/dataResults/experimentNoPenalization/meshLevel0/metricType_2_retractionType_1_DATA.txt};
						\addlegendentry{\ElasEuc}
						\addplot[color=magenta, very thick] table [x=iter,y=mshQua, col sep=comma]{matlab/dataResults/experimentNoPenalization/meshLevel0/metricType_3_retractionType_2_DATA.txt};
						\addlegendentry{\CompComp}
					\end{axis}
				\end{tikzpicture}
				\label{fig:experimentNoPenalizationMeshQuality}
			\end{center}
		\end{subfigure}
	\end{center}
	\caption{Results for the experiment described in \cref{subsection:unpenalizedAndIdeal}.}
	\label{fig:experimentNoPenalization}
\end{figure}
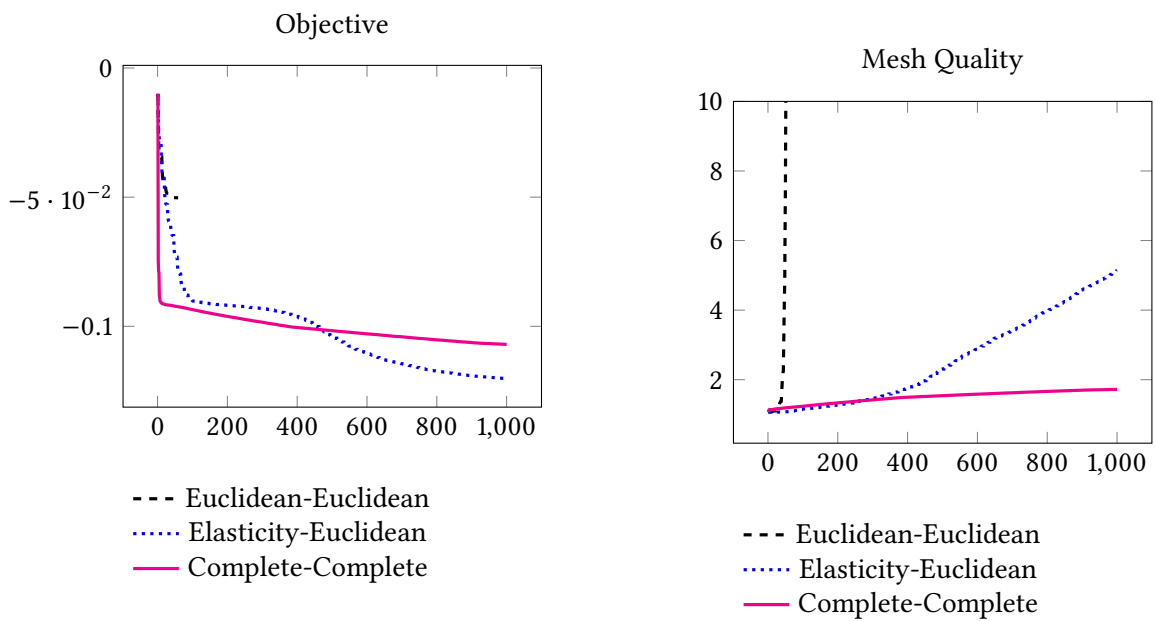

We compare the variants \EucEuc, \ElasEuc, and \CompComp.
For the latter, we use the parameters $\alphaMetric_1 = 10$, $\alphaMetric_2 = 1 $, $\alphaMetric_3 = 0$ and $\alphaMetric_4 = 0.01$.
The initial mesh for this first experiment is a coarse triangulation of the unit disc containing $N_V = 77$~vertices and $N_T = 128$~triangles.
The results are shown in \cref{fig:experimentNoPenalization} and \cref{table:experimentNoPenalization}.
The \EucEuc variant breaks down in iteration~$60$ with too small a trial step size and a disastrous value of the mesh quality measure~$\Theta$ from \eqref{eq:quality_measure} and it is thus evaluated as an unsuccessful experiment.
By contrast, the \ElasEuc and \CompComp variants produce meshes of comparable quality and similarly small values of the objective at iteration counts~$150$ and $15$, respectively.
As expected, both enter a phase of producing increasingly degenerate meshes afterwards before being stopped at iteration~$1000$.
However, we observe that the deterioration of the mesh quality is more pronounced for the \ElasEuc variant.

\begin{table}[htp]
	\begin{center}
	\begin{tabular}{c}
		\toprule
		variant
		\\
		\midrule
		\EucEuc
		\\
		\ElasEuc
		\\
		\CompComp
		\\
		\bottomrule
	\end{tabular}%
	\begin{filecontents*}{matlab/dataResults/experimentNoPenalization/meshLevel0/summaryExperimentOne.txt}
iter,counterArm,counterGeo,Obj,penObj,mshQua,normGrad,StopCrit
59,0,17,-0.0502331933042709,-0.0502331933042709,160.437402531471,0.129526167390427,9.88284894753999e-06
1000,0,0,-0.120261467646551,-0.120261467646551,5.14271156826822,0.0170106442500635,5.70868802818036e-05
1000,0,0,-0.107002815578143,-0.107002815578143,1.72096233313203,0.00389155381043758,3.07605940185057e-05
\end{filecontents*}
\pgfplotstabletypeset[%
	col sep = comma,
	every head row/.style={before row=\toprule,after row=\midrule},
	every last row/.style={after row=\bottomrule},
	columns={iter,Obj,mshQua},
	columns/iter/.style={int detect, column name={iter $(n)$}, column type = {r}},
	columns/Obj/.style={std, precision=4,column name={$j(Q^n)$}, column type = {r}},
	columns/mshQua/.style={std, precision=4,column name={$\Theta(Q^n)$}, column type = {r}}%
	]{matlab/dataResults/experimentNoPenalization/meshLevel0/summaryExperimentOne.txt}
	\end{center}
	\caption{Summary of the results obtained for the experiment described in \cref{subsection:unpenalizedAndIdeal}.}
	\label{table:experimentNoPenalization}
\end{table}

\begin{figure}[htp]
	\begin{center}
		\begin{subfigure}{0.45\textwidth}
			\centering
			\includegraphics[width=\textwidth]{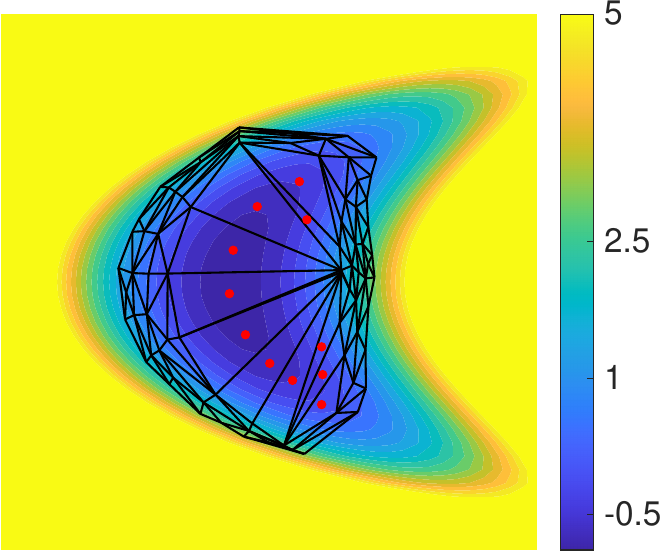}
			\caption{Variant \ElasEuc.}
		\end{subfigure}
		\hfill
		\begin{subfigure}{0.45\textwidth}
			\centering
			\includegraphics[width=\textwidth]{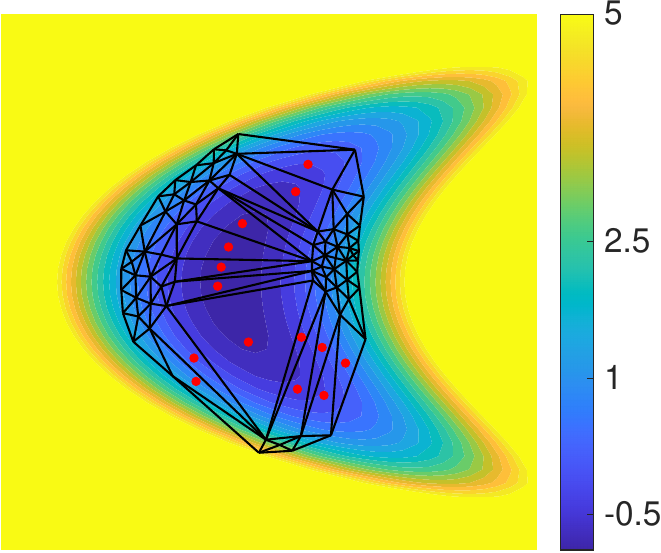}
			\caption{Variant \CompComp.}
		\end{subfigure}
		\caption{Location of the centers of the larger triangles at iterate~$1000$, superimposed on a contour plot of the right-hand side $\rhs$ for the experiment described in~\cref{subsection:unpenalizedAndIdeal}.}
		\label{fig:rhsMidpointsFinalMeshes}
	\end{center}
\end{figure}

As announced earlier, it is illustrative to study the meshes for the \ElasEuc and \CompComp variants at the final iteration~$1000$.
As shown in \cref{fig:rhsMidpointsFinalMeshes}, large triangles are produced where the values of the PDE's right-hand side function~$\rhs$ are smallest.
This is due to the discrete objective involving a quadrature formula for the evaluation of the element load vector, which evaluates the right-hand side only in the triangle centers, some of which are marked by red dots.
Given the opportunity, it thus can be concluded that the optimizer exploits the quadrature error.

A first conclusion at this point is that a gradient method, applied to an unpenalized problem without a solution, might be successful to produce a reasonably good approximation to the solution of the continuous shape optimization problem, provided that it is stopped sufficiently early.
As already noted, the \CompComp variant reaches this convenient stopping point at a much earlier iteration number.
However, the picture changes when comparing the respective run-times.

\Cref{table:executionTimesvariant4} shows the timings for each of the variants.
The column \texttt{state} summarizes the time devoted to solving the state equation at least once per iteration, depending on the number of Armijo backtracking steps.
The column \texttt{dObj} represents the time invested in assembling the derivative of the shape derivative.
Likewise, the column \texttt{backt} presents the time required to check whether the line search trial step sizes do not  satisfy \eqref{eq:modified_initial_s_euclidean_retraction} (in case of the Euclidean retraction) or do satisfy the Armijo condition \eqref{eq:Armijo_condition}.
The column \texttt{grad} shows the time needed in the transformation of the derivative to the gradient, \ie, for the solution of the linear system \eqref{eq:definition_riemannian_gradient}.
Finally, the column \texttt{retr} shows the time to evaluate the retraction.
This is not relevant for the Euclidean retraction, but only in case the geodesic equation associated with the metric \eqref{eq:phiCompleteMetric} is solved numerically.
The latter is achieved using the implementation of the Störmer--Verlet scheme detailed in \cite[Section~5]{HerzogLoayzaRomero:2022:1}.
We used $1024$~time steps for this purpose to ensure convergence of the fixed-point solver for the implicit sub-step.

As the timings clearly show, the numerical integration of the geodesic equation associated with the metric \eqref{eq:phiCompleteMetric} is prohibitively expensive in the \CompComp variant.
Therefore, we replace the \CompComp variant by \CompEuc for further experiments, \ie, we combine the metric \eqref{eq:phiCompleteMetric} with the Euclidean retraction.

\begin{table}[htp]
	\begin{center}
	\begin{tabular}{c}
	 	\toprule
	  Variant  \\
	 	\midrule
	 	\EucEuc   \\
	 	\ElasEuc  \\
	  \CompComp \\
		\bottomrule
	 \end{tabular}%
	\sisetup{round-mode = places, round-precision = 3, scientific-notation = fixed, fixed-exponent = 0}%
	\begin{filecontents*}{matlab/timesOfExecutionExperimentOne.txt}
iter,total,State,dObj,backt,grad,retr
60,0.83554,0.30815,0.086402,0.23711,0,0
150,1.7264,0.70648,0.16552,0.53427,0.18137,0
15,89.761,0.13892,0.031043,0.083803,0.030048,89.421
\end{filecontents*}
\pgfplotstabletypeset[%
	col sep = comma,
	every head row/.style={before row=\toprule,after row=\midrule},
	every last row/.style={after row=\bottomrule},
	columns={iter,total,perIter,State,dObj,backt,grad,retr},
	columns/iter/.style={int detect, column name={iter}},
	columns/total/.style={std, precision=2,column name={total},
	assign cell content/.code={%
	\pgfkeyssetvalue{/pgfplots/table/@cell content}%
	{\SI{##1}{\second}}%
	}},
	columns/perIter/.style={fixed, column name={per iter},
	assign cell content/.code={%
	\pgfkeyssetvalue{/pgfplots/table/@cell content}%
	{\SI{##1}{\second}}%
	}},
	columns/State/.style={std, precision=2,column name={state},
	assign cell content/.code={%
	\pgfkeyssetvalue{/pgfplots/table/@cell content}%
	{\SI{##1}{\second}}%
	}},
	columns/dObj/.style={std, precision=2,column name={dObj},
	assign cell content/.code={%
	\pgfkeyssetvalue{/pgfplots/table/@cell content}%
	{\SI{##1}{\second}}%
	}},
	columns/backt/.style={std, precision=2,column name={backt},
	assign cell content/.code={%
	\pgfkeyssetvalue{/pgfplots/table/@cell content}%
	{\SI{##1}{\second}}%
	}},
	columns/grad/.style={std, precision=2,column name={grad},
	assign cell content/.code={%
	\pgfkeyssetvalue{/pgfplots/table/@cell content}%
	{\SI{##1}{\second}}%
	}},
	columns/retr/.style={std, precision=2,column name={retr},
	assign cell content/.code={%
	\pgfkeyssetvalue{/pgfplots/table/@cell content}%
	{\SI{##1}{\second}}%
	}}
	]{matlab/timesOfExecutionExperimentOne.txt}
	\end{center}
	\caption{Execution times for the experiment described in \cref{subsection:unpenalizedAndIdeal}.}
	\label{table:executionTimesvariant4}
\end{table}

\subsection{Experiment~2: Solving the Penalized Problem}
\label{subsection:experiment_diff_penalization_parameters}

Our second experiment targets the penalized problem, for which the existence of a solution was proved in \cref{proposition:existenceSolutions}.
Due to the excessive time associated with the numerical integration of the geodesic equation associated with the metric \eqref{eq:phiCompleteMetric}, we consider only the Euclidean retraction \eqref{eq:euclideanRetraction} from now on.
We thus compare the variants \EucEuc, \ElasEuc and \CompEuc.
We solve the penalized problem with three different sets of parameters given in \cref{tab:parametersExpTwo}.
The initial mesh is again a coarse triangulation of the unit disc containing $N_V = 146$~vertices and $N_T = 258$~triangles.
The parameters for the metric \eqref{eq:phiCompleteMetric} $\alphaMetric_j$ are the same as in \cref{subsection:unpenalizedAndIdeal}, \ie, $\alphaMetric_1 = 10$, $\alphaMetric_2 = 1 $, $\alphaMetric_3 = 0$ and $\alphaMetric_4 = 0.01$.

\begin{table}[htp]
	\centering
	\begin{tabular}{c|rrrr}
		\toprule
		Parameter set & $\alphaPenalization_1$ & $\alphaPenalization_2$ & $\alphaPenalization_3$ & $\alphaPenalization_4$
		\\
		\midrule
		1 & 1 & 0.5 & 0.0 & 0.1
		\\
		2 & 0.1 & 0.01 & 0.0 & 0.001
		\\
		3 & 0.015 & 0.005 & 0.0 & 0.0005
		\\
		\bottomrule
	\end{tabular}
	\caption{Description of the parameter set for the experiment in \cref{subsection:experiment_diff_penalization_parameters}.}
	\label{tab:parametersExpTwo}
\end{table}

Since we know that the problem has a solution, we can use the stopping criterion in \eqref{eq:stoppingCriterionPenalized} with a tolerance of $\textup{tol} = 10^{-6}$.
The number of iterations and the final values of the objective and the penalty functionals are shown in \cref{table:experiment_diff_penalization_parameters}.
\Cref{fig:experiment_diff_penalization_parameters} shows the final iterates obtained for each variant, which are very similar to each other.

The first fact to highlight is that, variant \EucEuc, performs surprisingly well on the penalized problem, even for moderately small values of the penalty parameters $\alphaPenalization_j$ (parameter sets~1 and 2).
However, it does not quite converge within $1000$~iterations for parameter set~3.
Variants \ElasEuc and \CompEuc perform equally well, but the latter is faster; see \cref{table:executionTimesvariant3}.
Both variants are also comparable to each other and better compared to \EucEuc with respect to the values of the objective and the mesh quality, as shown in \cref{fig:experimentPenalization}.

We also mention that the evaluation of the derivative of the penalty function (column \texttt{dPen}), which might be a concern, does not require a major computational effort, at least not for the meshes of this size.

In conclusion, we find that the presence of the penalty terms helps preserve the mesh quality for all variants.
The variant \CompEuc performs fastest, at a numerical cost very close to that of \EucEuc.
This is partly due to the small cost of solving for the gradient, see \eqref{eq:definition_riemannian_gradient}.
Admittedly, the differences are small for the coarse mesh under consideration.
Therefore, we conduct a series of experiments in the following \cref{subsection:unpenalizedDiffMeshLevels} with finer meshes.

\begin{table}[htp]
	\centering
	\begin{tabular}{c|c}
		\toprule
		Parameter set & Variant \\
		\midrule
		\multirow{4}{*}{$1$}
		& \EucEuc    \\
		& \ElasEuc   \\
		& \CompEuc   \\
		\midrule
		\multirow{4}{*}{$2$}
		& \EucEuc    \\
		& \ElasEuc   \\
		& \CompEuc   \\
		\midrule
		\multirow{4}{*}{$3$}
		& \EucEuc  \\
		& \ElasEuc  \\
		& \CompEuc  \\
		\bottomrule
	\end{tabular}%
	\begin{filecontents*}{matlab/dataResults/experimentPenalization/summaryExperimentTwo.txt}
iter,counterArm,counterGeo,Obj,penObj,mshQua,normGrad,StopCrit
56,6,2,-0.0563708359533044,1.15792646384711,1.04200806038347,0.000285543447969675,5.08032640844647e-07
87,7,5,-0.0561937457823599,1.15793043980861,1.04197165274931,0.00158467857295704,8.68932340969408e-07
59,0,0,-0.0563225908292477,1.15792620959847,1.04201663956548,0.00039266237298329,9.47430604192334e-07
363,0,0,-0.0910185825733693,0.0192567108166872,1.05027455471518,0.00209190760609928,8.21982035145719e-07
261,0,0,-0.0910099324352734,0.0192910418019126,1.05037219523713,0.00164512006306142,9.69241011228306e-07
281,0,0,-0.0909761867506215,0.0192602005889455,1.05006377414496,0.00107223350877718,8.07253042703349e-07
1000,0,0,-0.0919308807596466,-0.0729003055798432,1.11648714257623,0.0019316459973981,0.000148995027410823
276,0,0,-0.0921313805239697,-0.0732629548680937,1.08948813000707,0.00270413454101796,8.10228572276084e-07
289,1,0,-0.0923432617494091,-0.0733954702800705,1.09494612129892,0.000408742297522931,8.3199522143973e-07
\end{filecontents*}
\pgfplotstabletypeset[%
	col sep = comma,
	every head row/.style={before row=\toprule,after row=\midrule},
	every last row/.style={after row=\bottomrule},
	every nth row={3}{before row=\midrule},
	columns={iter,Obj,penObj,mshQua},
	columns/iter/.style={int detect, column name={iter $(n)$}, column type = {r}},
	columns/Obj/.style={std, precision=4,column name={$j(Q^n)$}, column type = {r}},
	columns/penObj/.style={std, precision=4,column name={$j(Q^n) + \varphi(Q^n)$}, column type = {r}},
	columns/mshQua/.style={fixed, fixed zerofill = true, precision = 4,column name={$\Theta(Q^n)$}, column type = {r}}%
	]{matlab/dataResults/experimentPenalization/summaryExperimentTwo.txt}
	\caption{Summary of the results obtained for the experiment described in \cref{subsection:experiment_diff_penalization_parameters}.}
	\label{table:experiment_diff_penalization_parameters}
\end{table}

\begin{figure}[htp]
	\begin{center}
		\begin{subfigure}{0.75\textwidth}
			\begin{center}
				\includegraphics[scale=0.32]{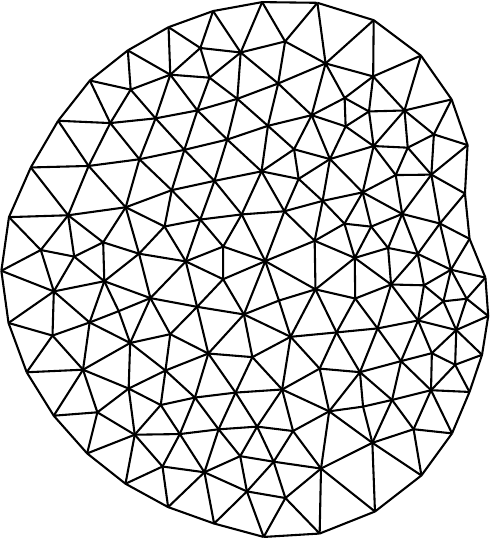}
				\quad
				\includegraphics[scale=0.32]{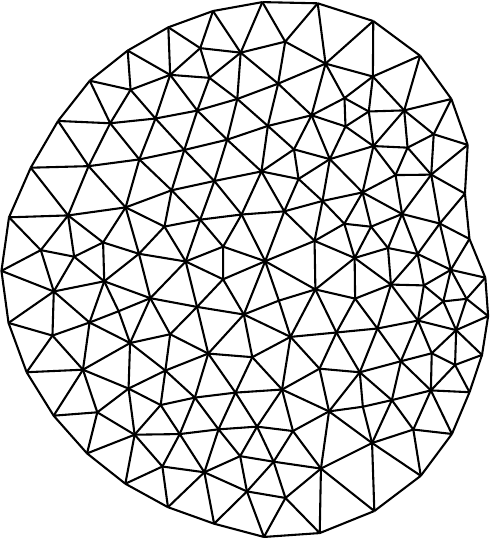}
				\quad
				\includegraphics[scale=0.32]{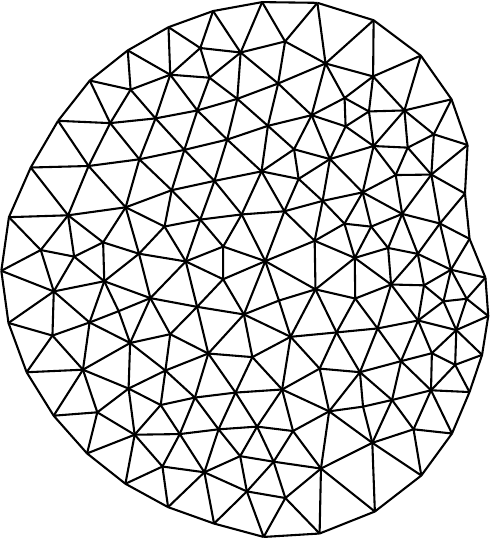}
				\caption{Parameter set $1$.}
			\end{center}
		\end{subfigure}
		\\
		\begin{subfigure}{0.75\textwidth}
			\begin{center}
				\includegraphics[scale=0.32]{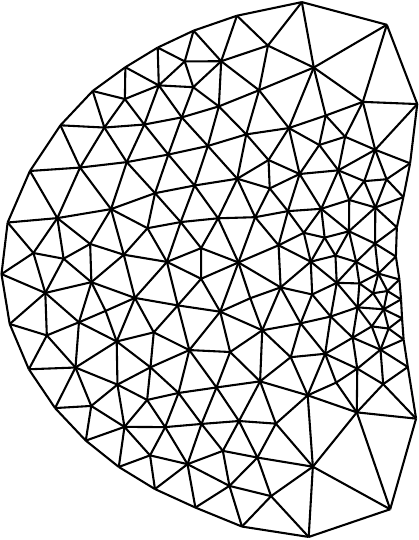}
				\quad
				\includegraphics[scale=0.32]{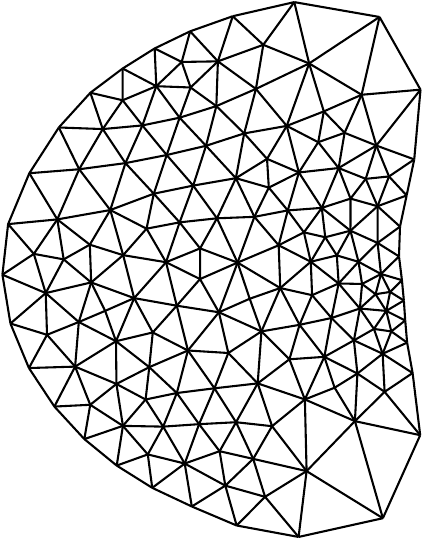}
				\quad
				\includegraphics[scale=0.32]{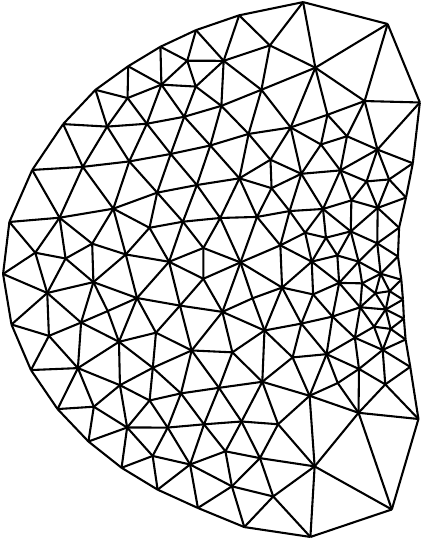}
				\caption{Parameter set $2$.}
			\end{center}
		\end{subfigure}
		\\
		\begin{subfigure}{0.75\textwidth}
			\begin{center}
				\includegraphics[scale=0.32]{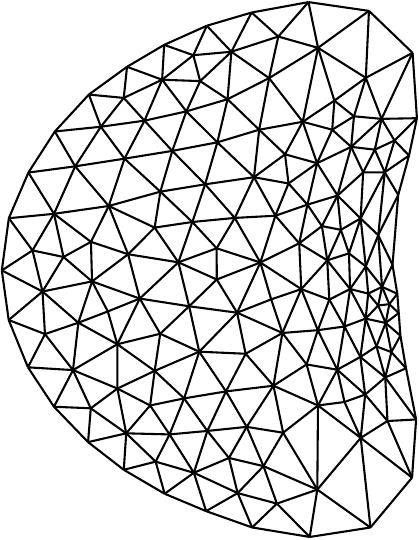}
				\quad
				\includegraphics[scale=0.32]{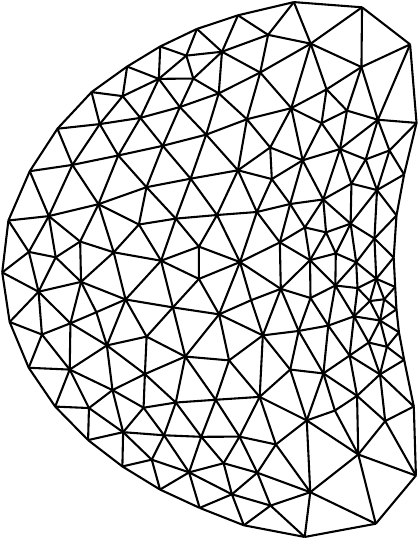}
				\quad
				\includegraphics[scale=0.32]{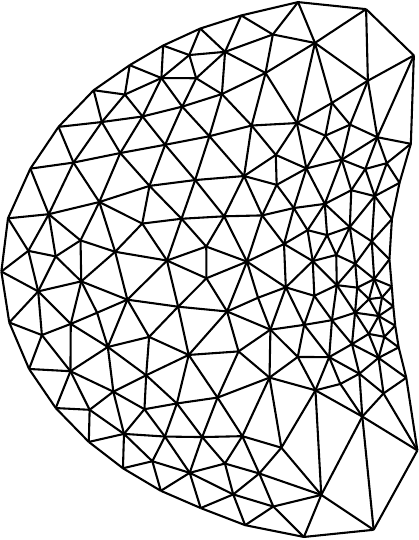}
				\caption{Parameter set $3$.}
			\end{center}
		\end{subfigure}
	\end{center}
	\caption{Final iterates obtained for the penalized problem with variants \EucEuc (left), \ElasEuc (middle) and \CompEuc (right) as described in \cref{subsection:experiment_diff_penalization_parameters}.}
	\label{fig:experiment_diff_penalization_parameters}
\end{figure}
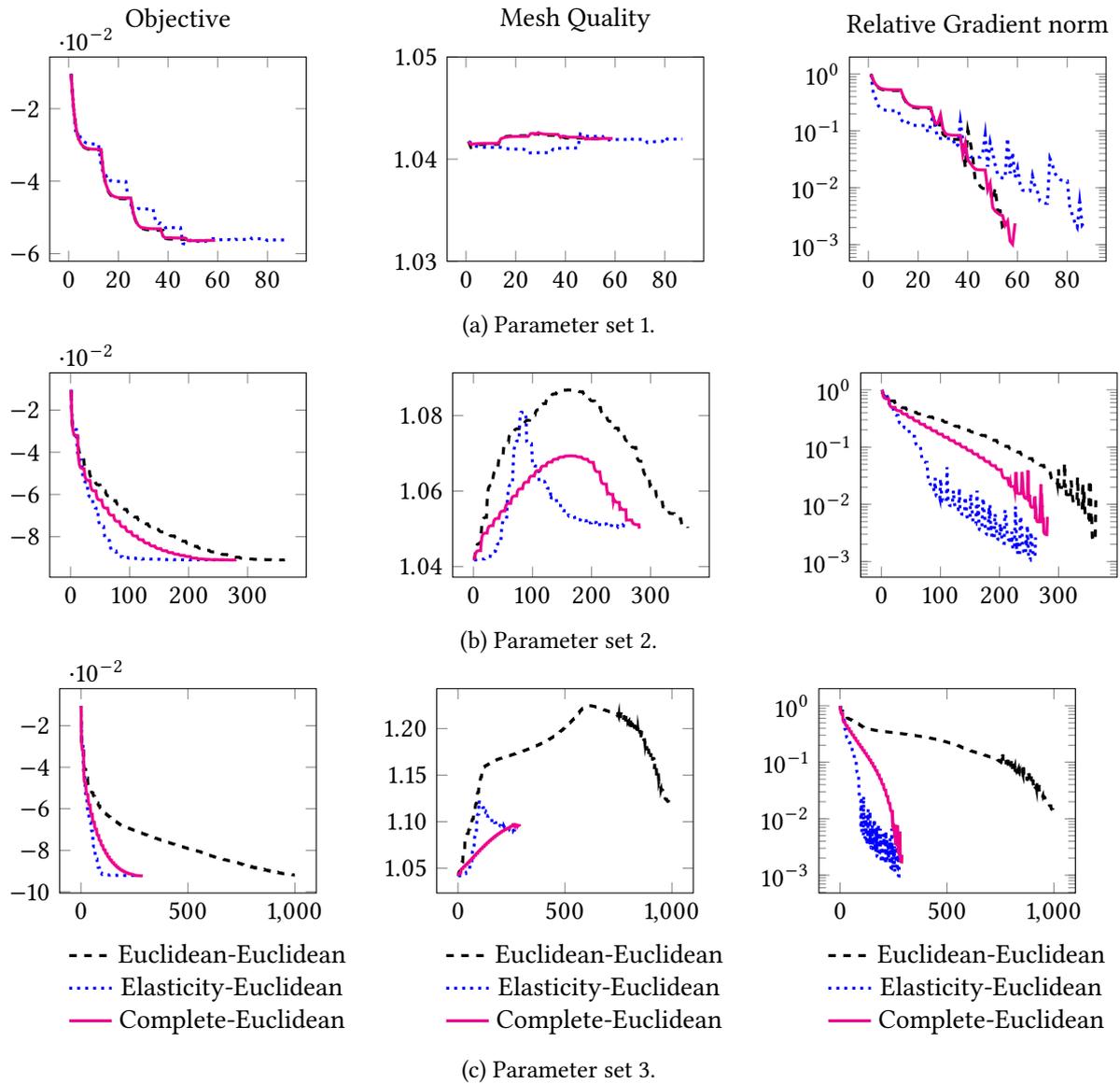
\begin{figure}[htp]
	\begin{subfigure}{\textwidth}
		\begin{center}
			\pgfplotsset{compat=1.14, width=0.33\textwidth}
			\begin{tikzpicture}
				\begin{axis}[
					legend style={draw=none, at={(0,-0.2)},anchor=north west},
					title = {Objective}
					]
					\addplot[color=black, dashed, very thick] table [x=iter,y=Obj, col sep=comma]{matlab/dataResults/experimentPenalization/parameters1/metricType_1_retractionType_1_DATA.txt};
					\addplot[color=blue, dotted, very thick] table [x=iter,y=Obj, col sep=comma]{matlab/dataResults/experimentPenalization/parameters1/metricType_2_retractionType_1_DATA.txt};
					\addplot[color=magenta, very thick] table [x=iter,y=Obj, col sep=comma]{matlab/dataResults/experimentPenalization/parameters1/metricType_3_retractionType_1_DATA.txt};
				\end{axis}
			\end{tikzpicture}
			\hfill
			\pgfplotsset{compat=1.14, width=0.33\textwidth}
			\begin{tikzpicture}
				\begin{axis}[
					legend style={ draw=none, at={(0,-0.2)},anchor=north west},
					title = {Mesh Quality},
					ymin = 1.03,
					ymax = 1.05,
					ytick ={1.03, 1.04, 1.05}
					]
					\addplot[color=black, dashed, very thick] table [x=iter,y=mshQua, col sep=comma]{matlab/dataResults/experimentPenalization/parameters1/metricType_1_retractionType_1_DATA.txt};
					\addplot[color=blue, dotted, very thick] table [x=iter,y=mshQua, col sep=comma]{matlab/dataResults/experimentPenalization/parameters1/metricType_2_retractionType_1_DATA.txt};
					\addplot[color=magenta, very thick] table [x=iter,y=mshQua, col sep=comma]{matlab/dataResults/experimentPenalization/parameters1/metricType_3_retractionType_1_DATA.txt};
				\end{axis}
			\end{tikzpicture}
	\hfill
	\pgfplotsset{compat=1.14, width=0.33\textwidth}
\pgfplotstableread[col sep = comma]{matlab/dataResults/experimentPenalization/parameters1/metricType_1_retractionType_1_DATA.txt}\datatableOne%
\pgfplotstablegetelem{0}{normGrad}\of{\datatableOne}%
\pgfplotstablecreatecol[%
create col/expr={\thisrow{normGrad}/\pgfplotsretval}%
]{relRedNormGrad}{\datatableOne}%
\pgfplotstablesave[col sep = comma]{\datatableOne}{temptable1.txt}%
\pgfplotstableread[col sep = comma]{matlab/dataResults/experimentPenalization/parameters1/metricType_2_retractionType_1_DATA.txt}\datatableTwo%
\pgfplotstablegetelem{0}{normGrad}\of{\datatableTwo}%
\pgfplotstablecreatecol[%
create col/expr={\thisrow{normGrad}/\pgfplotsretval}%
]{relRedNormGrad}{\datatableTwo}%
\pgfplotstablesave[col sep = comma]{\datatableTwo}{temptable2.txt}%
\pgfplotstableread[col sep = comma]{matlab/dataResults/experimentPenalization/parameters1/metricType_3_retractionType_1_DATA.txt}\datatableThree%
\pgfplotstablegetelem{0}{normGrad}\of{\datatableThree}%
\pgfplotstablecreatecol[%
create col/expr={\thisrow{normGrad}/\pgfplotsretval}%
]{relRedNormGrad}{\datatableThree}%
\pgfplotstablesave[col sep = comma]{\datatableThree}{temptable3.txt}%
	\begin{tikzpicture}
		\begin{semilogyaxis}[
			legend style={ draw=none, at={(0,-0.2)},anchor=north west},
			title = {Relative Gradient norm}
			]
			\addplot[color=black, dashed, very thick] table [x=iter,y=relRedNormGrad, col sep=comma]{temptable1.txt};
			\addplot[color=blue, dotted, very thick] table [x=iter,y=relRedNormGrad, col sep=comma]{temptable2.txt};
			\addplot[color=magenta, very thick] table [x=iter,y=relRedNormGrad, col sep=comma]{temptable3.txt};
		\end{semilogyaxis}
	\end{tikzpicture}
	\caption{Parameter set $1$.}
\end{center}
\end{subfigure}
	\\
	\begin{subfigure}{\textwidth}
		\begin{center}
			\pgfplotsset{compat=1.14, width=0.33\textwidth}
			\begin{tikzpicture}
				\begin{axis}
					\addplot[color=black, dashed, very thick] table [x=iter,y=Obj, col sep=comma]{matlab/dataResults/experimentPenalization/parameters2/metricType_1_retractionType_1_DATA.txt};
					\addplot[color=blue, dotted, very thick] table [x=iter,y=Obj, col sep=comma]{matlab/dataResults/experimentPenalization/parameters2/metricType_2_retractionType_1_DATA.txt};
					\addplot[color=magenta, very thick] table [x=iter,y=Obj, col sep=comma]{matlab/dataResults/experimentPenalization/parameters2/metricType_3_retractionType_1_DATA.txt};
				\end{axis}
			\end{tikzpicture}
			\hfill
			\pgfplotsset{compat=1.14, width=0.33\textwidth}
			\begin{tikzpicture}
				\begin{axis}
					\addplot[color=black, dashed, very thick] table [x=iter,y=mshQua, col sep=comma]{matlab/dataResults/experimentPenalization/parameters2/metricType_1_retractionType_1_DATA.txt};
					\addplot[color=blue, dotted, very thick] table [x=iter,y=mshQua, col sep=comma]{matlab/dataResults/experimentPenalization/parameters2/metricType_2_retractionType_1_DATA.txt};
					\addplot[color=magenta, very thick] table [x=iter,y=mshQua, col sep=comma]{matlab/dataResults/experimentPenalization/parameters2/metricType_3_retractionType_1_DATA.txt};
				\end{axis}
			\end{tikzpicture}
			\hfill
			\pgfplotsset{compat=1.14, width=0.33\textwidth}
			\pgfplotstableread[col sep = comma]{matlab/dataResults/experimentPenalization/parameters2/metricType_1_retractionType_1_DATA.txt}\datatableOne%
			\pgfplotstablegetelem{0}{normGrad}\of{\datatableOne}%
			\pgfplotstablecreatecol[%
			create col/expr={\thisrow{normGrad}/\pgfplotsretval}%
			]{relRedNormGrad}{\datatableOne}%
			\pgfplotstablesave[col sep = comma]{\datatableOne}{temptable1.txt}%
			\pgfplotstableread[col sep = comma]{matlab/dataResults/experimentPenalization/parameters2/metricType_2_retractionType_1_DATA.txt}\datatableTwo%
			\pgfplotstablegetelem{0}{normGrad}\of{\datatableTwo}%
			\pgfplotstablecreatecol[%
			create col/expr={\thisrow{normGrad}/\pgfplotsretval}%
			]{relRedNormGrad}{\datatableTwo}%
			\pgfplotstablesave[col sep = comma]{\datatableTwo}{temptable2.txt}%
			\pgfplotstableread[col sep = comma]{matlab/dataResults/experimentPenalization/parameters2/metricType_3_retractionType_1_DATA.txt}\datatableThree%
			\pgfplotstablegetelem{0}{normGrad}\of{\datatableThree}%
			\pgfplotstablecreatecol[%
			create col/expr={\thisrow{normGrad}/\pgfplotsretval}%
			]{relRedNormGrad}{\datatableThree}%
			\pgfplotstablesave[col sep = comma]{\datatableThree}{temptable3.txt}%
			\begin{tikzpicture}
				\begin{semilogyaxis}
					\addplot[color=black, dashed, very thick] table [x=iter,y=relRedNormGrad, col sep=comma]{temptable1.txt};
					\addplot[color=blue, dotted, very thick] table [x=iter,y=relRedNormGrad, col sep=comma]{temptable2.txt};
					\addplot[color=magenta, very thick] table [x=iter,y=relRedNormGrad, col sep=comma]{temptable3.txt};
				\end{semilogyaxis}
			\end{tikzpicture}
			\caption{Parameter set $2$.}
		\end{center}
	\end{subfigure}
	\\
	\begin{subfigure}{\textwidth}
		\begin{center}
			\pgfplotsset{compat=1.14, width=0.33\textwidth}
			\begin{tikzpicture}
				\begin{axis}[
					legend style={ draw=none, at={(0,-0.2)},anchor=north west},
					scaled y ticks = base 10:2
					]
					\addplot[color=black, dashed, very thick] table [x=iter,y=Obj, col sep=comma]{matlab/dataResults/experimentPenalization/parameters3/metricType_1_retractionType_1_DATA.txt};
					\addlegendentry{\EucEuc}
					\addplot[color=blue, dotted, very thick] table [x=iter,y=Obj, col sep=comma]{matlab/dataResults/experimentPenalization/parameters3/metricType_2_retractionType_1_DATA.txt};
					\addlegendentry{\ElasEuc}
					\addplot[color=magenta, very thick] table [x=iter,y=Obj, col sep=comma]{matlab/dataResults/experimentPenalization/parameters3/metricType_3_retractionType_1_DATA.txt};
					\addlegendentry{\CompEuc}
				\end{axis}
			\end{tikzpicture}
			\hfill
			\pgfplotsset{compat=1.14, width=0.33\textwidth}
			\begin{tikzpicture}
				\begin{axis}[
					legend style={ draw=none, at={(0,-0.2)}, anchor=north west},
					y tick label style={
						/pgf/number format/.cd,
						fixed,
						fixed zerofill,
						precision=2,
						/tikz/.cd
					}
					]
					\addplot[color=black, dashed, very thick] table [x=iter,y=mshQua, col sep=comma]{matlab/dataResults/experimentPenalization/parameters3/metricType_1_retractionType_1_DATA.txt};
					\addlegendentry{\EucEuc}
					\addplot[color=blue, dotted, very thick] table [x=iter,y=mshQua, col sep=comma]{matlab/dataResults/experimentPenalization/parameters3/metricType_2_retractionType_1_DATA.txt};
					\addlegendentry{\ElasEuc}
					\addplot[color=magenta, very thick] table [x=iter,y=mshQua, col sep=comma]{matlab/dataResults/experimentPenalization/parameters3/metricType_3_retractionType_1_DATA.txt};
					\addlegendentry{\CompEuc}
				\end{axis}
			\end{tikzpicture}
			\hfill
			\pgfplotsset{compat=1.14, width=0.33\textwidth}
			\pgfplotstableread[col sep = comma]{matlab/dataResults/experimentPenalization/parameters3/metricType_1_retractionType_1_DATA.txt}\datatableOne%
			\pgfplotstablegetelem{0}{normGrad}\of{\datatableOne}%
			\pgfplotstablecreatecol[%
			create col/expr={\thisrow{normGrad}/\pgfplotsretval}%
			]{relRedNormGrad}{\datatableOne}%
			\pgfplotstablesave[col sep = comma]{\datatableOne}{temptable1.txt}%
			\pgfplotstableread[col sep = comma]{matlab/dataResults/experimentPenalization/parameters3/metricType_2_retractionType_1_DATA.txt}\datatableTwo%
			\pgfplotstablegetelem{0}{normGrad}\of{\datatableTwo}%
			\pgfplotstablecreatecol[%
			create col/expr={\thisrow{normGrad}/\pgfplotsretval}%
			]{relRedNormGrad}{\datatableTwo}%
			\pgfplotstablesave[col sep = comma]{\datatableTwo}{temptable2.txt}%
			\pgfplotstableread[col sep = comma]{matlab/dataResults/experimentPenalization/parameters3/metricType_3_retractionType_1_DATA.txt}\datatableThree%
			\pgfplotstablegetelem{0}{normGrad}\of{\datatableThree}%
			\pgfplotstablecreatecol[%
			create col/expr={\thisrow{normGrad}/\pgfplotsretval}%
			]{relRedNormGrad}{\datatableThree}%
			\pgfplotstablesave[col sep = comma]{\datatableThree}{temptable3.txt}%
			\begin{tikzpicture}
				\begin{semilogyaxis}[
					legend style={ draw=none, at={(0,-0.2)}, anchor=north west}
					]
					\addplot[color=black, dashed, very thick] table [x=iter,y=relRedNormGrad, col sep=comma]{temptable1.txt};
					\addlegendentry{\EucEuc}
					\addplot[color=blue, dotted, very thick] table [x=iter,y=relRedNormGrad, col sep=comma]{temptable2.txt};
					\addlegendentry{\ElasEuc}
					\addplot[color=magenta, very thick] table [x=iter,y=relRedNormGrad, col sep=comma]{temptable3.txt};
					\addlegendentry{\CompEuc}
				\end{semilogyaxis}
			\end{tikzpicture}
			\caption{Parameter set $3$.}
		\end{center}
	\end{subfigure}
	\caption{Objective and mesh quality for the penalized problem described in \cref{subsection:experiment_diff_penalization_parameters}.}
	\label{fig:experimentPenalization}
\end{figure}

\begin{table}[htp]
	\begin{center}
		\begin{tabular}{c}
			\toprule
			Variant  \\
			\midrule
			\EucEuc   \\
			\ElasEuc  \\
			\CompEuc \\
			\bottomrule
		\end{tabular}%
		\sisetup{round-mode = places, round-precision = 3, scientific-notation = fixed, fixed-exponent = 0}%
		\begin{filecontents*}{matlab/timesOfExecutionExperimentTwo.txt}
iter,total,State,dObj,dPen,backt,grad
1000,10.159,4.262,1.2625,0.32715,3.4577,0
279,3.6176,1.3851,0.37143,0.091637,1.1161,0.43412
289,3.316,1.3358,0.3827,0.17545,1.0331,0.081148
\end{filecontents*}
\pgfplotstabletypeset[%
		col sep = comma,
		every head row/.style={before row=\toprule,after row=\midrule},
		every last row/.style={after row=\bottomrule},
		columns={iter,total,perIter,State,dObj, dPen, backt,grad},
		columns/iter/.style={int detect, column name={iter}},
		columns/total/.style={std, precision=2,column name={total},
			assign cell content/.code={%
				\pgfkeyssetvalue{/pgfplots/table/@cell content}%
				{\SI{##1}{\second}}%
		}},
		columns/perIter/.style={fixed, column name={per iter},
			assign cell content/.code={%
				\pgfkeyssetvalue{/pgfplots/table/@cell content}%
				{\SI{##1}{\second}}%
		}},
		columns/State/.style={std, precision=2,column name={state},
			assign cell content/.code={%
				\pgfkeyssetvalue{/pgfplots/table/@cell content}%
				{\SI{##1}{\second}}%
		}},
		columns/dObj/.style={std, precision=2,column name={dObj},
			assign cell content/.code={%
				\pgfkeyssetvalue{/pgfplots/table/@cell content}%
				{\SI{##1}{\second}}%
		}},
		columns/dPen/.style={std, precision=2,column name={dPen},
			assign cell content/.code={%
				\pgfkeyssetvalue{/pgfplots/table/@cell content}%
				{\SI{##1}{\second}}%
		}},
		columns/backt/.style={std, precision=2,column name={backt},
			assign cell content/.code={%
				\pgfkeyssetvalue{/pgfplots/table/@cell content}%
				{\SI{##1}{\second}}%
		}},
		columns/grad/.style={std, precision=2,column name={grad},
			assign cell content/.code={%
				\pgfkeyssetvalue{/pgfplots/table/@cell content}%
				{\SI{##1}{\second}}%
		}}%
		]{matlab/timesOfExecutionExperimentTwo.txt}
	\end{center}
	\caption{Execution times for the experiment described in \cref{subsection:experiment_diff_penalization_parameters} with $\alpha_1 = 0.015$, $\alpha_2 = 0.005$, $\alpha_3 = 0$, and $\alpha_4 = 0.0005$.}
	\label{table:executionTimesvariant3}
\end{table}

\subsection{Experiment~3: Unpenalized Problem with Finer Meshes}
\label{subsection:unpenalizedDiffMeshLevels}

The penalty approach may be criticized since it requires the user to make a somewhat arbitrary choice of the penalty parameters $\alphaPenalization_j$, $j = 1, \ldots, 4$.
Therefore we revisit here the unpenalized problem, aware of the fact that the discretized problem does not possess a solution any gradient method could converge to.
In contrast to the results of \cref{subsection:unpenalizedAndIdeal}, the meshes are now finer, and we only compare the two most promising gradient descent variants, \ElasEuc and \CompEuc.
We consider four mesh levels.
The first one contains $N_V = 541$~vertices and $N_T = 1016$~triangles.
The second one has $N_V = 775$~vertices and $N_T = 1468$~elements.
The third possesses $N_V = 2191$~vertices and $N_T = 4252$~triangles.
Finally, mesh level four has $N_V = 13455$~vertices and $N_T = 26588$~triangles.

\begin{figure}[htp]
	\begin{center}
		\pgfplotsset{compat=1.14, width=0.45\textwidth}
		\begin{tikzpicture}
			\begin{axis}[
				title = Objective,
				legend style={draw=none, at={(0,-0.2)},anchor=north west}
				]
				\addplot[color=blue, dotted, very thick] table [x=iter,y=Obj, col sep=comma]{matlab/dataResults/experimentNoPenalization/meshLevel2/metricType_2_retractionType_1_DATA.txt};
				\addlegendentry{\ElasEuc}
				\addplot[color=magenta, very thick] table [x=iter,y=Obj, col sep=comma]{matlab/dataResults/experimentNoPenalization/meshLevel2/metricType_3_retractionType_1_DATA.txt};
				\addlegendentry{\CompEuc}
			\end{axis}
		\end{tikzpicture}
		\hfill
		\pgfplotsset{compat=1.14, width=0.45\textwidth}
		\begin{tikzpicture}
			\begin{axis}[
				title = Mesh Quality,
				legend style={draw=none, at={(0,-0.2)},anchor=north west}
				]
				\addplot[color=blue, dotted, thick] table [x=iter,y=mshQua, col sep=comma]{matlab/dataResults/experimentNoPenalization/meshLevel2/metricType_2_retractionType_1_DATA.txt};
				\addlegendentry{\ElasEuc}
				\addplot[color=magenta, very thick] table [x=iter,y=mshQua, col sep=comma]{matlab/dataResults/experimentNoPenalization/meshLevel2/metricType_3_retractionType_1_DATA.txt};
				\addlegendentry{\CompEuc}
			\end{axis}
		\end{tikzpicture}
		\caption{Objective and mesh quality for the unpenalized problem at mesh level~$2$ described in \cref{subsection:unpenalizedDiffMeshLevels}.}
		\label{fig:ObjMeshQuaExperiment3}
	\end{center}
\end{figure}
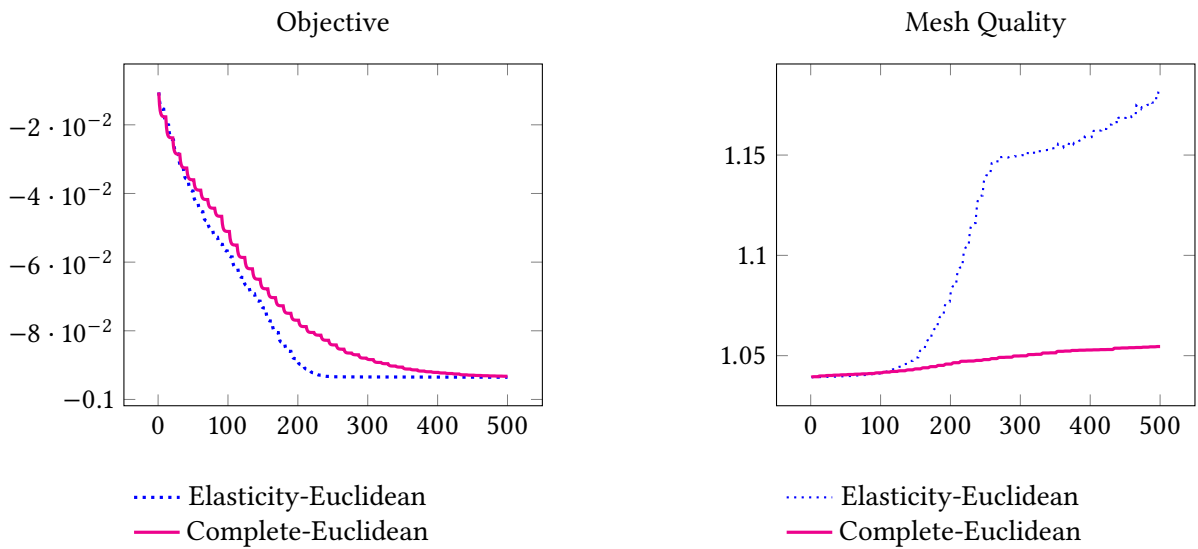

We allow the algorithm to run $500$~iterations, and we are mainly interested in comparing the values of the objective function and the mesh quality.
The results can be seen in \cref{fig:diffMeshLevels,fig:ObjMeshQuaExperiment3}.
We infer that both variants, \ElasEuc and \CompEuc, achieve a similar decrease of the objective function.
\ElasEuc needs fewer iterations to reach the plateau, but \CompEuc maintains a better mesh quality measure and has less numerical cost per iteration.
The latter is reflected in \cref{table:timesExecutionMeshLevels}.
Here we separately display the time required to \enquote{assemble} the matrices representing the Riemannian metric in column \texttt{assemG}.
More precisely, as in all experiments before, we only actually form this matrix in case of \ElasEuc, and employ a sparse direct solver to obtain the solution of the gradient equation \eqref{eq:definition_riemannian_gradient}.
In case of \CompEuc, we continue to work with matrix-vector products and the conjugate gradient solver.
In this case, the column \texttt{assemG} is dominated by the time to evaluate the first-order derivative of the penalty function.
We also observe that the time required to solve the gradient equation \eqref{eq:definition_riemannian_gradient} remains essentially constant in case of \CompEuc while the time for the direct solver in case of \ElasEuc grows with the problem size.

An inspection of the meshes at iteration~$500$ in \cref{fig:diffMeshLevels} shows triangles closer to equilateral when using \CompEuc and more elongated in case of \ElasEuc, as reflected by mesh quality plot in \cref{fig:ObjMeshQuaExperiment3}.

\begin{remark}
	\label{remark:relationMeshQualityMeasures}
	We recall that the mesh quality measure shown in \cref{fig:ObjMeshQuaExperiment3} is given in \eqref{eq:quality_measure}.
		Other quality measures, such as the aspect or radius ratios (see for instance \cite[Table~6, Rows~7 and~9]{Shewchuk:2002:1}), exhibit similar results.
	The complete metric \eqref{eq:phiCompleteMetric} with augmentation function \eqref{eq:penalization} gave rise to better values of all mesh quality measures considered, compared to the elasticity metric.
\end{remark}

Moreover, the triangles are smaller and the vertices more dense in regions which have deformed most compared to the initial circle mesh.
We can consider this behavior as a natural redistribution of the nodes promoted by the use of the complete metric.

\begin{figure}[htp]
	\captionsetup[subfigure]{labelformat=empty}
	\begin{subfigure}{0.45\textwidth}
		\centering
		\begin{tikzpicture}[zoomboxarray, zoomboxarray columns=1, zoomboxarray rows=1, zoomboxes below]
			\node [image node]{\includegraphics[width=0.4\textwidth]{Figures/experimentNoPenalization/meshLevel1/iteration_0500_mesh_metric_type_2_retraction_type_1.pdf} };
			\zoombox[magnification = 3,color code=blue]{0.85,0.5}
		\end{tikzpicture}
		\quad
		\begin{tikzpicture}[zoomboxarray, zoomboxarray columns=1, zoomboxarray rows=1, zoomboxes below]
			\node [image node]{\includegraphics[width=0.4\textwidth]{Figures/experimentNoPenalization/meshLevel1/iteration_0500_mesh_metric_type_3_retraction_type_1.pdf} };
			\zoombox[magnification = 3, color code=magenta]{0.85,0.5}
		\end{tikzpicture}
		\caption{Mesh Level $1$.}
		\label{fig:ElasEuc}
	\end{subfigure}
	\hfill
	\begin{subfigure}{0.45\textwidth}
		\centering
		\begin{tikzpicture}[zoomboxarray, zoomboxarray columns=1, zoomboxarray rows=1, zoomboxes below]
			\node [image node]{\includegraphics[width=0.4\textwidth]{Figures/experimentNoPenalization/meshLevel2/iteration_0500_mesh_metric_type_2_retraction_type_1.pdf} };
			\zoombox[magnification = 3,color code=blue]{0.85,0.5}
		\end{tikzpicture}
		\quad
		\begin{tikzpicture}[zoomboxarray, zoomboxarray columns=1, zoomboxarray rows=1, zoomboxes below]
			\node [image node]{\includegraphics[width=0.4\textwidth]{Figures/experimentNoPenalization/meshLevel2/iteration_0500_mesh_metric_type_3_retraction_type_1.pdf} };
			\zoombox[magnification = 3, color code=magenta]{0.85,0.5}
		\end{tikzpicture}
		\caption{Mesh Level $2$.}
		\label{fig:MeshLevel12}
	\end{subfigure}
	\caption{$500$th iterate in case of \ElasEuc (blue) and \CompEuc (magenta), for the experiment described in~\cref{subsection:unpenalizedDiffMeshLevels}.}
	\label{fig:diffMeshLevels}
\end{figure}

\begin{table}[htp]
\begin{center}
\begin{tabular}{cc}
	\toprule
Mesh Level         & Variant \\
\midrule
\multirow{2}{*}{1} &   \ElasEuc \\
                   &   \CompEuc      \\
									 \midrule
\multirow{2}{*}{2} &    \ElasEuc     \\
                   &    \CompEuc     \\
									 \midrule
\multirow{2}{*}{3} &   \ElasEuc      \\
                   &    \CompEuc     \\
									 \midrule
\multirow{2}{*}{4} &  \ElasEuc       \\
                   &     \CompEuc  \\
									  	\bottomrule
\end{tabular}%
	\sisetup{round-mode = places, round-precision = 3, scientific-notation = fixed, fixed-exponent = 0}%
	\begin{filecontents*}{matlab/timesOfExecutionExperimentThree.txt}
iter,total,dObj,backt,assemblyG,grad
500,10.979,1.3309,3.091,1.3945,1.0447
500,8.767,1.335,2.756,0.38091,0.11549
500,13.772,1.8188,3.9246,1.5558,1.3476
500,10.81,1.7645,3.3539,0.543,0.11747
500,35.869,5.4867,9.0416,4.2229,4.5789
500,27.605,5.3873,8.2337,1.416,0.21397
500,212.24,32.856,51.845,24.97,37.816
500,150.2,32.329,45.051,8.736,0.36437
\end{filecontents*}
\pgfplotstabletypeset[%
	col sep = comma,
	every head row/.style={before row=\toprule,after row=\midrule},
	every last row/.style={after row=\bottomrule},
	every nth row={2}{before row=\midrule},
	columns={iter,total,perIter,dObj, backt, assemblyG, grad},
	columns/iter/.style={int detect, column name={iter}},
	columns/total/.style={std, precision=2,column name={total},
	assign cell content/.code={%
	\pgfkeyssetvalue{/pgfplots/table/@cell content}%
	{\SI{##1}{\second}}%
	}},
	columns/total/.style={std, precision=2,column name={total},
	assign cell content/.code={%
	\pgfkeyssetvalue{/pgfplots/table/@cell content}%
	{\SI{##1}{\second}}%
	}},
	columns/perIter/.style={fixed, column name={per iter},
	assign cell content/.code={%
	\pgfkeyssetvalue{/pgfplots/table/@cell content}%
	{\SI{##1}{\second}}%
	}},
	columns/dObj/.style={std, precision=2,column name={dObj},
	assign cell content/.code={%
	\pgfkeyssetvalue{/pgfplots/table/@cell content}%
	{\SI{##1}{\second}}%
	}},
	columns/backt/.style={std, precision=2,column name={backt},
	assign cell content/.code={%
	\pgfkeyssetvalue{/pgfplots/table/@cell content}%
	{\SI{##1}{\second}}%
	}},
	columns/assemblyG/.style={std, precision=2,column name={assemblyG},
	assign cell content/.code={%
	\pgfkeyssetvalue{/pgfplots/table/@cell content}%
	{\SI{##1}{\second}}%
	}},
	columns/grad/.style={std, precision=2,column name={grad},
	assign cell content/.code={%
	\pgfkeyssetvalue{/pgfplots/table/@cell content}%
	{\SI{##1}{\second}}%
	}}%
	]{matlab/timesOfExecutionExperimentThree.txt}
	\end{center}
	\caption{Execution times for $500$~iterations for the experiments described in \cref{subsection:unpenalizedDiffMeshLevels}.}
	\label{table:timesExecutionMeshLevels}
\end{table}

\subsection{Experiment~4: Two-dimensional Optimal Bridge}
\label{subsection:experimentOptimalBridge}

We end this section by showing the performance of the variants \ElasEuc and \CompEuc while solving the penalized version of a well-known compliance minimization problem from structural mechanics.
We focus our attention on the design of a two-dimensional bridge.
The goal is to minimize the compliance of the system.
The area of the geometry is fixed, and, as suggested in \cite[Remark~1]{AllaireDapognyFrey:2014:1}, the area is added to the objective function using a fixed value of the Lagrange multiplier.
It is worth recalling that the penalized version of this problem has at least one globally optimal solution as per~\cref{corollary:existenceSolutionpenalizedprob}, since the compliance and the area are bounded from below and the objective is lower semi-continuous and not identically to $+\infty$.

Mathematically, this problem can be written as
\begin{equation}
	\label{eq:penalizedOptimalBridgeProblem}
	\begin{aligned}
		\text{Minimize}
		&
		\quad
		\int_{\Omega_Q} \bf \cdot \by \d x
		+
		\int_{\Gamma^N_Q} \bg \cdot \by \d s
		+
		\ell \, \int_{\Omega_Q} 1 \d x
		+
		\newpenalty(Q;\Qref)
		\\
		\text{\wrt}
		&
		\quad
		Q \in \planarmanifold
		,
		\;
		\by \in S^1_{\Gamma^D}(\Omega_Q) \times S^1_{\Gamma^D}(\Omega_Q)
		\\
		\text{\st}
		&
		\quad
		2 \muState \int_{\Omega_Q} \bvarepsilon(\by) \dprod \bvarepsilon(\bv) \d x
		+
		\lambdaState \int_{\Omega_Q} \trace(\bvarepsilon(\by)) \, \trace(\bvarepsilon(\bv)) \d x
		=
		\int_{\Omega_Q} \bf \cdot \bv \, \d x
		+
		\int_{\Gamma^N_Q} \bg \cdot \bv \, \d s
		\\
		\text{for all }
		&
		\quad
		\bv \in S^1_{\Gamma^D}(\Omega_Q) \times S^1_{\Gamma^D}(\Omega_Q),
	\end{aligned}
\end{equation}
where $\bf \in L^2(\R^2) \times L^2(\R^2)$ and $\bg \in H^1(\R^2) \times H^1(\R^2)$ are given volume and boundary loads, respectively.
Moreover, $S^1_{\Gamma^D}(\Omega_Q)$ denotes the space of piecewise linear, globally continuous finite element functions defined on the mesh $\Omega_Q$, with zero Dirichlet boundary condition on $\Gamma^D$.

We choose Young's modulus $\EState = 1$ and Poisson ratio $\nuState = 0.3$.
The associated Lamé parameters $\lambdaState$ and $\muState$ are computed using~\eqref{eq:Lame_parameters}.

\begin{figure}
	\begin{tikzpicture}
		\draw[thick] (0,0) -- (0,1) -- (2.5,4) -- (5,5) -- (7.5,4) -- (10,1) -- (10,0) -- (9,0) --  (5.5,0) -- (4.5,0) -- (1,0) --cycle;
		\draw[line width = 4pt, color = DirichletColor] (0,0)--(1,0);
		\draw[line width = 4pt, color = DirichletColor] (9,0)--(10,0);
		\draw[->, thick] (4.5,0) -- (4.5,-0.5);
		\draw[->, thick] (5,0) -- (5,-0.5);
		\draw[->, thick] (5.5,0) -- (5.5,-0.5);
		\draw[line width = 4pt, color = NeumannColor] (4.5,0)--(5.5,0);
		\node at (5,-0.7)  {$\bg = [0,-0.25]^\transp$};
		\draw[thick] (2,1) circle (0.5);
		\draw[thick] (3,3) circle (0.5);
		\draw[thick] (7,3) circle (0.5);
		\draw[thick] (8,1) circle (0.5);
		\node[left] at (0,0)  {$(0,0)$};
		\node[left] at (0,1)  {$(0,1)$};
		\node[left] at (2.5,4)  {$(2.5,4)$};
		\node[above] at (5,5)  {$(5,5)$};
		\node[right] at (7.5,4)  {$(7.5,4)$};
		\node[right] at (10,1)  {$(10,1)$};
		\node[right] at (10,0)  {$(10,0)$};
		\node[below] at (0.5,0)  {$\Gamma^D$};
		\node[below] at (9.5,0)  {$\Gamma^D$};
		\node[above] at (5,0)  {$\Gamma^N$};
	\end{tikzpicture}
	\caption{Informed initial shape for the experiment described in~\cref{subsection:experimentOptimalBridge}. The portion of the boundary depicted in blue corresponds to the homogeneous Dirichlet condition~$\Gamma^D$, and the green one to the inhomogeneous Neumann condition~$\Gamma^N$.}
	\label{fig:initialOptimalBridge}
\end{figure}
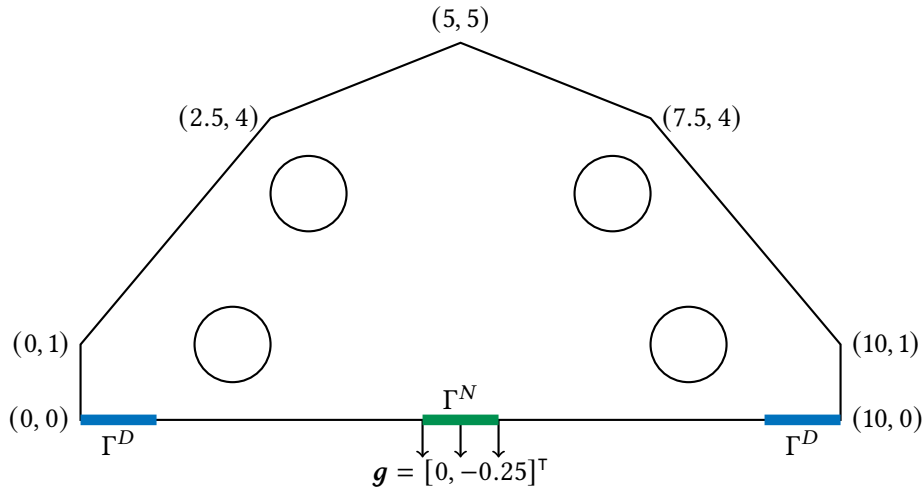

We assume that no body forces are applied, \ie $\bf = [0,0]^\transp$.
The boundary loads on the inhomogeneous Neumann portion of the boundary are given by $\bg = [0,-0.25]^\transp$.
When treating the problem as a topology problem, as was done in \cite[Section~6.2.1]{AllaireDapognyFrey:2014:1}, optimized shapes contain a number of holes.
We therefore use an informed initial shape depicted in~\cref{fig:initialOptimalBridge}.
In our experiment, we only allow the holes to be deformed and we fix the outer boundary by imposing the appropriate zero boundary conditions on the shape gradient \eqref{eq:definition_riemannian_gradient}.

The computation of the shape derivative is performed in a fully discrete framework, following the ideas presented in~\cite{Berggren:2010:1} and the Lagrangian formalism.
Due to the use of the compliance functional, no adjoint state is needed.
The fully discrete shape derivative for this problem in its volume formulation is obtained by a straightforward calculation and it is given by
\begin{equation}
	\label{eq:shape_derivative_optimal_bridge}
	\begin{split}
		\d_Q \cL[Q, \by][\bV]
		=
		2 \int_{\Omega_Q} (D\bf \, \by) \cdot \bV + (\bf \cdot \by) \div(\bV) \, \d x
		- \int_{\Omega_Q} A \bvarepsilon (\by) \dprod  \bvarepsilon (\by) \div(\bV) \, \d x \\
		+ \int_{\Omega_Q} A \, \textup{sym}(D\by \, D\bV) \dprod \bvarepsilon (\by) \, \d x
		+ \int_{\Omega_Q} A \, \bvarepsilon(\by) \dprod \textup{sym}(D \by \, D\bV) \, \d x \\
		+ \int_{\Gamma_N} (\bg \cdot \by) (\div(\bV)-n^\transp \textup{sym}(D\bV) \bn) \, \d s
		+ \ell \int_{\Omega_Q} \div(\bV) \, \d x
		+ \frac{\partial \varphi (Q)}{\partial (\vvec Q)^i}
		,
	\end{split}
\end{equation}
Here $\cL$ denotes the Lagrangian associated with problem \eqref{eq:penalizedOptimalBridgeProblem}.
Moreover, $A$ maps $\R^{2\times 2}$ into $\R^{2\times 2}$ such that $A G = 2 \muState G + \lambdaState \trace(G) \, \id $ holds, and $\textup{sym}(\cdot)$ denotes the symmetrization of a matrix.
The deformation field $\bV$ is any piecewise linear $\R^2$-valued function, described by its nodal values.

The initial mesh consists of $467$~nodes and $816$~elements.
The penalty parameters in \eqref{eq:penalization} are chosen as $\alphaPenalization_1 = 10^{-3}$ and $\alphaPenalization_2 = \alphaPenalization_3 = \alphaPenalization_4 = 0$.
Notice that we can set the parameters $\alphaPenalization_2,\alphaPenalization_3,\alphaPenalization_4$ to zero without affecting the existence of solution of the problem.
This is due to the fact that the outer boundary is fixed and no exterior self-intersections are possible.
Moreover, the mesh cannot shrink to a point, and translations of the mesh infinitely far away from the reference mesh $\Qref$ are also impossible.
The values for the elasticity metric (not be confused with the parameters associated to the linear elasticity model of the state equation) are $\EMetric = 1$, $\nuMetric = 0.4$ and $\deltaMetric = 0.2$.
On the other hand, the parameters for the complete metric given in~\eqref{eq:complete_metric_for_planar_triangular_meshes} are $\alphaMetric_1 = 13.25$ and $\alphaMetric_2 =  \alphaMetric_3 = \alphaMetric_4 = 0$.

\begin{figure}
	\begin{center}
	\begin{subfigure}{0.32\textwidth}
		\centering
		\includegraphics[width = 0.95\textwidth]{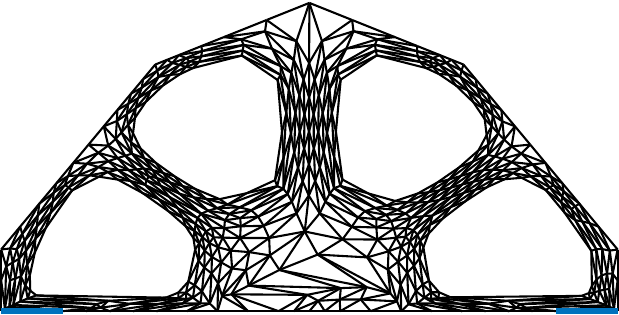}
		\caption{Iteration $1000$}
	\end{subfigure}
	\hfill
	\begin{subfigure}{0.32\textwidth}
		\centering
		\includegraphics[width = 0.95\textwidth]{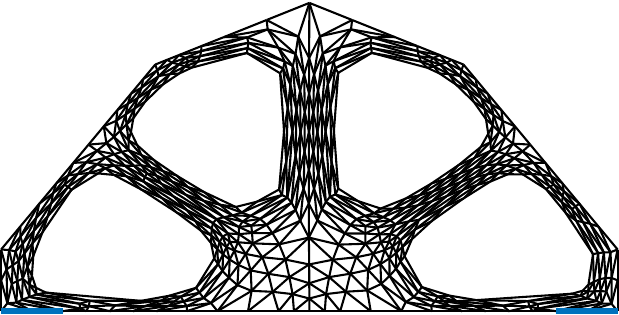}
		\caption{Iteration $700$}
	\end{subfigure}
	\hfill
	\begin{subfigure}{0.32\textwidth}
		\centering
		\includegraphics[width = 0.95\textwidth]{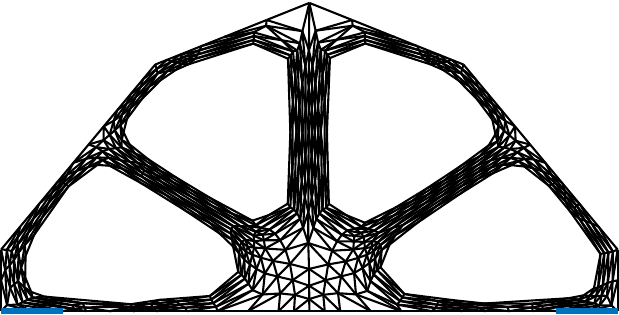}
		\caption{Iteration $1000$}
	\end{subfigure}
	\\
	\begin{subfigure}{0.32\textwidth}
		\centering
		\includegraphics[width = 0.95\textwidth]{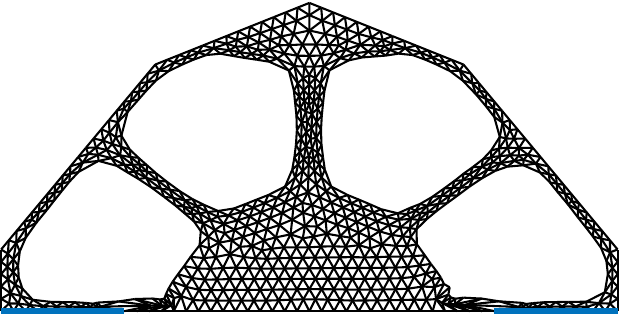}
		\caption{Iteration $200$}
	\end{subfigure}
	\hfill
	\begin{subfigure}{0.32\textwidth}
		\centering
		\includegraphics[width = 0.95\textwidth]{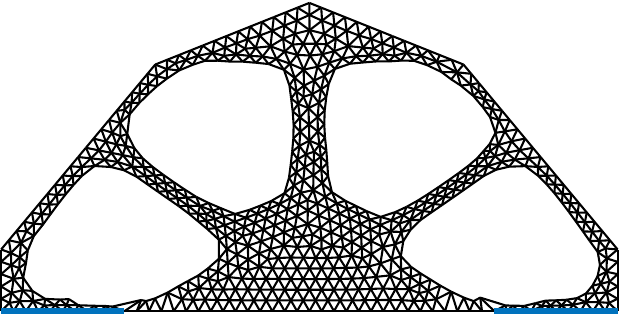}
			\caption{Iteration $50$}
	\end{subfigure}
	\hfill
	\begin{subfigure}{0.32\textwidth}
		\centering
		\includegraphics[width = 0.95\textwidth]{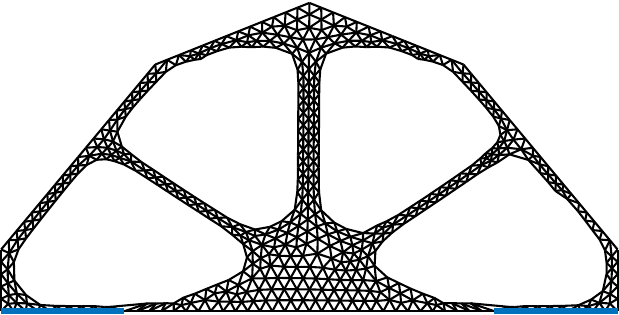}
			\caption{Iteration $100$}
	\end{subfigure}
\end{center}
	\caption{Iterates for the variant \CompEuc before (upper row) and after (lower row) remeshing from the optimal bridge experiment described in~\cref{subsection:experimentOptimalBridge}. The associated Lagrange multipliers are $\ell = 9,9 \cdot 10^{-2}$ (left), $\ell = 0.5$ (middle) and $\ell=1.5$ (right).}
	\label{fig:iterations_before_after_remeshing}
\end{figure}

We solve problem \eqref{eq:penalizedOptimalBridgeProblem} with three different values of the Lagrange multiplier $\ell$.
In each case, we let the gradient descent method run for a fixed number of iterations.
We use the values $\ell = 9.9 \cdot 10^{-2}$, $0.5$, and $1.5$, and we let the algorithm iterate for $1000$, $700$ and $1000$ steps, respectively.
In the upper row of~\cref{fig:iterations_before_after_remeshing}, we show the resulting meshes.
It is clear that due to the nature of the problem, no further improvement can be achieved by keeping the same mesh connectivity fixed to that of the initial mesh.
Thus, we remesh the shape using the current boundary.
After remeshing, we let the algorithm iterate, respectively, for $200$, $50$ and $100$ additional steps.
The final meshes are shown in the lower row of~\cref{fig:iterations_before_after_remeshing}.
It is worth mentioning that even though the addition of the penalization term to the objective function has an impact on the optimal shape, choosing appropriate values of the penalization parameters allows us to recover the expected (thin) structures which are typical in this context.

\begin{figure}
	\begin{center}
	\begin{subfigure}{0.32\textwidth}
		\centering
		\includegraphics[width = 0.95\textwidth]{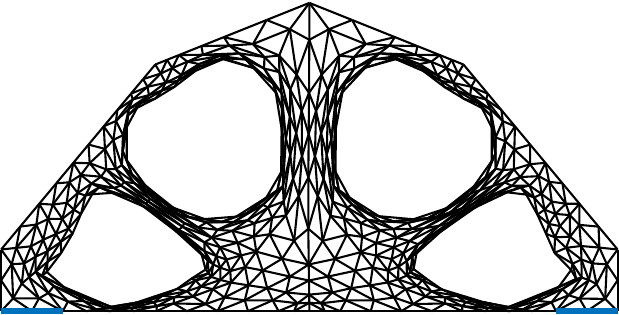}
		\caption{Iteration $1000$}
	\end{subfigure}
	\hfill
	\begin{subfigure}{0.32\textwidth}
		\centering
		\includegraphics[width = 0.95\textwidth]{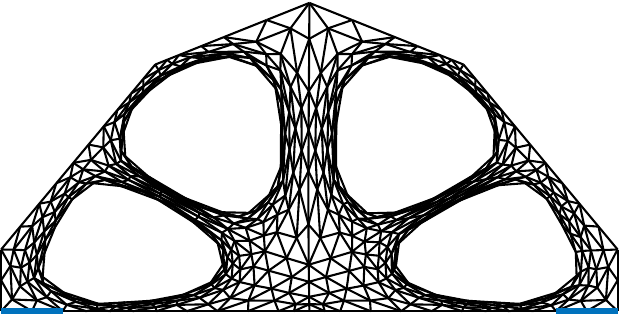}
		\caption{Iteration $700$}
	\end{subfigure}
	\hfill
	\begin{subfigure}{0.32\textwidth}
		\centering
		\includegraphics[width = 0.95\textwidth]{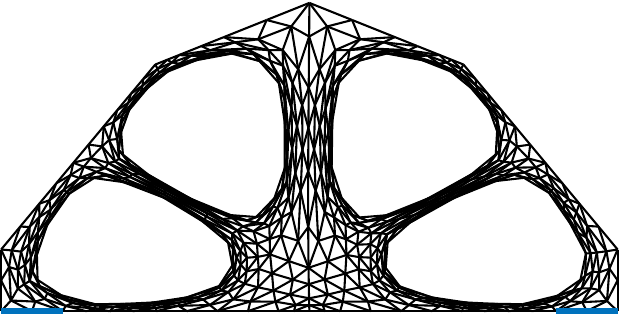}
		\caption{Iteration $1000$}
	\end{subfigure}
	\\
	\begin{subfigure}{0.32\textwidth}
		\centering
		\includegraphics[width = 0.95\textwidth]{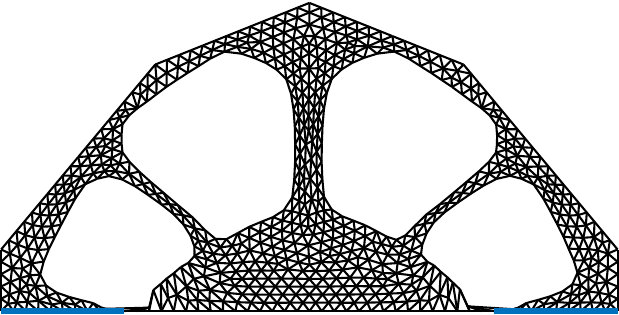}
		\caption{Iteration $200$}
	\end{subfigure}
	\hfill
	\begin{subfigure}{0.32\textwidth}
		\centering
		\includegraphics[width = 0.95\textwidth]{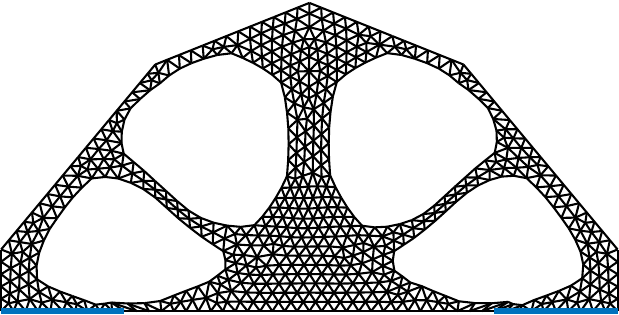}
			\caption{Iteration $50$}
	\end{subfigure}
	\hfill
	\begin{subfigure}{0.32\textwidth}
		\centering
		\includegraphics[width = 0.95\textwidth]{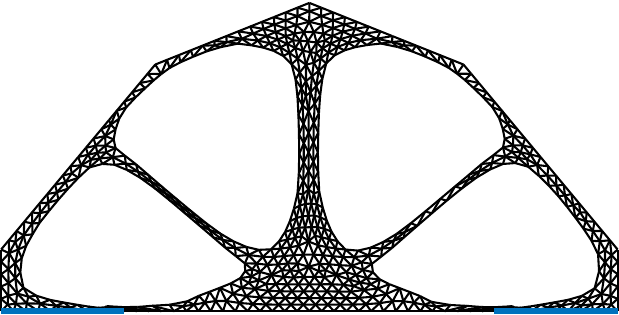}
			\caption{Iteration $100$}
	\end{subfigure}
\end{center}
	\caption{Iterates for the variant \ElasEuc before (upper row) and after (lower row) remeshing from the optimal bridge experiment described in~\cref{subsection:experimentOptimalBridge}. The associated Lagrange multipliers are $\ell = 9,9 \cdot 10^{-2}$ (left), $\ell = 0.5$ (middle) and $\ell=1.5$ (right).}
	\label{fig:iterations_before_after_remeshing_variant2}
\end{figure}

We display in \cref{table:timesExecutionOptimalBridge} the mesh sizes as well as the number of iterations of each experiment, together with the times of execution.
We can conclude that variant \CompEuc outperforms variant \ElasEuc for all experiments.
For comparison, we display the meshes obtained with the latter in \cref{fig:iterations_before_after_remeshing_variant2}.

\begin{figure}[htp]
	\begin{center}
		\pgfplotsset{compat=1.14, width=0.32\textwidth}
		\begin{tikzpicture}
			\begin{axis}[
				title = {Objective with $\ell = 9.9 \cdot 10 ^{-2}$},
				legend style={draw=none, at={(0,-0.2)},anchor=north west}
				]
				\addplot[color=blue, dotted, very thick] table [x=iter,y=Obj, col sep=comma]{matlab/dataResults/optimalBridge_new/ell99Em3_Initial/metricType_2_retractionType_1_DATA.txt};
				\addplot[color=blue, dotted, very thick] table [x expr=\thisrow{iter}+1000,y=Obj, col sep=comma]{matlab/dataResults/optimalBridge_new/ell99Em3_WarmStart_variant2/metricType_2_retractionType_1_DATA.txt};
				\addplot[color=magenta, very thick] table [x=iter,y=Obj, col sep=comma]{matlab/dataResults/optimalBridge_new/ell99Em3_Initial/metricType_3_retractionType_1_DATA.txt};
				\addplot[color=magenta, very thick] table [x expr=\thisrow{iter}+1000,y=Obj, col sep=comma]{matlab/dataResults/optimalBridge_new/ell99Em3_WarmStart_variant3/metricType_3_retractionType_1_DATA.txt};
			\end{axis}
		\end{tikzpicture}
		\hfill
		\pgfplotsset{compat=1.14, width=0.32\textwidth}
		\begin{tikzpicture}
			\begin{axis}[
				title = {Objective with $\ell = 5 \cdot 10 ^{-1}$},
				legend style={draw=none, at={(0,-0.2)},anchor=north west}
				]
				\addplot[color=blue, dotted, very thick] table [x=iter,y=Obj, col sep=comma]{matlab/dataResults/optimalBridge_new/ell5Em1_Initial/metricType_2_retractionType_1_DATA.txt};
				\addplot[color=blue, dotted, very thick] table [x expr=\thisrow{iter}+700,y=Obj, col sep=comma]{matlab/dataResults/optimalBridge_new/ell5Em1_WarmStart_variant2/metricType_2_retractionType_1_DATA.txt};
				\addplot[color=magenta,  very thick] table [x=iter,y=Obj, col sep=comma]{matlab/dataResults/optimalBridge_new/ell5Em1_Initial/metricType_3_retractionType_1_DATA.txt};
				\addplot[color=magenta, very thick] table [x expr=\thisrow{iter}+700,y=Obj, col sep=comma]{matlab/dataResults/optimalBridge_new/ell5Em1_WarmStart_variant3/metricType_3_retractionType_1_DATA.txt};
			\end{axis}
		\end{tikzpicture}
		\hfill
		\pgfplotsset{compat=1.14, width=0.32\textwidth}
		\begin{tikzpicture}
			\begin{axis}[
				title = {Objective with $\ell = 1.5$},
				legend style={draw=none, at={(0,-0.2)},anchor=north west}
				]
				\addplot[color=blue, dotted, very thick] table [x=iter,y=Obj, col sep=comma]{matlab/dataResults/optimalBridge_new/ell1p5_Initial/metricType_2_retractionType_1_DATA.txt};
				\addplot[color=blue, dotted, very thick] table [x expr=\thisrow{iter}+1000,y=Obj, col sep=comma]{matlab/dataResults/optimalBridge_new/ell1p5_WarmStart_variant2/metricType_2_retractionType_1_DATA.txt};
				\addplot[color=magenta,  very thick] table [x=iter,y=Obj, col sep=comma]{matlab/dataResults/optimalBridge_new/ell1p5_Initial/metricType_3_retractionType_1_DATA.txt};
				\addplot[color=magenta,  very thick] table [x expr=\thisrow{iter}+1000,y=Obj, col sep=comma]{matlab/dataResults/optimalBridge_new/ell1p5_WarmStart_variant3/metricType_3_retractionType_1_DATA.txt};
			\end{axis}
		\end{tikzpicture}
		\\
		\pgfplotsset{compat=1.14, width=0.32\textwidth}
		\begin{tikzpicture}
			\begin{axis}[
				title = {Mesh quality with $\ell = 9.9 \cdot 10 ^{-2}$},
				legend style={draw=none, at={(0,-0.2)},anchor=north west}
				]
				\addplot[color=blue, dotted, very thick] table [x=iter,y=mshQua, col sep=comma]{matlab/dataResults/optimalBridge_new/ell99Em3_Initial/metricType_2_retractionType_1_DATA.txt};
				\addplot[color=blue, dotted, very thick] table [x expr=\thisrow{iter}+1000,y=mshQua, col sep=comma]{matlab/dataResults/optimalBridge_new/ell99Em3_WarmStart_variant2/metricType_2_retractionType_1_DATA.txt};
				\addplot[color=magenta,  very thick] table [x=iter,y=mshQua, col sep=comma]{matlab/dataResults/optimalBridge_new/ell99Em3_Initial/metricType_3_retractionType_1_DATA.txt};
				\addplot[color=magenta,  very thick] table [x expr=\thisrow{iter}+1000,y=mshQua, col sep=comma]{matlab/dataResults/optimalBridge_new/ell99Em3_WarmStart_variant3/metricType_3_retractionType_1_DATA.txt};
			\end{axis}
		\end{tikzpicture}
		\hfill
		\pgfplotsset{compat=1.14, width=0.32\textwidth}
		\begin{tikzpicture}
			\begin{axis}[
				title = {Mesh quality with $\ell = 5 \cdot 10 ^{-1}$},
				legend style={draw=none, at={(0,-0.2)},anchor=north west}
				]
				\addplot[color=blue, dotted, very thick] table [x=iter,y=mshQua, col sep=comma]{matlab/dataResults/optimalBridge_new/ell5Em1_Initial/metricType_2_retractionType_1_DATA.txt};
				\addplot[color=blue, dotted, very thick] table [x expr=\thisrow{iter}+700,y=mshQua, col sep=comma]{matlab/dataResults/optimalBridge_new/ell5Em1_WarmStart_variant2/metricType_2_retractionType_1_DATA.txt};
				\addplot[color=magenta, very thick] table [x=iter,y=mshQua, col sep=comma]{matlab/dataResults/optimalBridge_new/ell5Em1_Initial/metricType_3_retractionType_1_DATA.txt};
				\addplot[color=magenta, very thick] table [x expr=\thisrow{iter}+700,y=mshQua, col sep=comma]{matlab/dataResults/optimalBridge_new/ell5Em1_WarmStart_variant3/metricType_3_retractionType_1_DATA.txt};
			\end{axis}
		\end{tikzpicture}
		\hfill
		\pgfplotsset{compat=1.14, width=0.32\textwidth}
		\begin{tikzpicture}
			\begin{axis}[
				title = {Mesh quality with $\ell = 1.5$},
				legend style={draw=none, at={(0,-0.2)},anchor=north west}
				]
				\addplot[color=blue, dotted, very thick] table [x=iter,y=mshQua, col sep=comma]{matlab/dataResults/optimalBridge_new/ell1p5_Initial/metricType_2_retractionType_1_DATA.txt};
				\addplot[color=blue, dotted, very thick] table [x expr=\thisrow{iter}+1000,y=mshQua, col sep=comma]{matlab/dataResults/optimalBridge_new/ell1p5_WarmStart_variant2/metricType_2_retractionType_1_DATA.txt};
				\addplot[color=magenta,  very thick] table [x=iter,y=mshQua, col sep=comma]{matlab/dataResults/optimalBridge_new/ell1p5_Initial/metricType_3_retractionType_1_DATA.txt};
				\addplot[color=magenta,  very thick] table [x expr=\thisrow{iter}+1000,y=mshQua, col sep=comma]{matlab/dataResults/optimalBridge_new/ell1p5_WarmStart_variant3/metricType_3_retractionType_1_DATA.txt};
			\end{axis}
		\end{tikzpicture}
		\caption{Objective and mesh quality for the penalized optimal bridge problem described in \cref{subsection:experimentOptimalBridge} for the variants \ElasEuc (dotted blue) and \CompEuc (solid magenta).}
		\label{fig:Objective_and_MeshQuality_OptimalBridge}
	\end{center}
\end{figure}
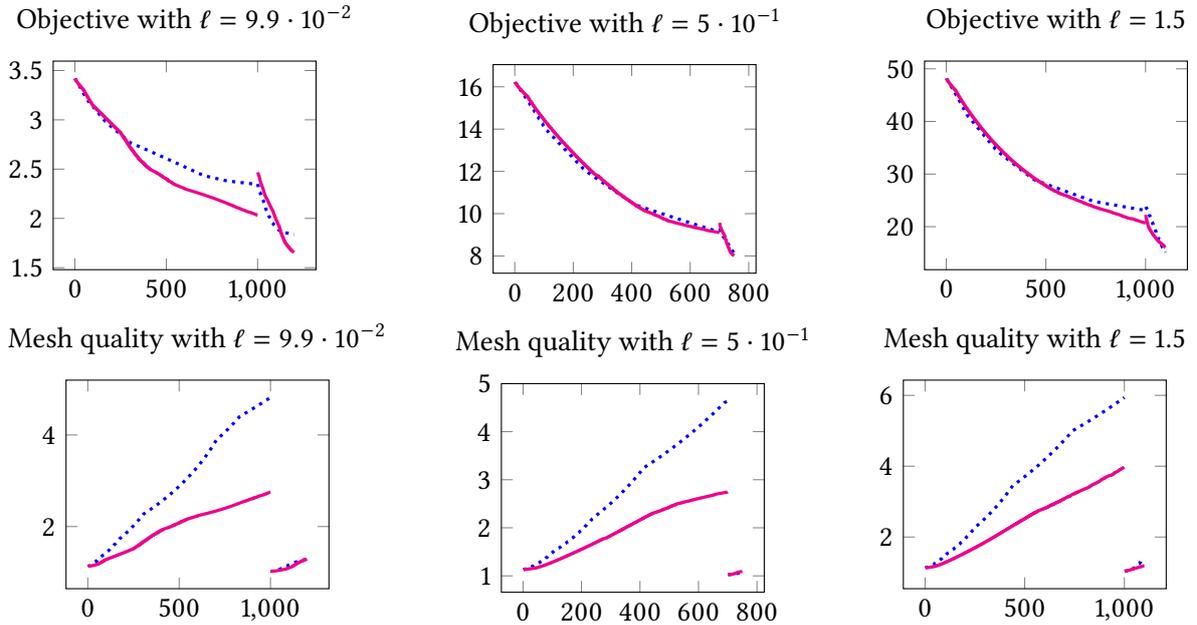

A quantitative comparison of the decay of the objective function and the evolution of the mesh quality for the variants \ElasEuc and \CompEuc are shown in \cref{fig:Objective_and_MeshQuality_OptimalBridge}.
We emphasize the discontinuity of the objective function produced during the remeshing.

\begin{figure}
	\begin{center}
	\begin{subfigure}{0.32\textwidth}
		\centering
		\includegraphics[width = 0.95\textwidth]{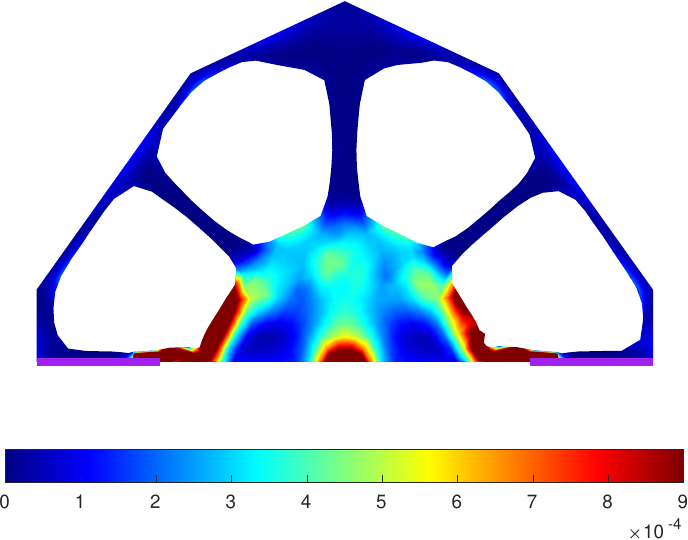}
		\caption{$\ell = 9.9\cdot 10^{-2}$}
	\end{subfigure}
	\hfill
	\begin{subfigure}{0.32\textwidth}
		\centering
		\includegraphics[width = 0.95\textwidth]{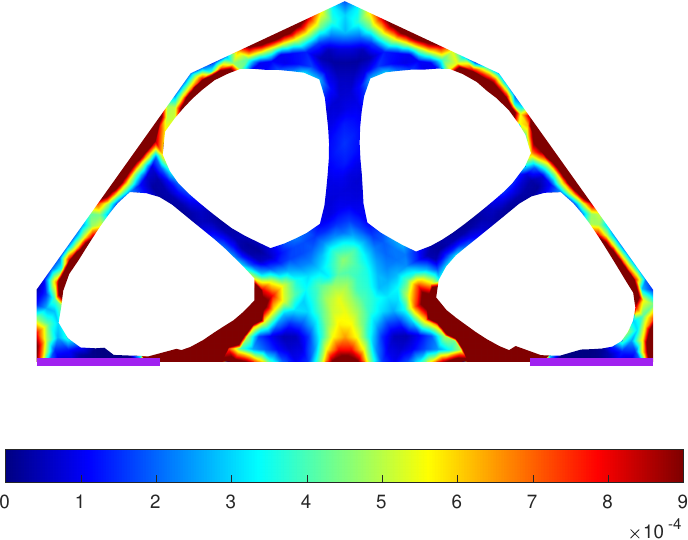}
		\caption{$\ell = 5\cdot 10^{-1}$}
	\end{subfigure}
	\hfill
	\begin{subfigure}{0.32\textwidth}
		\centering
		\includegraphics[width = 0.95\textwidth]{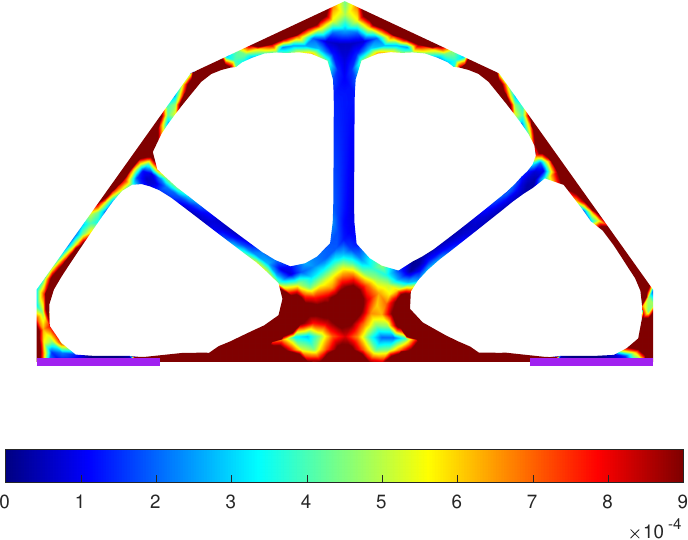}
		\caption{$\ell = 1.5$}
	\end{subfigure}
	\\
	\begin{subfigure}{0.32\textwidth}
		\centering
		\includegraphics[width = 0.95\textwidth]{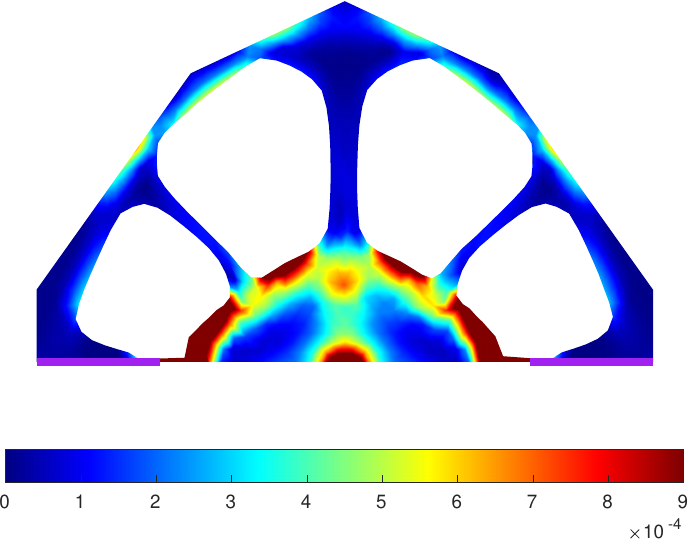}
		\caption{$\ell = 9.9\cdot 10^{-2}$}
	\end{subfigure}
	\hfill
	\begin{subfigure}{0.32\textwidth}
		\centering
		\includegraphics[width = 0.95\textwidth]{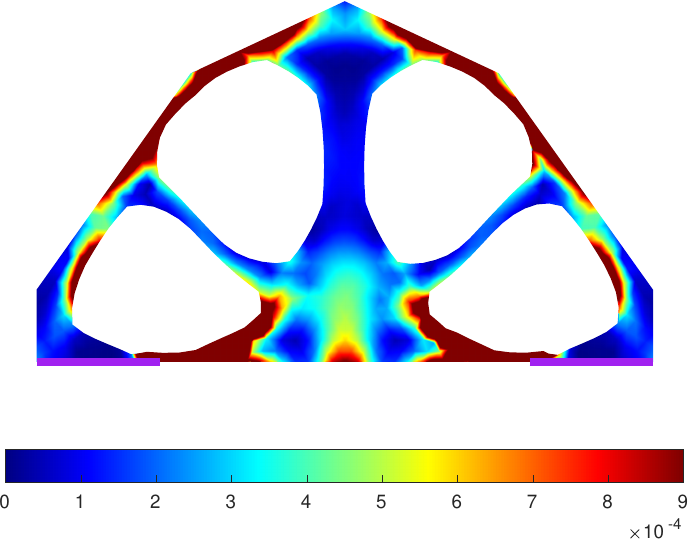}
		\caption{$\ell = 5\cdot 10^{-1}$}
	\end{subfigure}
	\hfill
	\begin{subfigure}{0.32\textwidth}
		\centering
		\includegraphics[width = 0.95\textwidth]{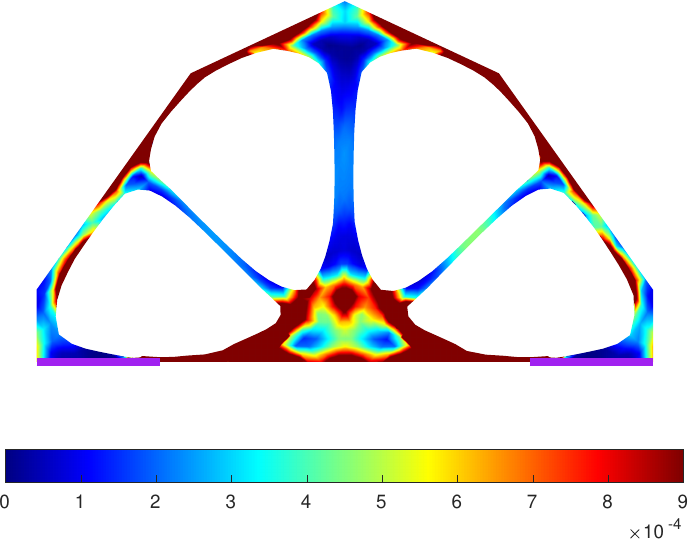}
		\caption{$\ell = 1.5$}
	\end{subfigure}
\end{center}
	\caption{Associated cell-wise elastic energy for the variant \CompEuc (upper row) and \ElasEuc (lower row) from the optimal bridge experiment described in~\cref{subsection:experimentOptimalBridge}}
	\label{fig:ElasticDef}
\end{figure}

Finally, we display in \cref{fig:ElasticDef} the local contributions of each mesh cell~$T$ to the elastic energy, \ie,
\begin{equation*}
	\muState \int_T \bvarepsilon(\by) \dprod \bvarepsilon(\by) \d x
	+
	\frac{\lambdaState}{2} \int_T \trace(\bvarepsilon(\by)) \, \trace(\bvarepsilon(\by)) \d x
\end{equation*}

\begin{table}
\begin{center}
\sisetup{round-mode = places, round-precision = 3}
		\begin{filecontents*}{matlab/dataResults/optimalBridge_new/summaryOptimalBridge.txt}
ell,Nq,Nt,Iter,timeVar2,timeVar3
0.099,467,816,1000,87.82167,66.482806
0.099,702,1137,200,34.201715,24.148767
0.099,714,1154,200,36.244948,25.341103
0.5,467,816,700,58.785719,44.616238
0.5,642,1006,50,7.060599,5.045124
0.5,652,1027,50,7.042115,5.292769
1.5,467,816,1000,83.883444,64.180465
1.5,594,898,100,11.788991,8.778348
1.5,581,856,100,11.402449,8.604922
\end{filecontents*}
\pgfplotstabletypeset[%
		col sep = comma,
		every head row/.style={before row=\toprule,after row=\midrule},
		every last row/.style={after row=\bottomrule},
		every nth row={3}{before row=\midrule},
		columns/ell/.style={
		assign cell content/.code={%
			\ifnum\pgfplotstablerow=0
				\pgfkeyssetvalue{/pgfplots/table/@cell content}%
				{\multirow{3}{*}{\pgfmathprintnumber[sci, sci zerofill,retain unit mantissa=false]{##1}}}%
			\else
				\ifnum\pgfplotstablerow=3
					\pgfkeyssetvalue{/pgfplots/table/@cell content}%
					{\multirow{3}{*}{\pgfmathprintnumber[sci, sci zerofill,retain unit mantissa=false]{##1}}}
				\else
					\ifnum\pgfplotstablerow=6
						\pgfkeyssetvalue{/pgfplots/table/@cell content}%
						{\multirow{3}{*}{\pgfmathprintnumber[sci, sci zerofill,retain unit mantissa=false]{##1}}}%
					\else
						\pgfkeyssetvalue{/pgfplots/table/@cell content}{}%
					\fi
				\fi
			\fi
		},
		column name=$\ell$},
		columns/Nq/.style={
		int detect,
		column name={$N_V$}},
		columns/Nt/.style={
		int detect,
		column name=$N_T$},
		columns/Iter/.style={
		int detect,
		column name= iter $(n)$},
		columns/timeVar2/.style={std, precision=2,column name={\ElasEuc},
		assign cell content/.code={%
		\pgfkeyssetvalue{/pgfplots/table/@cell content}%
		{\SI{##1}{\second}}%
		}},
		columns/timeVar3/.style={std, precision=2,column name={\CompEuc},
		assign cell content/.code={%
		\pgfkeyssetvalue{/pgfplots/table/@cell content}%
		{\SI{##1}{\second}}%
		}}
		]{matlab/dataResults/optimalBridge_new/summaryOptimalBridge.txt}
\end{center}
\caption{Execution times for the experiment described in \cref{subsection:experimentOptimalBridge}.}
\label{table:timesExecutionOptimalBridge}
\end{table}

\section{Conclusions and Outlook}
\label{section:conclusions_outlook}

We studied a discretized, PDE-constrained shape optimization problem, in which the shape is represented by a triangular mesh and the vertex positions serve as optimization variables.
The PDE under consideration is discretized using a standard finite element approach.
This is a common approach in computational shape optimization.

We clarified that the set of admissible vertex positions, which maintain the connectivity and orientation of the initial mesh, forms an open, connected submanifold of the vector space of all vertex positions.
Unfortunately, the minimization of a typical shape optimization objective often does not have a solution even if the objective is bounded below.
The reason is that the infimum is not attained by an admissible vertex configuration.
This results in optimizers tending to points on the boundary of the manifold, which correspond to degenerate meshes, thereby exploiting, \eg, the quadrature error, in order to reduce the objective below the optimal value of the continuous problem.
To the best of our knowledge, this fact has not been reported explicitly in the literature.

The situation is somewhat reminiscent of a class of ill-posed inverse problems, although the ill-posedness appears only in the discretized problem.
It can be dealt with by, \eg, explicit regularization.
This approach corresponds to the addition of an appropriate penalty function in shape optimization.
For this purpose, we proposed a novel penalty function, whose properness ensures the existence of optimal shapes within the manifold of admissible vertex positions.
If this is not desired, then early stopping of the optimizer can still provide a reasonable approximation of the continuous solution.

We also proposed a novel Riemannian metric for discretized shape optimization problems, which governs the formation of the gradient direction in a steepest descent method.
This metric is derived from the same penalty function which was used to ensure the existence of minimizers.
Its advantage over other metrics in use, such as the Lamé system, is that it is represented by a rank-$1$ perturbation of the identity matrix and thus the conversion of the derivative to the gradient can be achieved very efficiently and in a matrix-free way by performing two conjugate gradient iterations.
Our numerical experiments show that the new metric admits a gradient algorithm which compares favorably, both for penalized and unpenalized problems, with a gradient method based on the elasticity metric.
Also, it is sufficient to employ a cheap (Euclidean) update of the vertex positions.

In follow-up work, it would be interesting to study the proposed metric as the base metric in a quasi-Newton scheme to accelerate convergence.
Another open question is whether the presence of the proposed penalty terms is compatible with the continuous limit problem, or else whether minimizers to the penalized discrete problems fail to converge to a minimizer of the continuous problem.

\appendix

\section{Proof of \texorpdfstring{\cref{proposition:edgelengths_heights_bounded_sublevelset}}{Proposition~\ref{proposition:edgelengths_heights_bounded_sublevelset}}}
\label{appendix:proofs_varphi}

This appendix is devoted to the study of some of the properties of the penalty function~$\newpenalty$ and, in particular, the proof of \cref{proposition:edgelengths_heights_bounded_sublevelset}.
In particular we prove that in any sublevel set of $\newpenalty$, all the edge lengths, and heights of any mesh are bounded.
Our proof leverages the $2$-path connectedness of connectivity complexes~$\Delta$; see \cref{definition:connectivity_complex}.

\subsubsection*{Proof of~\cref{proposition:edgelengths_heights_bounded_sublevelset}}
We recall from \eqref{eq:penalizedProblem} the definition of the penalty function~$\newpenalty$,
\begin{equation*}
	\begin{aligned}
			\newpenalty(Q;\Qref)
			&
			\coloneqq
			\sum_{k=1}^{N_T} \frac{1}{N_T} \frac{\alpha_1}{\psi_Q(i_0^k,i_1^k,i_2^k)}
			+
			\frac{\alpha_2}{\sum_{k=1}^{N_T}\area[auto]{Q}[i_0^k,i_1^k,i_2^k] }
			\\
			&
			\quad
			+
			\sum_{[j_0,j_1] \in E_\partial}
			\sum_{\substack{i_0 \in V_\partial \\ i_0 \neq j_0,j_1}} \frac{1}{\#E_\partial \#V_\partial} \frac{\alpha_3}{\regularizedDistance{Q}[i_0][[j_0,j_1]]}
			+
			\frac{\alpha_4}{2} \norm{Q-\Qref}_F^2
			.
		\end{aligned}
\end{equation*}
Since we keep $Q$ fixed throughout the proof, we simplify the notation and drop the dependence on $Q$.
Thus we write $\edgelengthTriangle{k}{\ell}$ in place of $\edgelength{Q}{\ell}(i_0^k,i_1^k,i_2^k)$ for the edge lengths, $\areaTriangle{k}$ in place of $\area[auto]{Q}[i_0^k,i_1^k,i_2^k]$, and $\psi_k$ in place of $\psi_Q(i_0^k,i_1^k,i_2^k)$.

We recall from \eqref{eq:sublevel_set} that $\cN_b$ denotes a non-empty sublevel set of $\newpenalty$.
Let us consider $Q \in \cN_b$ arbitrary but fixed.
The proof of the proposition is broken down into several steps.

\begin{enumeratearabic}
	\item
		\label[step]{step:bound_one_edge}
		We find upper and lower bounds for the length of the longest edge of one particular triangle, denoted as $\edgelengthTriangle{\overline{k}}{\overline{\ell}}$.

	\item
		\label[step]{step:bound_remaining_edges_bar_k}
		Using the bounds from \cref{step:bound_one_edge} we find upper and lower bounds for the remaining edges of the $\overline{k}$-th triangle.

	\item
		\label[step]{step:bound_heights}
		We compute lower bounds for the heights~$\heightTriangle{\overline{k}}{\ell}$ of the $\overline{k}$-th triangle using the results of \cref{step:bound_remaining_edges_bar_k}.

	\item
		\label[step]{step:bounds_for_remaining_triangles}
		We consider an arbitrary triangle $k$ different from $\overline{k}$.
		Based on the $2$-path connectedness of $\Delta$ we use the bounds from \cref{step:bound_one_edge,step:bound_remaining_edges_bar_k} to find a lower bound for all edge lengths of the $k$-th triangle.
\end{enumeratearabic}
We point out that all bounds are going to be independent of~$Q$ but they only depend on~$b$.

Since $Q \in \cN_b$ holds, we immediately obtain
\begin{equation}
	\sum_{k=1}^{N_T} \areaTriangle{k}
	\ge
	\frac{\alpha_2}{b}
	.
\end{equation}
Since the areas are all positive, there exists at least one triangle $\overline{k}$ such that
\begin{equation*}
	\areaTriangle{\overline{k}}
	\ge
	\frac{1}{N_T} \frac{\alpha_2}{b}
	.
\end{equation*}
We now use the so-called isoperimetric inequality for triangles, see \cite[Theorem~25, p.42]{AgricolaFriedrich:2008:1}, which states
\begin{equation}
	\label{eq:isoperimetric_inequality}
	\areaTriangle{\overline{k}}
	\le
	\frac{\paren[Big](){%
		\edgelengthTriangle{\overline{k}}{0} + \edgelengthTriangle{\overline{k}}{1} + \edgelengthTriangle{\overline{k}}{2}
		}^2
	}%
	{12 \sqrt{3}}
	.
\end{equation}
Denoting $\edgelengthTriangle{\overline{k}}{\overline{\ell}} \coloneqq \max_{\ell=0,1,2} \edgelengthTriangle{\overline{k}}{\ell}$, we obtain
\begin{equation*}
	\edgelengthTriangle{\overline{k}}{\overline{\ell}}
	\ge
	\frac{2}{3^{\nicefrac{1}{4}}}
	\paren[auto](){\frac{1}{N_T} \frac{\alpha_2}{b}}^{\nicefrac{1}{2}}
	>
	0
	.
\end{equation*}
Notice moreover, that $Q \in \cN_b$ implies $\norm{Q-\Qref}_F^2 \le 2b / \alpha_4$, which in turn implies $\norm{Q}_F \le \sqrt{2b / \alpha_4} + \norm{\Qref}_F$.
We denote by $\overline{i}_0$ and $\overline{i}_1$ the vertices of triangle~$\overline{k}$ which form the edge whose edge length is $\edgelengthTriangle{\overline{k}}{\overline{\ell}}$.
Then we can estimate $\edgelengthTriangle{\overline{k}}{\overline{\ell}} = \norm{q_{\overline{i}_0} - q_{\overline{i}_1}} \le \norm{q_{\overline{i}_0}} + \norm{q_{\overline{i}_1}}\le \sqrt{2} \, \norm{Q}_F$.
Thus, $\edgelengthTriangle{\overline{k}}{\overline{\ell}} \le 2 \, \sqrt{b / \alpha_4} + \sqrt{2} \, \norm{\Qref}_F$.
Altogether we found
\begin{equation}
	\label{eq:bounds_bar_k}
	\frac{2}{3^{\nicefrac{1}{4}}}
	\paren[auto](){\frac{1}{N_T} \frac{\alpha_2}{b}}^{\nicefrac{1}{2}}
	\le
	\edgelengthTriangle{\overline{k}}{\overline{\ell}}
	\le
	2 \, \sqrt{b / \alpha_4} + \sqrt{2} \, \norm{\Qref}_F
\end{equation}
This concludes \cref{step:bound_one_edge}.

Now, we proceed to find upper and lower bounds for $\edgelengthTriangle{\overline{k}}{\overline{\ell} \oplus 1}$ and $\edgelengthTriangle{\overline{k}}{\overline{\ell} \oplus 2}$ for $j = 1,2$.
Recall that $\edgelengthTriangle{\overline{k}}{\overline{\ell}}$ denotes the length of the longest edge.
Thus, it holds $\edgelengthTriangle{\overline{k}}{\overline{\ell} \oplus 1},\edgelengthTriangle{\overline{k}}{\overline{\ell} \oplus 2} \le  2 \, \sqrt{b / \alpha_4} + \sqrt{2} \, \norm{\Qref}_F $.
On the other hand, from $Q \in \cN_b$ and the definition of $\psi_k$ given in \eqref{eq:Shewchuk64_qualitymeasure} it follows that
\begin{equation*}
	b
	\ge
	\newpenalty(Q;\Qref)
	\ge
	\frac{\alpha_1}{N_T}
	\frac{
		\paren[Big](){\edgelengthTriangle{\overline{k}}{\overline{\ell}}}^2
		+
		\paren[Big](){\edgelengthTriangle{\overline{k}}{\overline{\ell} \oplus 1}}^2
		+
		\paren[Big](){\edgelengthTriangle{\overline{k}}{\overline{\ell} \oplus 2}}^2
	}{4 \, \sqrt{3} \, \areaTriangle{\overline{k}}}
	.
\end{equation*}
For the triangle area we have $\areaTriangle{\overline{k}} = \edgelengthTriangle{\overline{k}}{\overline{\ell}} \, \heightTriangle{\overline{k}}{\overline{\ell}} / 2$.
Moreover, it is easy to see that $\heightTriangle{k}{\ell}\le \edgelengthTriangle{k}{\overline{\ell} \oplus 1}$ and $\heightTriangle{k}{\ell}\le \edgelengthTriangle{k}{\ell \oplus 2}$ holds for all $k=1,\ldots,N_T$ and all $\ell=0,1,2$.
We will focus here on the bounds for the edge length $\edgelengthTriangle{\overline{k}}{\overline{\ell} \oplus 1}$.
The bounds for $\edgelengthTriangle{\overline{k}}{\overline{\ell} \oplus 2}$ can be obtained using the same arguments.
From the previous estimates we obtain $\areaTriangle{\overline{k}} \le \edgelengthTriangle{\overline{k}}{\overline{\ell}} \edgelengthTriangle{\overline{k}}{\overline{\ell} \oplus 1} / 2$.
Using the lower bound for $\edgelengthTriangle{\overline{k}}{\overline{\ell}}$ in \eqref{eq:bounds_bar_k}, we obtain the following estimate
\begin{equation*}
	b
	\ge
	\newpenalty(Q;\Qref)
	\ge
	\frac{\alpha_1}{N_T}
	\frac{(2 / 3^{\nicefrac{1}{4}}) \paren[auto](){\alpha_2 / (N_T \, b)}^{\nicefrac{1}{2}}}{2 \, \sqrt{3} \, \edgelength{\overline{k}}{\overline{\ell} \oplus 1}}
	.
\end{equation*}
This implies $\edgelengthTriangle{\overline{k}}{\overline{\ell} \oplus 1} \ge \alpha_1 \, \alpha_2^{\nicefrac{1}{2}} / \paren[auto](){{3^{\nicefrac{3}{4}}(N_T \, b)^{\nicefrac{3}{2}}}}$.
Summarizing, the edge lengths of the $\overline{k}$-th triangle satisfy
\begin{equation}
	\label{eq:bounds_edgelength_triangle_bark}
		\min\paren[auto]\{\}{
			\frac{2 \alpha_2^{\nicefrac{1}{2}}}{3^{\nicefrac{1}{4}}(N_T \, b)^{\nicefrac{1}{2}}}
			,
			\frac{\alpha_1 \, \alpha_2^{\nicefrac{1}{2}}}{3^{\nicefrac{3}{4}}(N_T \, b)^{\nicefrac{3}{2}}}
		}
		\le
		\edgelengthTriangle{\overline{k}}{\ell}
		\le
		2\sqrt{b / \alpha_4} + \sqrt{2} \norm{\Qref}_F
\end{equation}
for all $\ell=0,1,2$.
Moreover, to simplify notation let us denote as
\begin{align}
	\label{eq:c0_definition}
	c_0
	&
	\coloneqq
	\min
	\paren[auto]\{\}{
		\frac{2 \alpha_2^{\nicefrac{1}{2}}}{3^{\nicefrac{1}{4}}(N_T \, b)^{\nicefrac{1}{2}}}
		,
		\frac{\alpha_1 \, \alpha_2^{\nicefrac{1}{2}}}{3^{\nicefrac{3}{4}}(N_Tb)^{\nicefrac{3}{2}}}
	}
	,
	\\
	\label{eq:C0_definition}
	C_0
	&
	\coloneqq
	2\sqrt{b / \alpha_4}
	+
	\sqrt{2} \, \norm{\Qref}_F
	.
\end{align}
Thus we have concluded \cref{step:bound_remaining_edges_bar_k}.

The bounds from \cref{step:bound_heights} are immediately obtained from $\areaTriangle{\overline{k}} = \edgelengthTriangle{\overline{k}}{\overline{\ell}} \, \heightTriangle{\overline{k}}{\overline{\ell}} / 2 = \edgelengthTriangle{\overline{k}}{\overline{\ell} \oplus 1} \, \heightTriangle{\overline{k}}{\overline{\ell} \oplus 1} / 2$ and using the bounds from \cref{step:bound_one_edge,step:bound_remaining_edges_bar_k}.
Thus, we conclude
\begin{equation*}
	\frac{1}{\heightTriangle{\overline{k}}{\ell}}
	\le
	\frac{2\sqrt{3} N_T \, b}{\alpha_1 \, c_0}
\end{equation*}
for all $\ell = 0,1,2$.

Finally, we focus on \cref{step:bounds_for_remaining_triangles}.
Having found the constants for the $\overline{k}$-th triangle, we will use it as a pivot to compute the constants for the remaining triangles, based on the $2$-path connectedness of $\Delta$.
To this end, we consider an arbitrary triangle~$k$ different from $\overline{k}$.
From all possible paths joining these triangles, guaranteed to exist by the $2$-path connectedness, we choose a shortest one.
Since $\Delta$ is a finite collection of simplices, there is an upper bound~$L \le N_T$ on the lengths of the shortest paths between any two triangles.

\begin{figure}[htp]
	\centering
	\includegraphics[width=0.5\textwidth]{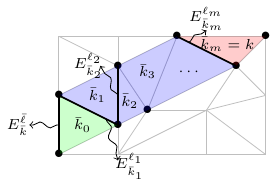}
	\caption{Illustration of the path of triangles used in \cref{step:bounds_for_remaining_triangles} of the proof of \cref{proposition:edgelengths_heights_bounded_sublevelset} in \cref{appendix:proofs_varphi}.}
	\label{fig:bounded_edgelengths_heights_complex_mesh}
\end{figure}

Suppose that the path joining the $\overline{k}$-th and the $k$-th triangles has $m + 1 \le L$ elements; see \cref{fig:bounded_edgelengths_heights_complex_mesh} for an illustration.
We denote the triangles involved by $\overline{k} = \overline{k}_0$, then $\overline{k}_1$ etc.\ up to $\overline{k}_m = k$.
For the initial triangle we know the bounds given in~\eqref{eq:bounds_edgelength_triangle_bark} hold, \ie,
\begin{equation}
	\label{eq:bounds_edgelength_T0}
	c_0
	\le
	\edgelengthTriangle{\overline{k}_0}{\ell}
	\le
	C_0
\end{equation}
for all $\ell=0,1,2$.
Triangles $\overline{k}_0$ and $\overline{k}_1$ share an edge, and denote its length as seen from the first triangle as $\edgelengthTriangle{\overline{k}_1}{\ell_1}$, for which \eqref{eq:bounds_edgelength_T0} also holds.
Using the same techniques as before one can prove
\begin{equation}
	\label{eq:bounds_edgelength_T1}
	\frac{\alpha_1 \, c_0}{2 \sqrt{3} N_T \, b}
	\le
	\edgelengthTriangle{\overline{k}_1}{\ell_1 \oplus 1}
	,
	\;
	\edgelengthTriangle{\overline{k}_1}{\ell_1 \oplus 2}
	\le
	\frac{2\sqrt{3}N_T \, b C_0}{\alpha_1}
\end{equation}
We denote by $c_1 \coloneqq \min\{c_0,(\nicefrac{\alpha_1}{2\sqrt{3} N_T \, b}) \, c_0\}$ and $C_1 \coloneqq \max \{C_0, (2 \, \sqrt{3} N_T \, b/\alpha_1) \, C_0\}$.
	Then, it holds $c_1 \le \edgelengthTriangle{\overline{k_1}}{\ell} \le C_1$ for all $\ell = 0,1,2$.
In the same manner we can bound the heights as follows,
\begin{equation*}
	\frac{1}{\heightTriangle{\overline{k}_1}{\ell}}
	\le
	\frac{2\sqrt{3} N_T \, b}{\alpha_1 \, c_1}
	.
\end{equation*}

By repeating this process until we reach last element of the path, \ie, the triangle~$\overline{k}_m = k$.
We obtain the following bounds
\begin{equation*}
	c_m
	\le
	\edgelengthTriangle{\overline{k}_m}{\ell}
	\le
	C_m
	\quad
	\text{and}
	\quad
	\frac{1}{\heightTriangle{\overline{k}_m}{\ell}}
	\le
	\frac{2\sqrt{3}N_T \, b}{\alpha_1 c_m}
\end{equation*}
where
\begin{align*}
	c_m
	&
	=
	\min
	\paren[auto]\{\}{
		c_{m-1}
		,
		\paren[auto](){
			\frac{\alpha_1}{2\sqrt{3}N_T \, b}
		}
		c_{m-1}
	}
	,
	\\
	C_m
	&
	=
	\max
	\paren[auto]\{\}{
		C_{m-1}
		,
		\paren[auto](){
			\frac{2\sqrt{3}N_T \, b }{\alpha_1}
		}
		C_{m-1}
	}
	.
\end{align*}
Notice that by the way the constants $c_m$ and $C_m$ are defined, it holds $c_m \le c_{m-1} \le \cdots \le c_0$ and $C_m \geq C_{m-1}\geq \cdots \geq C_0$.
Recalling we have denoted by $L\le N_T$ the length of the longest paths we can conclude that for all $k = 1,\ldots, N_T$ and all $\ell = 0,1,2$
\begin{equation*}
	c_L
	\le
	\edgelengthTriangle{k}{\ell}
	\le
	C_L
	\quad
	\text{and}
	\quad
	\frac{1}{\height{k}{\ell}}
	\le
	\frac{2\sqrt{3} N_T \, b}{\alpha_1 c_L}
\end{equation*}
holds, where the constants $c_L$ and $C_L$ do not depend on $Q$, the chosen pivot triangle $\overline{k}$ nor on the edge $\overline{\ell}$.
They do depend, however, on $\alpha_1$, $\alpha_2$ and $\alpha_4$ as well as the connectivity of the mesh encoded in~$\Delta.$

\section*{Acknowledgments}

The second author was partially funded by the Deutsche Forschungsgemeinschaft (DFG, German Research Foundation) under Germany's Excellence Strategy EXC~2044--390685587, Mathematics Münster: Dynamics--Geometry--Structure.
We wish to thank two anonymous reviewers, whose comments helped to improve the manuscript.
In particular, we are indebted to one reviewer who suggested considering the compliance minimization problem.